\documentclass{report}
\input{diagrams}
\usepackage{epsfig}

\newtheorem{thm}{Theorem}[section]
\newtheorem{defn}[thm]{Definition}
\newtheorem{propn}[thm]{Proposition}
\newtheorem{lotsofegs}[thm]{Examples}
\newenvironment{eg}[1]
	{\begin{lotsofegs}\label{#1}\end{lotsofegs}\begin{enumerate}}
	{\end{enumerate}}

\newcommand{\mcm}[3]{\newcommand{#1}[#2]{{\ensuremath{#3}}}}
\usepackage{latexsym}
\usepackage{amstex}
\mcm{\emptybk}{0}{\:\:}
\mcm{\blank}{0}{(\emptybk)}
\mcm{\dashbk}{0}{\mbox{---}}
\mcm{\hyph}{0}{\mbox{-}}
\mcm{\diagspace}{0}{\mbox{\hspace{2em}}}
\newcommand{\done}{\hfill\ensuremath{\Box}}
\newcommand{\ucontents}[2]{\addcontentsline{toc}{#1}{\numberline{}{#2}}}
\mcm{\mc}{1}{\mathcal{#1}}
\mcm{\mr}{1}{\mathrm{#1}}
\mcm{\mi}{1}{\mathit{#1}}
\mcm{\mb}{1}{\mathbf{#1}}
\mcm{\cat}{1}{\mc{#1}}
\mcm{\scat}{1}{\Bbb{#1}}
\mcm{\fcat}{1}{\mb{#1}}
\newcommand{\url}[1]{\mbox{\tt #1}}
\mcm{\elt}{0}{\in}
\mcm{\such}{0}{\:|\:}
\mcm{\blob}{0}{\scriptscriptstyle{\bullet}}
\mcm{\of}{0}{\raisebox{0.2mm}{\ensuremath{\scriptstyle\circ}}}
\mcm{\op}{0}{\mr{op}}
\mcm{\dom}{0}{\mr{dom}}
\mcm{\cod}{0}{\mr{cod}}
\mcm{\Hom}{0}{\mr{Hom}}
\mcm{\End}{0}{\mr{End}}
\mcm{\id}{0}{\mi{id}}
\mcm{\Set}{0}{\fcat{Set}}
\mcm{\Cat}{0}{\fcat{Cat}}
\mcm{\Multicat}{0}{\fcat{Multicat}}
\mcm{\Graph}{0}{\fcat{Graph}}
\mcm{\nat}{0}{\mathbb{N}}	
\mcm{\pr}{2}{\tuplebts{#1,#2}}
\mcm{\triple}{3}{\tuplebts{#1,#2,#3}}
\mcm{\range}{2}{#1,\,\ldots\,,#2}
\mcm{\tuplebts}{1}{(#1)}
\mcm{\tuple}{3}{\tuplebts{\range{#1,#2}{#3}}}
\mcm{\bftuple}{2}{\tuplebts{\range{#1}{#2}}}
\mcm{\atuplebts}{1}{\langle #1 \rangle}
\mcm{\abftuple}{2}{\atuplebts{\range{#1}{#2}}}
\mcm{\ftrcat}{2}{[#1,#2]}
\mcm{\homset}{3}{#1(#2,#3)}
\mcm{\multihom}{3}{#1(#2;#3)} 
\mcm{\go}{0}{\rTo}
\mcm{\og}{0}{\lTo}
\mcm{\goby}{1}{\rTo^{#1}}
\mcm{\ogby}{1}{\lTo^{#1}}
\mcm{\goiso}{0}{\goby{\diso}}
\mcm{\goesto}{0}{\,\longmapsto\,}
\mcm{\spn}{3}{#2 \og #1 \go #3}
\mcm{\spaan}{5}{#2 \ogby{#4} #1 \goby{#5} #3}
\mcm{\parpair}{2}{\stackrel{\rTo^{#1}}{\rTo_{#2}}}
\mcm{\vslob}{3}
	{\left.
	\begin{diagram}[height=1.5em]
	#1		\\
	\dTo>{\,#2}	\\
	#3		\\
	\end{diagram}
	\right.}
\newarrow{Get}....>
\newarrow{NT}===={=>}
\newarrow{Equals}=====
\newarrow{Edge}---->		
\newenvironment{slopeydiag}
	{\begin{diagram}[size=2em]}
	{\end{diagram}}

\newenvironment{opetope}
	{\begin{diagram}[size=1em,abut,tight,noPS]}
	{\end{diagram}}
\newcommand{\pullshape}
	{\setlength{\unitlength}{1em}
	\begin{picture}(2,5)(-1,-5)
	\put(0,-5){\line(1,1){1}}
	\put(0,-5){\line(-1,1){1}}
	\end{picture}}
\newcommand{\Spbk}{\overprint{\raisebox{-2.5em}{\pullshape}}}
\mcm{\diso}{0}{\sim}
\newcommand{\cnr}{}	
\newarrow{Monic}{vee}---{vee}
\mcm{\monic}{0}{\rMonic}
\mcm{\Struc}{0}{\fcat{Struc}}
\mcm{\integers}{0}{\mathbb{Z}}
\newarrow{Mod}--+->
\mcm{\iso}{0}{\,\cong\,}
\mcm{\eqv}{0}{\,\simeq\,}
\newenvironment{lopetope}
	{\begin{diagram}[size=1.5em,tight]}
	{\end{diagram}}
\mcm{\sub}{0}{\,\subseteq\,}
\newenvironment{ntdiag}
	{\begin{diagram}[size=1.5em]}	
	{\end{diagram}}
\mcm{\ob}{1}{\mr{ob}\,#1}
\newcommand{\pf}{\noindent\textbf{Proof}}
\mcm{\node}{0}{\bullet}
\newenvironment{tree}
	{\begin{diagram}[height=1em,width=.75em,abut,noPS,tight]}	
	{\end{diagram}}
\newcommand{\lt}[1]{\ldLine(#1,2)}
\newcommand{\rt}[1]{\rdLine(#1,2)}
\newcommand{\dn}{\dLine}
\mcm{\nl}{1}{\stackrel{\textstyle #1}{\node}}
\mcm{\bktdvslob}{3}
	{\left(
	\begin{diagram}[height=1.5em]
	#1		\\
	\dTo>{\,#2}	\\
	#3		\\
	\end{diagram}
	\right)}

\mcm{\ust}{0}{{}^{*}}
\mcm{\stbk}{0}{\blank^{*}}
\mcm{\comp}{0}{\mi{comp}}
\mcm{\ids}{0}{\mi{ids}}
\mcm{\Cpn}{1}{\pr{\Set/S_{#1}}{T_{#1}}}
\mcm{\PD}{1}{\fcat{PD}_{#1}}
\mcm{\TR}{1}{\fcat{TR}(#1)}
\mcm{\Eee}{0}{\cat{E}}
\mcm{\Eeep}{0}{\cat{E'}}
\mcm{\Cartpr}{0}{\pr{\Eee}{T}}
\mcm{\Cartprp}{0}{\pr{\Eeep}{T'}}
\mcm{\gph}{2}{\spn{#1}{T #2}{#2}}
\mcm{\graph}{4}{\spaan{#1}{T #2}{#2}{#3}{#4}}
\mcm{\Gph}{0}{\fcat{Graph}}
\mcm{\fc}{0}{\fcat{fc}}
\mcm{\fm}{0}{\fcat{fm}}
\mcm{\Zeropr}{0}{\pr{\Set}{\id}}
\mcm{\Onepr}{0}{\pr{\Gph}{\fc}}
\mcm{\Bim}{1}{\fcat{Bim}(#1)}
\mcm{\Bee}{0}{\cat{B}}
\mcm{\Ab}{0}{\fcat{Ab}}
\mcm{\Span}{0}{\fcat{Span}}
\newcommand{\bref}[1]{(\ref{#1})}
\mcm{\ehom}{3}{#1[#2,#3]}
\mcm{\eend}{2}{#1[#2]}
\mcm{\etrip}{3}{[ #1,#2,#3 ]}
\mcm{\esing}{1}{[ #1 ]}
\mcm{\Prof}{0}{\fcat{Prof}}

\mcm{\commacat}{2}{(#1\!\downarrow\! #2)}
\mcm{\pointyspn}{3}{%
\left(
\begin{slopeydiag}
	&	&#1	&	&	\\
	&\ldTo	&	&\rdTo	&	\\
#2	&	&	&	&#3	\\
\end{slopeydiag}
\right)}
\mcm{\utree}{0}{\bullet}
\mcm{\atsr}{0}{\Box}
\mcm{\Mnd}{0}{\triple{T}{\eta}{\mu}}
\mcm{\relhom}{5}{#1_{#2}(\range{#3}{#4};#5)}
\mcm{\tropset}{0}{\ftrcat{\fcat{TR}^{\op}}{\Set}}
\mcm{\dopset}{0}{\ftrcat{\Delta^{\op}}{\Set}}
\mcm{\odblbkt}{0}{[\![}
\mcm{\cdblbkt}{0}{]\!]}
\newenvironment{minidiagram}
	{\begin{diagram}[size=1.5em]}
	{\end{diagram}}
\newcommand{\epsln}{\varepsilon}
\newcommand{\obt}{[ }
\newcommand{\cbt}{] }
\mcm{\sqbftuple}{2}{\obt\range{#1}{#2}\cbt}
\mcm{\END}{0}{\fcat{End}}
\mcm{\Hpn}{1}{\pr{\Eee_{#1}}{P_{#1}}}
\mcm{\piccy}{1}{\begin{array}{c}\epsfig{file=Q#1.ps}\end{array}}

\newlength{\gwidth}	
\newlength{\gvert}	
\newlength{\gdrop}	
\newlength{\gbaredrop}	
\newlength{\goffset}	
\newlength{\gtemp}	
\newcommand{\present}[1]{%
\makebox[1\gwidth]{%
\rule[-1\gdrop]{0ex}{1\gvert}%
\raisebox{-1\gbaredrop}{#1}}}

\newcommand{\ginitdims}[2]{%
\setlength{\unitlength}{1em}%
\setlength{\goffset}{.25\unitlength}%
\setlength{\gwidth}{#1\unitlength}%
\setlength{\gvert}{#2\unitlength}%
\setlength{\gdrop}{.5\gvert}%
\addtolength{\gdrop}{-1\goffset}%
\setlength{\gbaredrop}{1\gdrop}%
\addtolength{\gvert}{.6\unitlength}%
\addtolength{\gdrop}{.3\unitlength}}	
\newcommand{\cinitdims}[2]{%
\setlength{\unitlength}{1em}%
\setlength{\goffset}{.35\unitlength}%
\setlength{\gwidth}{#1\unitlength}%
\setlength{\gvert}{#2\unitlength}%
\setlength{\gdrop}{.5\gvert}%
\addtolength{\gdrop}{-1\goffset}%
\setlength{\gbaredrop}{1\gdrop}%
\addtolength{\gvert}{.6\unitlength}%
\addtolength{\gdrop}{.3\unitlength}}	
\newcommand{\abovepic}[1]{%
\settoheight{\gtemp}{\ensuremath{#1}}%
\addtolength{\gvert}{1\gtemp}%
\settodepth{\gtemp}{\ensuremath{#1}}%
\addtolength{\gvert}{1\gtemp}}
\newcommand{\belowpic}[1]{%
\settoheight{\gtemp}{\ensuremath{#1}}%
\addtolength{\gvert}{1\gtemp}%
\addtolength{\gdrop}{1\gtemp}%
\settodepth{\gtemp}{\ensuremath{#1}}%
\addtolength{\gvert}{1\gtemp}%
\addtolength{\gdrop}{1\gtemp}}
\newcommand{\cell}[4]{\put(#1,#2){\makebox(0,0)[#3]{\ensuremath{#4}}}}
\mcm{\zmark}{0}{\scriptstyle{\bullet}}
\newcommand{\prectwo}[3]%
{\begin{picture}(4.2,3.4)(-0.1,-0.2)%
\cell{2}{3.2}{b}{#1}%
\cell{2}{-0.2}{t}{#2}%
\cell{2.2}{1.5}{l}{#3}%
\qbezier(0,2)(2,4)(4,2)%
\qbezier(0,1)(2,-1)(4,1)%
\put(4,2){\vector(1,-1){0}}%
\put(4,1){\vector(1,1){0}}%
\put(2,2.5){\vector(0,-1){2}}%
\end{picture}}
\mcm{\ctwo}{3}{%
\cinitdims{4.2}{3.4}%
\abovepic{#1}%
\belowpic{#2}%
\present{\prectwo{#1}{#2}{#3}}}
\newcommand{\pretopezs}[2]{%
\begin{picture}(2.6,2.3)(-1.3,-2.2)%
\cell{0}{-2.2}{t}{#1}%
\cell{0}{-1.2}{c}{#2}%
\qbezier(0,0)(-2,-2)(0,-2)%
\qbezier(0,0)(2,-2)(0,-2)%
\end{picture}}
\mcm{\topezs}{2}{%
\ginitdims{2.6}{2.3}%
\belowpic{#1}%
\present{\pretopezs{#1}{#2}}}
\newlength{\volt}
\newcommand{\transistor}[5]
{\setlength{\unitlength}{1\volt}
\begin{picture}(18,12)(-5,-6)
\put(0,6){\line(0,-1){12}}
\put(0,-6){\line(3,2){9}}
\put(0,6){\line(3,-2){9}}
\put(9,0){\line(1,0){2}}
\put(-2,4){\line(1,0){2}}
\put(-2,2){\line(1,0){2}}	
\put(-2,-4){\line(1,0){2}}
\put(12,-0.5){\ensuremath{#5}}
\put(-5,3.5){\ensuremath{#2}}
\put(-5,1.5){\ensuremath{#3}}		
\put(-5,-4.5){\ensuremath{#4}}
\thicklines
\put(-1.5,0){\line(1,0){.1}}	
\put(-1.5,-1){\line(1,0){.1}}		
\put(-1.5,-2){\line(1,0){.1}}		
\thinlines
\put(2,-0.5){\ensuremath{#1}}
\end{picture}}
\newcommand{\ctransistor}[5]
	{\raisebox{-6\volt}{\transistor{#1}{#2}{#3}{#4}{#5}}}

\begin{document}

\title{Generalized Enrichment for Categories and Multicategories}
\author{Tom Leinster\\ \\
	\normalsize{Department of Pure Mathematics, University of
	Cambridge}\\ 
	\normalsize{Email: leinster@@dpmms.cam.ac.uk}\\
	\normalsize{Web: http://www.dpmms.cam.ac.uk/$\sim$leinster}}
\date{28 January, 1999\\ 
	\vspace{8mm}
	\normalsize
	\textbf{Abstract}\\
	\vspace{3mm} \raggedright \setlength{\rightskip}{0pt}
In this paper we answer the question: `what kind of a structure can a general
multicategory be enriched in?' (Here `general multicategory' is used in the
sense of \cite{GOM}, \cite{Bur} or \cite{Her}.) The answer is, in a sense to
be made precise, that a multicategory of one type can be enriched in a
multicategory of the type one level up. In particular, we will be able to
speak of a $T_n$-multicategory enriched in a $T_{n+1}$-multicategory, where
$T_n$ is the monad expressing the pasting-together of $n$-opetopes, as
constructed in \cite{SHDCT}.\newline\\
\vspace{-2mm}
The answer for general multicategories reduces to something surprising in the
case of ordinary categories: a category may be enriched in an
`\fc-multicategory', a very general kind of 2-dimensional structure
encompassing monoidal categories, plain multicategories, bicategories and
double categories. It turns out that \fc-multicategories also provide a
natural setting for the bimodules construction. We also explore enrichment
for some multicategories other than just categories. An extended application
is given: the relaxed multicategories of Borcherds and Soibelman are
explained in terms of enrichment.\newline\\}

\maketitle


\tableofcontents
\chapter*{Introduction}
\ucontents{chapter}{Introduction}

At first, the idea of generalizing enrichment from categories to
multicategories might seem rather bland---it's easy to write down
the definition of a plain multicategory enriched in abelian groups, for
instance. But the situation is much more interesting than that, and the
reason why lies in the question `what kind of structure do we enrich in?'

At the most simple level, it is clear enough that we can speak coherently not
just about categories enriched in a monoidal category, but also about
categories enriched in a (plain) multicategory: the essential point is that
the composition morphisms
\[
\homset{\scat{C}}{B}{C} \otimes \homset{\scat{C}}{A}{B}
\go
\homset{\scat{C}}{A}{C}
\]
only have a $\otimes$ in the domain. 

Next, what might a multicategory be enriched in? To answer this, we use the
machinery of general multicategories developed in~\cite{GOM}
and~\cite{SHDCT}, which we assume in this paper. The basic idea there was
that for a suitable monad $T$ on a category \Eee, there is a species of
multicategory (the `$T$-multicategories') in which the shape of the domain of
an arrow is determined by what $T$ is---e.g.\ if $T$ is the free monoid monad
then the domain of an arrow is a \emph{sequence} of objects. In particular,
for each $n$ we constructed the set $S_n$ of $n$-opetopes and the monad $T_n$
on $\Set/S_n$, and in a $T_n$-multicategory the domain of an arrow is a
pasting-together of labelled $n$-opetopes. A $T_0$-multicategory is just a
category, a $T_1$-multicategory is an ordinary multicategory, and beyond that
the $T_n$-multicategories are less familiar structures.

The answer to the question will turn out to be: an ordinary multicategory can
be enriched in a $T_2$-multicategory. So we can, for $n=0$ and 1, speak of a
$T_n$-multicategory enriched in a $T_{n+1}$-multicategory. In fact, this
pattern persists for all $n$---so we end up with a hierarchy of types of
multicategory, in which a multicategory of one type can be enriched in a
multicategory of the next type.

The story doesn't end there. In~\ref{sec:def-enr} we give a definition of
enriched $T$-multicategory for general $T$. When $T$ is taken to be the
identity monad on \Set\ (recalling that a $T$-multicategory is then just a
category) this definition gives a notion of a category enriched in an
`\fc-multicategory'. It turns out that `\fc-multicategories' are a very
general kind of two-dimensional structure, encompassing not only ordinary
multicategories but also bicategories and double categories; we will
therefore be able to speak of categories enriched in any of these
structures. As an offshoot, \fc-multicategories also appear to be the natural
context in which to take (bi)modules, an operation usually performed just on
bicategories.

The paper is laid out as follows. In Chapter~\ref{sec:prelims} we mention
the concepts we'll need from previous papers, and add in a few new
pieces. One of them is the actual definition of enriched
$T$-multicategory. This is not meant to be a self-evidently appropriate
definition; the best justification I can give for it is to show what it means
for various particular $T$'s. Even the most simple case is complex: this is
the subject of Chapter~\ref{sec:cats}, `Enriched Categories', where we unwind
the definition for $\Cartpr=\Zeropr$.

Also included in Chapter~\ref{sec:cats} is the bimodules construction on
\fc-multicategories. This is partly because this chapter is where we first
meet \fc-multicategories, and partly because of the result that any category
enriched in $V$ yields a category enriched in \Bim{V}. The bimodules
construction also produces an interesting structure, $\Bim{\Span}$, with
categories, functors and profunctors all incorporated.

In Chapter~\ref{sec:other} we finally explore enrichment for some
multicategories other than just categories. Enriched plain multicategories
provide the first case, and we spend some time looking at them. In
particular, we consider `operads in a symmetric monoidal category', a notion
used in topology, and see how they fit into our scheme. We also discuss the
above-mentioned hierarchy of types of multicategory.

The final chapter, `Relaxed Multicategories', stands slightly apart; it can
be viewed as an extended application of the earlier theory. Also called
pseudo-monoidal categories, relaxed multicategories have been used by various
authors in the area of quantum algebra and quantum field theory (\cite{BeDr},
\cite{Bor}, \cite{Soi}). We show that relaxed multicategories are just
multicategories enriched in a certain $T_2$-multicategory, which itself
arises naturally from the general theory of multicategories. We then examine
various related structures, including relaxed monoidal categories and
`relaxed categories'.

Throughout the paper we assume the existence of a free multicategory on any
graph of the relevant kind. Conscience obliges me to include at least a
sketch of a construction, but the details are so unimportant that it is
relegated to an appendix.

A final note on size might be appropriate. We basically ignore the
distinction between sets and classes, small and large, and so on: for
instance, an ordinary multicategory officially has a \emph{set} of objects
and a \emph{set} of arrows, but we will frequently talk about the
multicategory of sets or abelian groups. In fact, the structures in
which we enrich are almost always large. I have not attempted to provide any
justification for this, and hope that the reader will not find the issue any
more disturbing than is usual in category theory. 

\subsection*{Acknowledgements}

I am very grateful to Martin Hyland for the encouragement and advice he has
given me in this project. I would also like to thank Craig Snydal for
numerous useful conversations on relaxed multicategories, and especially for
his crucial example \ref{egs:rel-cats}\bref{eg:craig} of a relaxed
category. Richard Borcherds too has given me the benefit of his wisdom on
this topic, and I thank him for this.

This work was supported by a PhD grant from the EPSRC. The document was
prepared in \LaTeX, using Paul Taylor's package for some of the
diagrams. Craig Snydal and Tim Perkins have given me valuable technical
assistance in creating the more complicated diagrams.

\chapter{Preliminaries}	\label{sec:prelims}

In this chapter we gather together the basic concepts on which the paper is
built. Some of them are already written up in Chapters I and IV of
\cite{SHDCT} (and some of them are also in \cite{GOM}): in that case we
simply list the ideas. Others are new, and we go through them
properly. Putting some of these basic concepts together gives us the
definition of enriched multicategory (\ref{sec:def-enr}) which, despite its
brevity, will take most of the rest of the paper to unwind.

\section{Multicategories}	\label{sec:pre-mti}

We recall the basic notions of multicategories from \cite{GOM} or
\cite[I]{SHDCT}, or \cite{Bur} or \cite{Her}, and make some notational
changes from \cite{GOM} and \cite{SHDCT}.
\begin{itemize}
\item
We work with a cartesian monad $T$ on a cartesian category \Eee, and just
write `\Cartpr\ is cartesian'. In \cite{GOM} and \cite{SHDCT}, we wrote
\pr{\cat{S}}{\ust} rather than \Cartpr.

\item	\label{p:gphs}
Given cartesian \Cartpr, there's a bicategory $\Span\Cartpr$ in which a
0-cell is an object of \Eee, a 1-cell $A\go A'$ is a diagram
$(\spn{B}{TA}{A'})$ in \Eee, and 1-cell composition is defined by pullback. A
1-cell of the form $(\spn{B}{TA}{A})$ is called an \Cartpr-graph (on $A$), and
\Cartpr-graphs form a category $\Cartpr\hyph\Gph$. Where possible we drop the
$\Eee$, so we'll speak of $T$-graphs and $T\hyph\Gph$.

\item
An \Cartpr-multicategory (or $T$-multicategory) is a monad in $\Span\Cartpr$;
the category of \Cartpr-multicategories and maps (functors) between them is
called $\Cartpr\hyph\Multicat$ (or $T\hyph\Multicat$). 

\item
A $T$-multicategory $C$ is a diagram $(\spn{C_1}{TC_0}{C_0})$ in \Eee\ with
arrows $C_1\of C_1 \goby{\comp} C_1$ and $C_0 \goby{\ids} C_1$ satisfying
some axioms. We say that $C$ is a multicategory on $C_0$. A $T$-operad is a
$T$-multicategory on 1.

\item 
If $\Cartpr=\Zeropr$ then $T\hyph\Multicat = \Cat$ and a $T$-operad is a
monoid. 

\item
Suppose $\Cartpr=\pr{\Set}{\mr{free\ monoid}}$. Then $T$-multicategories are
called plain multicategories, and similarly operads; these are the original
(non-symmetric) multicategories and operads. An arrow in a plain multicategory
is represented by any one of the pictures 
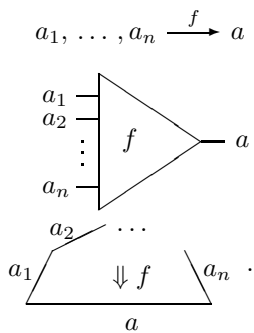
\[
\range{a_1}{a_n} \goby{f} a
\]
\[
\ctransistor{f}{a_1}{a_2}{a_n}{a}
\]
\[
\begin{opetope}
	&	&	&\cnr	&\ldots	&	&	&	\\
	&\cnr	&\ruLine(2,1)<{a_2}&&	&	&\cnr	&	\\
\ruLine(1,2)<{a_1}&&	&	&\Downarrow f&	&	&\rdLine(1,2)>{a_n}\\
\cnr	&	&	&	&\rLine_{a}&	&	&\cnr	\\
\end{opetope}
\ \ \ \ .
\]
When $n=0$, the first version looks like
\[
\cdot \goby{f} a,
\]
the second has no legs on the left-hand (`input') side, and the third is
drawn as
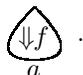
\[
\topezs{a}{\Downarrow\!\!f}\ .
\]

\end{itemize}

We will also need to use transformations between functors between
$T$-multicategories. Generalizing directly from categories:
\begin{defn}	\label{defn:transf}
Let \Cartpr\ be cartesian, let $C$ and $D$ be $T$-multicategories, and let
$C \parpair{F}{G} D$ be functors between them. A \emph{transformation}
$\alpha:F\go G$ is an arrow $\alpha:C_0\go D_1$ such that
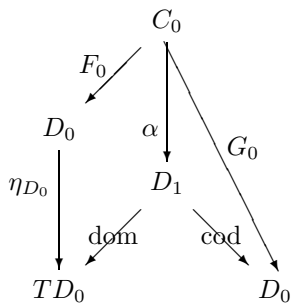
\begin{enumerate}
\item
\begin{slopeydiag}
		&		&C_0		&		&	\\
		&\ldTo<{F_0}	&		&\rdTo(2,5)>{G_0}&	\\
D_0		&		&\dTo<{\alpha}	&		&	\\
\dTo<{\eta_{D_0}}&		&D_1		&		&	\\
		&\ldTo~{\dom}	&		&\rdTo~{\cod}	&	\\
TD_0		&		&		&		&D_0	\\
\end{slopeydiag}
commutes
\item
$D_1 \of D_1 \goby{\comp} D_1$ coequalizes the two maps $C_1 \parpair{}{}
D_1\of D_1$ given by the two diagrams
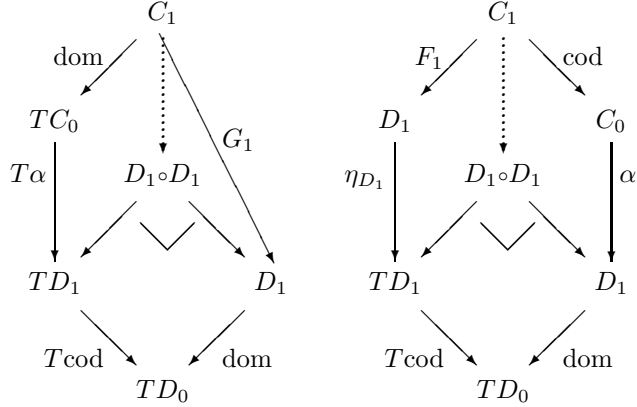
\[
\begin{slopeydiag}
		&		&C_1		&		&	\\
		&\ldTo<{\dom}	&		&\rdTo(2,5)>{G_1}&	\\
TC_0		&		&\dGet		&		&	\\
\dTo<{T\alpha}	&		&D_1\of D_1 \Spbk&		&	\\
		&\ldTo		&		&\rdTo		&	\\
TD_1		&		&		&		&D_1	\\
		&\rdTo<{T\cod}	&		&\ldTo>{\dom}	&	\\
		&		&TD_0		&		&	\\
\end{slopeydiag}
\diagspace
\begin{slopeydiag}
		&		&C_1		&		&	\\
		&\ldTo<{F_1}	&		&\rdTo>{\cod}	&	\\
D_1		&		&\dGet		&		&C_0	\\
\dTo<{\eta_{D_1}}&		&D_1\of D_1 \Spbk&	&\dTo>{\alpha}	\\
		&\ldTo		&		&\rdTo		&	\\
TD_1		&		&		&		&D_1	\\
		&\rdTo<{T\cod}	&		&\ldTo>{\dom}	&	\\
		&		&TD_0		&		&	\\
\end{slopeydiag}
\ \ \ \ .
\]
\end{enumerate}
\end{defn}

\begin{figure}
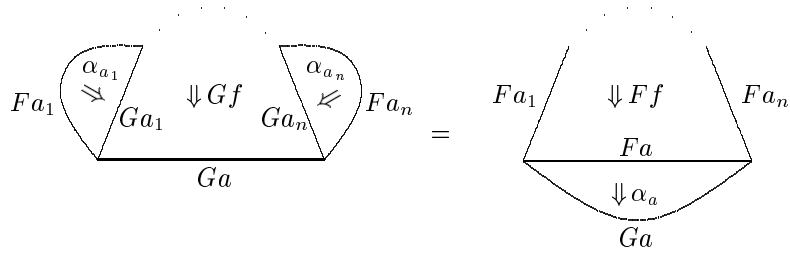

\begin{center}
\piccy{nat}
\end{center}
\caption{Naturality of a transformation}
\label{fig:transf}
\end{figure}

Figure~\ref{fig:transf} illustrates the second axiom in the case of plain
multicategories; the first axiom merely states what the domain and codomain
of $\alpha$ at an object of $C$ are, so is also implicit in the figure.
Later (\ref{egs:bims}\bref{eg:bim-cart-span}) we will find a less
arbitrary-looking definition of transformation.

Transformations can be composed in the usual ways, so that $T\hyph\Multicat$
becomes a 2-category.

Finally for this section, we describe two elementary ways of generating
multicategories. These, together with the free multicategory construction,
are what enable us to define enrichment.

\begin{propn}	\label{propn:basic-constrs}
Let \Cartpr\ be cartesian.
\begin{enumerate}
\item	\label{prop:indiscrete}
There is a functor
\[
I: \Eee \go T\hyph\Multicat
\]
which sends an object $A$ of \Eee\ to a $T$-multicategory with graph
\[
\spaan{TA\times A}{TA}{A}{\mr{pr}_1}{\mr{pr}_2}.
\]
\item	\label{prop:algtomti}
There is a functor
\[
M: \Eee^T \go T\hyph\Multicat
\]
which sends a $T$-algebra $(TA\goby{h}A)$ to a $T$-multicategory with graph
\[
\spaan{TA}{TA}{A}{1}{h}.
\]
\end{enumerate}
\end{propn}

\emph{Remark:} Part~\bref{prop:indiscrete} generalizes the familiar
construction of the indiscrete category on a set. For plain multicategories,
it sends a set $A$ to the multicategory whose objects are the elements of
$A$, and with one arrow $\range{a_1}{a_n} \go a$ for each $\range{a_1}{a_n},
a \in A$. Part~\bref{prop:algtomti} generalizes the construction of a strict
monoidal category (hence a multicategory) from a monoid $A$, in which the
objects of the monoidal category are the elements of $A$, the only arrows are
identities, and the tensor is the multiplication of $A$. 

\pf
\begin{enumerate}
\item
A $T$-multicategory on $A$ is a monoid in the monoidal category of
$T$-graphs on $A$ (this being a one-object sub-bicategory of
\Span\Cartpr). But \spaan{TA\times A}{TA}{A}{\mr{pr}_1}{\mr{pr}_2} is
the terminal $T$-graph on $A$, and a terminal object in a monoidal category
is always a monoid in a unique way, so this graph has a unique multicategory
structure. Extending this to morphisms is straightforward.
\item
Again, there is a unique multicategory structure on the given graph: $TA\of
TA = T^2 A$, and one is forced to put $\comp=\mu_A$ and $\ids=\eta_A$. The
associativity and identity axioms for a multicategory then say exactly that
$(TA\goby{h}A)$ is an algebra. Extending to morphisms is, again,
straightforward.
\done
\end{enumerate}

\section{Free Multicategories}	\label{sec:free-mti}

Just as one can form the free category on a graph, one can form the free
\Cartpr-multicategory on a \Cartpr-graph, as long as \Eee\ and $T$ are
suitably pleasant. In the appendix we give an exact definition of what it
means for \Cartpr\ to be \emph{suitable}, and prove the following result:
\begin{thm}	\label{thm:free-main}
Let \Cartpr\ be suitable. Then the forgetful functor
\[
\Cartpr\hyph\Multicat \go \Eeep = \Cartpr\hyph\Gph
\]
has a left adjoint, the adjunction is monadic, and if $T'$ is the resulting
monad on \Eeep\ then \Cartprp\ is suitable.
\end{thm}

When one takes the free category on an ordinary graph, the collection of
objects (or vertices) is unchanged, and the corresponding fact for
multicategories is encapsulated in a variant of the theorem. If $S$
is an object of \Eee\ then we write $\Cartpr\hyph\Multicat_S$ for the
subcategory of $\Cartpr\hyph\Multicat$ whose objects $C$ have $C_0=S$, and
whose morphisms $F$ have $F_0=1_S$; similarly, we write
$\Eeep_S$ for the category of \Cartpr-graphs on $S$.
\begin{thm}	\label{thm:free-fixed}
Let \Cartpr\ be suitable and let $S\in\Eee$. Then the forgetful functor
\[
\Cartpr\hyph\Multicat_S \go \Eeep_S
\]
has a left adjoint, the adjunction is monadic, and if $T'_S$ is the resulting
monad on $\Eeep_S$ then \pr{\Eeep_S}{T'_S} is suitable.
\end{thm}

All the examples of \Cartpr\ in this paper will be suitable. This can be seen
from these two theorems and the fact that \Zeropr\ is suitable.

\section{The Definition of Enriched Multicategory}	\label{sec:def-enr}

Let $T$ be a monad on a category \Eee. If \Cartpr\ is suitable then we obtain
\Cartprp\ (as in~\ref{sec:free-mti}) which is also suitable. Thus the notion
of `\Cartprp-multicategory' makes sense. The goal of this section is to
define, for any $T'$-multicategory $V$, the category of $T$-multicategories
enriched in $V$.

Applying Proposition~\ref{propn:basic-constrs}, we obtain a composite functor
\[
\Eee \goby{I} \Cartpr\hyph\Multicat \eqv \Eeep^{T'} \goby{M} 
\Cartprp\hyph\Multicat,
\]
for any suitable \Cartpr.

\begin{defn}
Let \Cartpr\ be suitable, let \Cartprp\ be as above, and let $V$ be an
\Cartprp-multicategory. Then a \emph{$T$-multicategory enriched in $V$}
consists of an object $C_0$ of \Eee\ together with a map $MI(C_0) \go V$
in \Cartprp\hyph\Multicat.
\end{defn}
We usually denote the pair \pr{C_0}{MI(C_0)\go V} by the single letter $C$,
and also call $C$ a \emph{$V$-enriched $T$-multicategory}. 

\emph{Warning:}\/ Despite the terminology, there is apparently no
`underlying' $T$-multicategory of an enriched $T$-multicategory in general,
as there is for categories enriched in a monoidal category.

\begin{defn}	\label{defn:enr-ftr}
Let $C$ and $D$ be $V$-enriched \Cartpr-multicategories. Then a
\emph{$V$-enriched functor} from $C$ to $D$ consists of a map $f: C_0 \go
D_0$ in \Eee, together with a transformation
\begin{diagram}[width=2em,height=1em]
MI(C_0)	&	&\rTo^{MI(f)}	&	&MI(D_0)	\\
	&\rdTo(2,5)&		&\ \ldTo(2,5)&		\\
	&	&\ruNT>{\phi}	&	&		\\
	&\	&		&	&		\\
	&	&		&	&		\\
	&	&V.		&	&		\\
\end{diagram}
\end{defn}

We can compose enriched functors by pasting together the transformations, and
thus we obtain a category $\Cartpr_{V}\hyph\Multicat$ of $V$-enriched
\Cartpr-multicategories.

\section{Structured Categories}	\label{sec:struc}

In \cite[4.3]{GOM} and \cite[I.4]{SHDCT} we defined \emph{\Cartpr-structured
categories}, which are to \Cartpr-multicategories as strict monoidal
categories are to plain multicategories. We write $\Cartpr\hyph\Struc$ for
the category of \Cartpr-structured categories, and, as for multicategories,
omit the `\Eee' when we can. Any $T$-structured category has an underlying 
$T$-multicategory, and the forgetful functor $U$ thus defined has a left
adjoint: 
\begin{diagram}[width=2cm]
T\hyph\Struc	&
\pile{\rTo^{U}	\\ \top	\\ \lTo_{F}}	&
T\hyph\Multicat.
\end{diagram}
In the case of plain multicategories and strict monoidal categories, an
object (respectively, arrow) in $FC$ is a sequence of objects (respectively,
arrows) in $C$.

In fact,%
\label{weak-mon-to-mti} 
a monoidal category does not have to be \emph{strict} in order to have an
underlying plain multicategory: any monoidal category will do. If $D$ is the
monoidal category then we define a plain multicategory $C$ with the same
objects as $D$ and with
\[
\multihom{C}{\range{a_1}{a_n}}{a} =
\homset{D}{a_1 \otimes\cdots\otimes a_n}{a}.
\]
In order for this to make sense, $D$ must have $n$-fold tensor products for
all $n$, not just $n=0$ and $n=2$. There are two attitudes we can take to
this. One is to abandon the usual definition of monoidal category, and work
instead with an `unbiased' definition in which $n$-fold tensors are part of
the structure (cf.\ \cite[p.\ 8]{SHDCT}). The other is to use the traditional
definition, but to derive $n$-fold tensors by, for instance, associating to
the left. In both cases there is a canonical isomorphism between any two ways
of tensoring a string of $n$ objects, which means that we can define
composition and identities in $C$, making it into a multicategory. In the
second case there are many ways to choose $n$-fold tensors, but different
choices give isomorphic $C$'s.

Our final observation is that if $D$ and $D'$ are two monoidal categories,
then a lax%
\label{p:lax-ftr-multiftr}
monoidal functor $D\go D'$ is just the same as a map $U(D)\go
U(D')$ of their underlying multicategories. (When $D=1$, this says that a map
$1\go U(D')$ of multicategories is a monoid in $D'$.) Moreover, the
definition of a lax monoidal functor from $D$ to $D'$ makes perfect sense
when $D'$ is any plain multicategory, just by replacing the tensor of $D'$ by
commas, and we still have the result that a multicategory map $U(D)\go D'$ is
the same as a lax monoidal functor from $D$ to $D'$.

\section{Opetopes}	\label{sec:opetopes}

We will need to use the following concepts, as discussed on pages 63--72 of
\cite{SHDCT}:
\begin{itemize}
\item The set $S_n$ of $n$-opetopes and the monad $T_n$ on $\Set/S_n$,
defined recursively and satisfying the clauses on p.~\pageref{p:Cpn-property}
below
\item The geometric representation of opetopes and pasting diagrams
\item $T_n$-structured categories and $T_n$-multicategories
\item The $T_n$-structured category \PD{n} of $n$-pasting diagrams, being the
free $T_n$-structured category on the terminal $T_n$-multicategory
\item $\PD{1}=\Delta$, the simplicial category; $\PD{2}=\fcat{TR}$, so called
because 2-pasting diagrams correspond one-to-one with trees
\item The composite $\tau\of\bftuple{\tau_1}{\tau_n}$ of trees, where
$\tau\in\TR{n}$ is an $n$-leafed tree; the identity tree
$\utree\in\TR{1}$. 
\end{itemize}

\section{Some Two-Dimensional Structures}	\label{sec:two-dims}

In this section we collect together various bicategories and double
categories which we will need later. In a bicategory, $*$ will denote the
horizontal composition of 2-cells.  For the basics of double categories, see
\cite{KS} or even \cite[II.6]{SHDCT}.

\begin{enumerate}
\item 	\label{sec:spans}
We have already met the bicategory $\Span\Cartpr$ in~\ref{sec:pre-mti}, for
cartesian \Cartpr. We write $\Span(\Eee)$ for $\Span\pr{\Eee}{\id}$ and
$\Span$ for $\Span(\Set)$.

\item 	\label{sec:rel}
There is a 2-category \fcat{Rel} of sets and relations, with:
	\begin{description}
	\item[0-cells] sets
	\item[1-cells] relations (i.e.\ a 1-cell $A\rTo B$ is a subset of
	$A\times B$)
	\item[2-cells] inclusions	
	\item[1-cell composition] usual composition of relations.
	\end{description}

\item 	\label{sec:glue}
Let \fcat{Glue} be the sub-bicategory of \Span\ in which all 1-cells are of
the form
\begin{diagram}
\cdot	&\lMonic	&\cdot	&\rMonic	&\cdot	\\
\end{diagram}
(i.e.\ both projections of the span are injective). This is also a
sub-2-category of \fcat{Rel}; a 1-cell $X\go Y$ in \fcat{Glue} is a partial
bijection from $X$ to $Y$. 

\item 	\label{sec:hty-bi}
Any topological space $X$ has a \emph{homotopy bicategory} associated to it.
The objects are points of $X$, the 1-cells are paths
in $X$, and the 2-cells are homotopy classes of path homotopies.

\item 	\label{dbl-cats}
Let $W$ be a 2-category. Then we can construct from $W$ two double
categories, $V$ and $V'$. In both cases, a 0-cell is a 0-cell of $W$; a
horizontal 1-cell is a 1-cell in $W$; a vertical 1-cell is also a 1-cell in
$W$; and a 2-cell inside
\begin{minidiagram}
X		&\rTo^{q}	&X'		\\
\dTo<{p}	&		&\dTo>{p'}	\\
Y		&\rTo_{r}	&Y'		\\
\end{minidiagram}
is either a 2-cell $r\of p \go p'\of q$ in $W$ (in the case of $V$), or a
2-cell $p'\of q \go r\of p$ in $W$ (in the case of $V'$). In both cases, the
composition is obvious.

\item 
There's a one-to-one correspondence between bicategories with precisely one
0-cell and monoidal categories. Given such a bicategory, \Bee, one defines a
monoidal category whose objects are the 1-cells of \Bee\ and whose morphisms
are the 2-cells, and with $b\otimes b' = b\of b'$ and $\beta\otimes\beta' =
\beta*\beta'$, where $b$, $b'$ are 1-cells of \Bee\ and $\beta$, $\beta'$ are
2-cells.

We could equally well have chosen the opposite orientation, so that $b\otimes
b' = b'\of b$ and $\beta\otimes\beta' = \beta'*\beta$. However, we try to
stick to our choice, in Chapters~\ref{sec:cats} and~\ref{sec:other} at
least. The consequence is that `$\otimes$ and $\of$ go in the same
direction': for example, this accounts for the apparently odd reversal of $R$
and $R'$ in \ref{egs:degens}\bref{eg:weak-double}
(page~\pageref{eg:weak-double}).

\item 	\label{deg-mon-cat} 
A strict monoidal category with one object consists of a set $V$ with two
separate monoid structures which share a unit and commute with each
other. This implies that the two monoid structures are equal and commutative:
hence a one-object strict monoidal category is just a commutative monoid.

If we work through this argument for (weak) monoidal categories, we find that
a one-object monoidal category is precisely a commutative monoid equipped
with a distinguished invertible element. The way this comes about is that $V$
has \emph{a priori} three distinguished invertible elements: the
associativity isomorphism $a$, and the left and right unit isomorphisms $l$
and $r$. But the axioms force $l=r$ and $a=1$, so we are just left with $V$
and the invertible $r\in V$, satisfying no further axioms.

\end{enumerate}

\chapter{Enriched Categories}	\label{sec:cats}

The most basic example of enriched multicategories comes when we take
$\Cartpr=\Zeropr$. Since a \Zeropr-multicategory is just a (small) category,
an enriched \Zeropr-multicategory will be called an \emph{enriched
category}. We will see later that this encompasses the traditional definition
of a category enriched in a monoidal category.

If $\Cartpr=\Zeropr$ then the structure $V$ in which we are enriching is a \Onepr-multicategory, where \Gph\
is the category of functors from 
$
\begin{diagram}[width=2em]
(\bullet	&\pile{\rTo\\ \rTo}	&\bullet)	\\
\end{diagram}
$ 
to \Set\ and \fc\ is the free category monad on \Gph. The first task, then,
is to see what an \fc-multicategory is.

\section{\fc-Multicategories}

The graph structure of an \fc-multicategory $V$ looks like
\begin{diagram}
	&	&Y=(Y_1	&\pile{\rTo\\ \rTo}	&Y_0)	&	&	\\
	&	&\ldTo	&			&\rdTo	&	&	\\
\fc(X)=(X_1^*&\pile{\rTo\\ \rTo}&X_0)&		&
X=(X_1 &\pile{\rTo\\ \rTo}&X_0),\\
\end{diagram}
where $X$ and $Y$ are graphs, the $X_i$ and $Y_j$ are sets, $X_1^*$ is the
set of paths in $X$, the horizontal arrows are set maps, and the diagonal
arrows are graph maps. Think of elements of $X_0$ as 0-cells, elements of
$X_1$ as horizontal 1-cells, elements of $Y_0$ as vertical 1-cells, and
elements of $Y_1$ as 2-cells, as in the picture
\begin{equation}	\label{eq:two-cell}
\begin{diagram}
x_0	&\rTo^{\xi_1}	&x_1	&\rTo^{\xi_2}	&\ 	&\cdots	
&\ 	&\rTo^{\xi_n}	&x_n	\\
\dTo<{y}&		&	&		&\Downarrow\,\eta&	
&	&		&\dTo>{y'}\\
x	&		&	&		&\rTo_{\xi}	&	
&	&		&x'	\\
\end{diagram}
\end{equation}
($n\geq 0$, $x_i, x, x' \in X_0$, $\xi_i, \xi \in X_1$, $y, y' \in Y_0$,
$\eta\in Y_1$). The multicategory structure on \spn{Y}{\fc(X)}{X} makes
\spn{Y_0}{X_0}{X_0} into a category (that is, we can compose vertical
1-cells) and gives a composition
\begin{equation}	\label{eq:pasted-two-cells}
\begin{diagram}[width=.5em,height=1em]
\blob&\rTo^{\xi_1^1}&\cdots&\rTo^{\xi_1^{k_1}}&
\blob&\rTo^{\xi_2^1}&\cdots&\rTo^{\xi_2^{k_2}}&\blob&
\ &\cdots&\ &
\blob&\rTo^{\xi_n^1}&\cdots&\rTo^{\xi_n^{k_n}}&\blob\\
\dTo<{y_0}&&\Downarrow\eta_1&&
\dTo&&\Downarrow\eta_2&&\dTo&
\ &\cdots&\ &
\dTo&&\Downarrow\eta_n&&\dTo>{y_n}\\
\blob&&\rTo_{\xi_1}&&
\blob&&\rTo_{\xi_2}&&\blob&
\ &\cdots&\ &
\blob&&\rTo_{\xi_n}&&\blob\\
\dTo<{y}&&&&&&&&\Downarrow\eta &&&&&&&&\dTo>{y'}\\
\blob&&&&&&&&\rTo_{\xi}&&&&&&&&\blob\\
\end{diagram}
\end{equation}
$\goesto$
\begin{diagram}[width=.5em,height=1em]
\blob&\rTo^{\xi_1^1}&\ &&
&&&&\cdots&
&&&
&&\ &\rTo^{\xi_n^{k_n}}&\blob\\
&&&&&&&&&&&&&&&&\\
\dTo<{y\of y_0}&&&&&&&&\Downarrow\eta\of\tuple{\eta_1}{\eta_2}{\eta_n}
&&&&&&&&\dTo>{y'\of y_n}\\
&&&&&&&&&&&&&&&&\\
\blob&&&&&&&&\rTo_{\xi}&&&&&&&&\blob\\
\end{diagram}
($n\geq 0, k_i\geq 0$, with \blob's representing elements of $X_0$) and
identities
\[
\begin{diagram}[width=1em,height=2em]
x&\rTo^{\xi}&x'\\
\end{diagram}
\diagspace\goesto\diagspace
\begin{diagram}[width=1em,height=2em]
x&\rTo^{\xi}&x'\\
\dTo<{1_x}&\Downarrow 1_{\xi}&\dTo>{1_{x'}}\\
x&\rTo_{\xi}&x'.\\
\end{diagram}
\]
The composition and identities obey associativity and identity laws.

The pictures in the nullary case are worth a short comment.
When $n=0$, the 2-cell of diagram~\bref{eq:two-cell} is drawn as
\begin{diagram}
x_0		&\rEquals		&x_0		\\
\dTo<{y}	&\Downarrow\,\eta	&\dTo>{y'}	\\
x		&\rTo_{\xi}		&x',		\\
\end{diagram}
and the diagram~\bref{eq:pasted-two-cells} of pasted-together 2-cells is
drawn as
\begin{diagram}[size=1.5em]
w_0		&\rEquals		&w_0		\\
\dTo<{y_0}	&=			&\dTo>{y_0}	\\
x_0		&\rEquals		&x_0		\\
\dTo<{y}	&\Downarrow\,\eta	&\dTo>{y'}	\\
x		&\rTo_{\xi}		&x'.		\\
\end{diagram}
The composite%
\label{p:null-notation}
of this last diagram will be written as $\eta\of y_0$.

As such, \fc-multicategories are not familiar, but various degenerate cases
are. These are set out below, and summarized in Table~\ref{table:degens}.
\begin{table}
\hspace{-16mm}%
\begin{tabular}{r|c|c|c}
			&not `representable'	&`representable'	
&`uniformly representable'\\
\hline
no vertical degeneracy	&\fc-multicategory	&weak double category
&double category	\\
\hline
vertically discrete	&vertically discrete	&bicategory
&2-category	\\
 			&\fc-multicategory	&
&		\\
\hline
vertically trivial	&plain multicategory	&monoidal category	
&strict monoidal \\
			&			&
&category
\end{tabular}
\caption{Some of the possible degeneracies of an \fc-multicategory. The entries
are explained on pages \pageref{egs:degens}--\pageref{end-degens}.}
\label{table:degens}
\end{table}

\begin{eg}{egs:degens}

\item	\label{eg:strict-double}
Any double category (see \ref{sec:two-dims}) has an underlying
\fc-multicategory, defined by saying that a 2-cell
\begin{diagram}
x_0	&\rTo^{\xi_1}	&x_1	&\rTo^{\xi_2}	&\ 	&\cdots	
&\ 	&\rTo^{\xi_n}	&x_n	\\
\dTo<{y}&		&	&		&\Downarrow\,\eta&	
&	&		&\dTo>{y'}\\
x	&		&	&		&\rTo_{\xi}	&	
&	&		&x'	\\
\end{diagram}
in the \fc-multicategory is just a 2-cell
\begin{diagram}
x_0	&\rTo^{\xi_n \of \cdots \of \xi_1}	&x_n\\
\dTo<{y}&\Downarrow			&\dTo>{y'}\\
x	&\rTo_{\xi}			&x'\\
\end{diagram}
in the double category.

\item	\label{eg:weak-double}
In fact,~\bref{eg:strict-double} works even when the double category is
`horizontally weak'. A typical example of such a structure has rings (not
necessarily commutative) as its 0-cells, bimodules as its horizontal 1-cells,
ring homomorphisms as its vertical 1-cells, and `homomorphisms of bimodules
with respect to the vertical changes of base' as 2-cells. In other words, a
2-cell looks like
\begin{diagram}
R	&\rMod^{M}	&R'	\\
\dTo<{f}&\Downarrow\phi	&\dTo>{f'}\\
S	&\rMod_{N}	&S',	\\
\end{diagram}
where $R$, $R'$, $S$, $S'$ are rings, $M$ is an \pr{R'}{R}-bimodule (i.e.\
simultaneously a left $R'$-module and a right $R$-module) and $N$ similarly,
$f$ and $f'$ are ring homomorphisms, and $\phi: M \go N$ is an abelian group
homomorphism such that
\[
\phi(r'.m.r) = f'(r').\phi(m).f(r).
\]
(The crossed arrow is used here to indicate a bimodule.)  Composition of
horizontal 1-cells is tensor, composition of vertical 1-cells is the usual
composition of ring homomorphisms, and composition of 2-cells is defined in
an evident way. The essential point is that although the 0-cells and vertical
1-cells form a category, the same cannot be said of the horizontal structure:
tensor only obeys the associative and unit laws up to coherent
isomorphism. We do not bother to write down the full definition of a weak
double category, since it is just an easy extension of the definition of a
bicategory.

(In order to have a 1-cell `$\xi_n\of\cdots\of\xi_1$', as in the second
diagram of~\bref{eg:strict-double}, we must either define weak double
category in an unbiased manner, or else choose a particular $n$-fold
composition. This is exactly the same issue as was discussed on
page~\pageref{weak-mon-to-mti} for monoidal categories and plain
multicategories.)

Another example has small categories as 0-cells, profunctors as horizontal
1-cells, functors as vertical 1-cells, and `morphisms of profunctors with
respect to the vertical functors' as 2-cells.  We will explore both these
examples further, and give proper definitions, in section~\ref{sec:bim}.

\item	\label{eg:vdfcm}
Suppose that all vertical 1-cells are the identity, i.e.\ $Y_0=X_0$ and
\[
(\spn{Y_0}{X_0}{X_0}) = (\spaan{X_0}{X_0}{X_0}{1}{1}), 
\]
so that the whole graph looks like
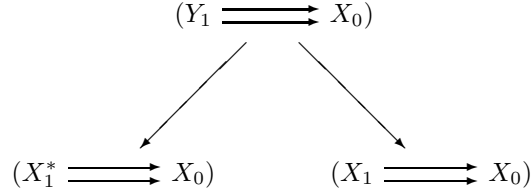
\begin{diagram}
	&	&(Y_1	&\pile{\rTo\\ \rTo}	&X_0)	&	&	\\
	&	&\ldTo	&			&\rdTo	&	&	\\
(X_1^*&\pile{\rTo\\ \rTo}&X_0)&		&(X_1	&\pile{\rTo\\ \rTo}&X_0)\\
\end{diagram}
or in another format,
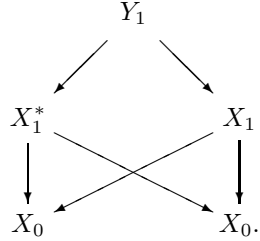
\begin{slopeydiag}
	&	&Y_1	&	&	\\
	&\ldTo	&	&\rdTo	&	\\
X_1^*	&	&	&	&X_1	\\
\dTo	&\rdTo(4,2)&	&\ldTo(4,2)&\dTo\\
X_0	&	&	&	&X_0.	\\
\end{slopeydiag}
Then we call the \fc-multicategory \emph{vertically discrete} (since the
category of 0-cells and vertical 1-cells is discrete). It consists of some
objects $x, x', \ldots$, some 1-cells $\xi, \xi',\ldots$, and some 2-cells
looking like
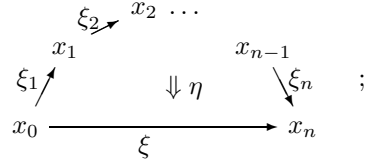
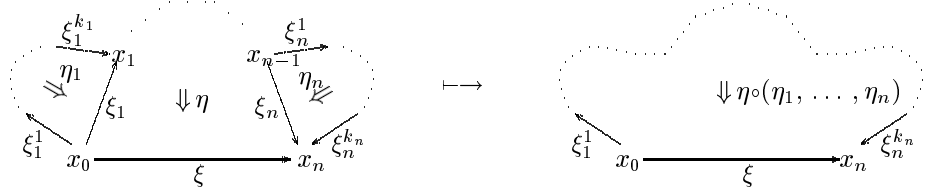
\begin{lopetope}
		&	&	&x_2	&\ldots	&	&	&	\\
		&x_1	&\ruEdge(2,1)^{\xi_2}&&	&	&x_{n-1}&	\\
\ruEdge(1,2)<{\xi_1}&	&	&	&\Downarrow \eta&&	
&\rdEdge(1,2)>{\xi_{n}\ \ \ \ \ ;}\\
x_0		&	&	&\rEdge_{\xi}&	&	&	&x_n	\\
\end{lopetope}
there is a composition 
\[
\piccy{vdcomp}
\]
and an identity function
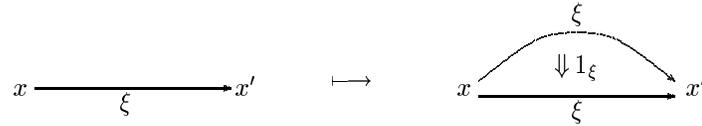
\[
\piccy{vdid}
\]
which obey the inevitable associativity and identity laws.

\item
Any bicategory gives rise to a vertically discrete \fc-multicategory, in the
same way that any weak double category gives rise to an \fc-multicategory
(and with similar observations on $n$-fold composites).

\item
Suppose we are back in the situation of~\bref{eg:vdfcm}, with the extra
condition that there's only one object. Thus $X_0=Y_0=1$, and the graph looks
like 
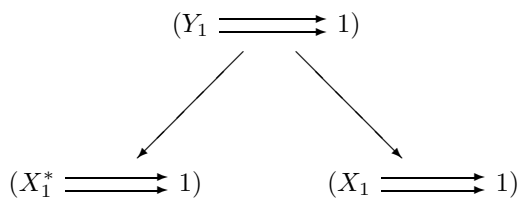
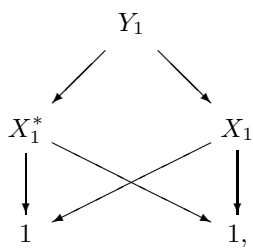
\begin{diagram}
	&	&(Y_1	&\pile{\rTo\\ \rTo}	&1)	&	&	\\
	&	&\ldTo	&			&\rdTo	&	&	\\
(X^*_1&\pile{\rTo\\ \rTo}&1)&		&(X_1	&\pile{\rTo\\ \rTo}&1)\\
\end{diagram}
or
\begin{slopeydiag}
	&	&Y_1	&	&	\\
	&\ldTo	&	&\rdTo	&	\\
X^*_1	&	&	&	&X_1	\\
\dTo	&\rdTo(4,2)&	&\ldTo(4,2)&\dTo\\
1	&	&	&	&1,	\\
\end{slopeydiag}
where $X_1^*$ is the free monoid on $X_1$. The \fc-multicategory then
consists of some `horizontal 1-cells' $\xi, \xi', \ldots$ and some `2-cells'
$\eta, \eta', \ldots$, looking like
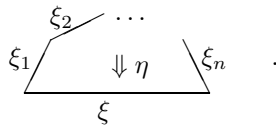
\[
\begin{opetope}
		&	&	&\cnr 	&\ldots	&	&	&	\\
		&\cnr	&\ruLine(2,1)^{\xi_2}&&	&	&\cnr 	&	\\
\ruLine(1,2)<{\xi_1}&	&	&	&\Downarrow \eta&&	
&\rdLine(1,2)>{\xi_{n}\ \ \ \ \ .}\\
\cnr		&	&	&\rLine_{\xi}&	&	&	&\cnr	\\
\end{opetope}
\]
In other words, it's just a plain multicategory, where the $\xi$'s are now
objects and the $\eta$'s are arrows. So a plain multicategory is a special
kind of \fc-multicategory.

\item
Specializing the previous example further, we finally see that a monoidal
category is a rather special kind of \fc-multicategory, since any (strict or
weak) monoidal category has an underlying plain multicategory
(see \ref{sec:struc}). So it will make sense in our language to speak of
a category enriched in a monoidal category. Alternatively, a monoidal
category is a one-object bicategory, and any bicategory is an
\fc-multicategory, so a monoidal category becomes an \fc-multicategory in
this way too. (It makes no difference, up to isomorphism, whether we obtain an
\fc-multicategory from a given monoidal category by going via plain
multicategories or bicategories.)

\item
A one-object monoidal category is a commutative monoid with a distinguished
invertible element (see page~\pageref{deg-mon-cat}). We will therefore be
able to speak of a category enriched in such a structure. This is, of course,
encompassed in the traditional definition of enriched category.%
\label{end-degens}

\end{eg}

\section{Enriched Categories: Elementary Description}

Having seen what a general \fc-multicategory is, we can now describe
what $MI(A)$ is, for a set $A$. Firstly, $I(A)=(\spn{A\times A}{A}{A})$ and
\[
\fc(I(A)) = (\spaan{\coprod_{n\geq 1}A^n}{A}{A}{\mr{first}}{\mr{last}}). 
\]
So, the graph of $MI(A)$ is
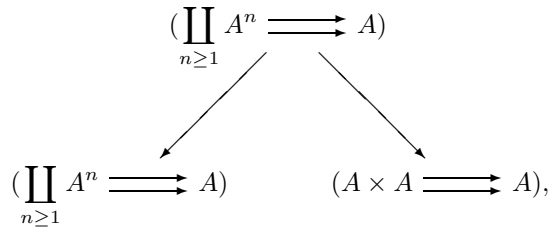
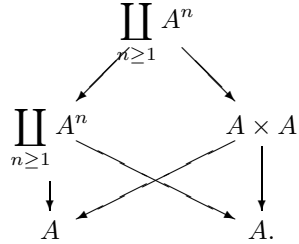
\begin{diagram}
	&	&(\coprod_{n\geq 1}A^n	&\pile{\rTo\\ \rTo}	&A)	&
	&	\\
	&	&\ldTo	&			&\rdTo	&	&	\\
(\coprod_{n\geq 1}A^n&\pile{\rTo\\ \rTo}&A)&		
&(A\times A	&\pile{\rTo\\ \rTo}&A),\\
\end{diagram}
or
\begin{slopeydiag}
			&	&\coprod_{n\geq 1}A^n	&	&	\\
			&\ldTo	&			&\rdTo	&	\\
\coprod_{n\geq 1}A^n	&	&			&	&A\times A\\
\dTo			&\rdTo(4,2)&			&\ldTo(4,2)&\dTo\\
A			&	&			&	&A.	\\
\end{slopeydiag}
The 0-cells of $MI(A)$ are the elements of $A$; the only vertical 1-cells are
identities (so it is vertically discrete); there is one horizontal 1-cell
$a\go b$ for each $a, b\in A$; and for each $a_0, \ldots, a_n$ ($n\geq 0$)
there is a unique 2-cell
\begin{diagram}
a_0	&\rTo		&a_1	&\rTo		&\ 	&\cdots	
&\ 	&\rTo		&a_n	\\
\dTo<{1}&		&	&		&\Downarrow	&	
&	&		&\dTo>{1}\\
a_0	&		&	&		&\rTo		&	
&	&		&a_n.\\
\end{diagram}
The composition and identities are uniquely determined.

We can now give an elementary description of a category enriched in an
\fc-multicategory $V$. It consists of a set $C_0$ (`of objects') together
with a map $MI(C_0) \go V$ of \fc-multicategories: that is, for each $a\in
C_0$ a 0-cell \eend{C}{a} of $V$, for each $a,b \in C_0$ a horizontal 1-cell
\[
\ehom{C}{a}{b}: \eend{C}{a} \go \eend{C}{b}
\]
in $V$, and for each sequence $a_0, \ldots, a_n$ ($n\geq 0$) in $A$ a 2-cell
\begin{diagram}
\eend{C}{a_0}&\rTo^{\ehom{C}{a_0}{a_1}}&\eend{C}{a_1}	&\rTo^{\ehom{C}{a_1}{a_2}}	
&\ 	&\cdots	&\ 	&\rTo^{\ehom{C}{a_{n-1}}{a_n}}	&\eend{C}{a_n}	\\
\dTo<{1}&		&	&		
&\Downarrow\,\kappa_{a_0,\ldots a_n}	&	&	&	&\dTo>{1}\\
\eend{C}{a_0}&	&	&	&\rTo_{\ehom{C}{a_0}{a_n}}	&	
&	&		&\eend{C}{a_n}	\\
\end{diagram}
in $V$ such that
\label{p:kappas-closed}
\begin{eqnarray*}
&
\begin{diagram}[width=.5em,height=1em]
\blob&\rTo&\cdots&\rTo&
\blob&\rTo&\cdots&\rTo&\blob&
\ &\cdots&\ &
\blob&\rTo&\cdots&\rTo&\blob\\
\dTo&&\Downarrow\,\kappa&&
\dTo&&\Downarrow\,\kappa&&\dTo&
\ &\cdots&\ &
\dTo&&\Downarrow\,\kappa&&\dTo\\
\blob&&\rTo&&
\blob&&\rTo&&\blob&
\ &\cdots&\ &
\blob&&\rTo&&\blob\\
\dTo&&&&&&&&\Downarrow\,\kappa &&&&&&&&\dTo\\
\blob&&&&&&&&\rTo&&&&&&&&\blob\\
\end{diagram}
\\
=&
\begin{diagram}[width=.5em,height=1em]
\blob&\rTo&\ &&
&&&&\cdots&
&&&
&&\ &\rTo&\blob\\
&&&&&&&&&&&&&&&&\\
\dTo&&&&&&&&\Downarrow\,\kappa
&&&&&&&&\dTo\\
&&&&&&&&&&&&&&&&\\
\blob&&&&&&&&\rTo&&&&&&&&\blob\\
\end{diagram}
\end{eqnarray*}
and
\[
\begin{diagram}[width=2em,height=2em]
\eend{C}{a_0}	&\rTo^{\ehom{C}{a_0}{a_1}}	&\eend{C}{a_1}	\\
\dTo<{1}	&\Downarrow\,1			&\dTo>{1}	\\
\eend{C}{a_0}	&\rTo_{\ehom{C}{a_0}{a_1}}	&\eend{C}{a_1}	\\
\end{diagram}
\diagspace = \diagspace
\begin{diagram}[width=2em,height=2em]
\eend{C}{a_0}	&\rTo^{\ehom{C}{a_0}{a_1}}	&\eend{C}{a_1}	\\
\dTo<{1}	&\Downarrow\,\kappa_{a_0,a_1}	&\dTo>{1}	\\
\eend{C}{a_0}	&\rTo_{\ehom{C}{a_0}{a_1}}	&\eend{C}{a_1}.	\\
\end{diagram}
\]
An easy induction shows that it is equivalent to specify just the binary and
nullary%
\label{p:finitary}
composites in $C$. Thus a category enriched in $V$ consists of:
\begin{itemize}
\item a set $C_0$
\item for each $a\in C_0$, a 0-cell $\eend{C}{a}$ of $V$
\item for each $a,b\in C_0$, a horizontal 1-cell $\eend{C}{a} 
\goby{\ehom{C}{a}{b}} \eend{C}{b}$ in $V$
\item for each $a,b,c\in C_0$, a 2-cell
\begin{diagram}
\eend{C}{a}&\rTo^{\ehom{C}{a}{b}}&\eend{C}{b}	&\rTo^{\ehom{C}{b}{c}}
&\eend{C}{c}	\\
\dTo<{1}&		&\Downarrow\,\kappa_{a,b,c}	&	&\dTo>{1}\\
\eend{C}{a}&		&\rTo_{\ehom{C}{a}{c}}		&	&\eend{C}{c}\\
\end{diagram}
\item for each $a\in C_0$, a 2-cell
\begin{diagram}
\eend{C}{a}	&\rEquals		&\eend{C}{a}	\\
\dTo<{1}	&\Downarrow\,\kappa_{a}	&\dTo>{1}	\\
\eend{C}{a}	&\rTo_{\ehom{C}{a}{a}}	&\eend{C}{a}	\\
\end{diagram}
\end{itemize}
obeying the usual associativity and identity laws (as for a monad).

When we ask what a category enriched in an \fc-multicategory $V$ is, we do
not have to consider any of the vertical 1-cells of $V$ apart from the
identities, since $MI(C_0)$ is vertically discrete for any set $C_0$. Thus if
$\bar{V}$ is the underlying vertically discrete \fc-multicategory of
$V$---that is, $V$ with all non-identity vertical 1-cells discarded---then a
category enriched in $V$ is the same thing as a category enriched in
$\bar{V}$. The vertical 1-cells in $V$ will, however, play a
significant role when we move on to discuss \emph{maps} between $V$-enriched
categories.

Before giving specific examples of enriched categories, we state what the
definition of a $V$-enriched category reduces to when $V$ is degenerate in
some way such as those listed above. By the foregoing comments, we may
assume immediately that $V$ is vertically discrete.

\begin{eg}{egs:enr-deg}

\item	\label{eg:enr-bi}
Let $V$ be a bicategory.
Then a category enriched in $V$ consists of:
	\begin{itemize}
	\item a set $C_0$
	\item for each $a\in C_0$, an object $\eend{C}{a}$ of $V$
	\item for each $a,b\in C_0$, a 1-cell $\eend{C}{a}
	\goby{\ehom{C}{a}{b}} \eend{C}{b}$ in $V$
	\item for each $a,b,c\in C_0$, a 2-cell
	\[
	\eend{C}{a} \ctwo{\ehom{C}{b}{c}\of\ehom{C}{a}{b}}{\ehom{C}{a}{c}}{}
	\eend{C}{c}
	\]
	\item for each $a\in C_0$, a 2-cell
	\[
	\eend{C}{a} \ctwo{1}{\ehom{C}{a}{a}}{} \eend{C}{a}
	\]
	\end{itemize}
such that these 2-cells satisfy associativity and identity axioms. 

This notion of a category enriched in a bicategory is just the same as that
of Walters et al (see \cite{BCSW}, \cite{CKW}, \cite{Wal}).
It is also the same as B\'{e}nabou's notion of a polyad. (A \emph{polyad} in a
bicategory $V$ consists of a set $C_0$ and a lax functor $\cat{I}(C_0) \go
V$, where $\cat{I}(C_0)$ is the 2-category whose object-set is $C_0$ and each
of whose hom-categories is \fcat{1}. See \cite{Ben} for details.) Polyads
were so called because a one-object polyad is a monad (in the bicategory
$V$); we can also see directly from our definition that a category enriched
in $V$ is just a monad in $V$ when $C_0 = 1$.

\item
When $V$ is a monoidal category, our definition of a category enriched in $V$
coincides with the traditional one. This can be seen by regarding $V$ as a
bicategory with one object and using the data and axioms in~\bref{eg:enr-bi}.

\item	\label{eg:deg-plain}
Degenerating in a different direction, let $V$ be a plain multicategory. Then
a $V$-enriched category consists of
	\begin{itemize}
	\item a set $C_0$
	\item for each $a,b,c\in C_0$, an arrow
	\begin{diagram}
	\ehom{C}{a}{b},\ehom{C}{b}{c}&\rTo^{\kappa_{a,b,c}}&\ehom{C}{a}{c}\\
	\end{diagram}
	in $V$
	\item for each $a\in C_0$, an arrow
	\begin{diagram}
	\cdot &\rTo^{\kappa_{a}} &\ehom{C}{a}{a}\\
	\end{diagram}
	in $V$ (where $\cdot$ denotes the empty sequence)
	\end{itemize}
such that the associativity and identity axioms hold:
\begin{center}
$\kappa_{acd} \of \pr{\kappa_{abc}}{1} =
\kappa_{abd} \of \pr{1}{\kappa_{bcd}}$	\\
$\kappa_{aab} \of \pr{\kappa_a}{1} = 1 =
\kappa_{abb} \of \pr{1}{\kappa_b}$	\\
\end{center}
for all $a,b,c,d\in C_0$ (Figure~\ref{fig:ass-and-id}).

\begin{figure}
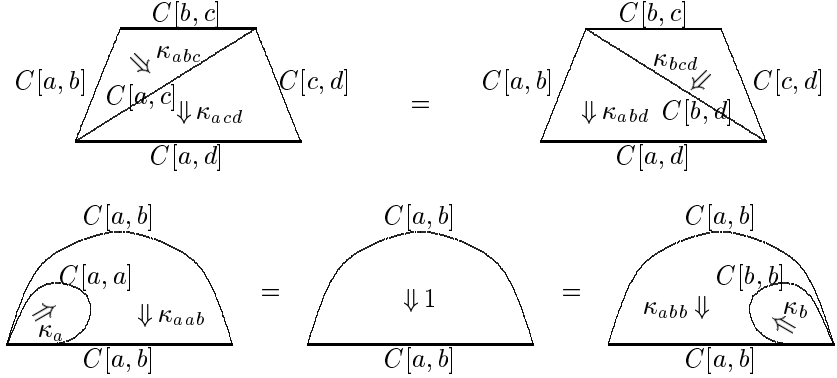

\centerline{\piccy{assid}}
\caption{Associativity and identity axioms for a category enriched in a plain
multicategory}
\label{fig:ass-and-id}
\end{figure}

Of course, enrichment in a monoidal category is a special case of this.

\item
A plain operad $V$ is a one-object plain multicategory, so we can speak of a
category enriched in $V$. Such consists of 
	\begin{itemize}
	\item a set $A$
	\item for each $a,b,c\in A$, an element \etrip{a}{b}{c} of $V(2)$
	\item for each $a\in A$, an element \esing{a} of $V(0)$
	\end{itemize}
satisfying the equations
\begin{center}
$\etrip{a}{c}{d} \of (\etrip{a}{b}{c},1) = 
\etrip{a}{b}{d} \of (1,\etrip{b}{c}{d})$	\\
$\etrip{a}{a}{b} \of (\esing{a},1) = 1 = \etrip{a}{b}{b} \of (1,\esing{b})$
\end{center}
for all $a,b,c,d \in A$.

\item	\label{eg:comm-mon}
As observed on page~\pageref{deg-mon-cat}, a one-object monoidal
category is a commutative monoid $V$ with a distinguished invertible element
$r$. Either by treating this as a special case of a monoidal category, or by
looking at its underlying operad, we can see that a category enriched in
\pr{V}{r} consists of a set $A$ and functions
\begin{center}
$\etrip{\dashbk}{\dashbk}{\dashbk} : A\times A\times A \go V$	\\
$\esing{\dashbk} : A \go V$	
\end{center}
satisfying the equations
\begin{center}
$\etrip{a}{c}{d} + \etrip{a}{b}{c} =
\etrip{a}{b}{d} + \etrip{b}{c}{d}$	\\
$\etrip{a}{a}{b} + \esing{a} = r = \etrip{a}{b}{b} + \esing{b}$
\end{center}
for all $a,b,c,d \in A$.

\item
Any category can be regarded as a locally discrete 2-category, so in this
sense we can speak of a category enriched in a category. That is: fix a
category $W$, and let $V$ be the 2-category whose 0- and 1-cells form the
category $W$, and all of whose 2-cells are identities. Then a category
enriched in $V$ works out to be a set $C_0$ together with a functor $I(C_0)
\go W$, where $I(C_0)$ is the indiscrete category on $C_0$. 

\item
Similarly, any category $W$ can be regarded as a locally indiscrete
2-category $V$---thus the 0- and 1-cells of $V$ form the category $W$, and
there is precisely one 2-cell between any two parallel 1-cells. A category
enriched in $V$ is a set $C_0$ together with a graph map $I(C_0) \go W$. 

\end{eg}

\section{Enriched Categories: Examples} \label{sec:specific-enr-cats}

We now proceed to some specific examples of categories enriched in
\fc-multicategories. Along the way, we define some \fc-multicategories which
for now will be vertically discrete, but later will be redefined with more
interesting vertical structure. As observed before, the vertical structure is
immaterial when we are just considering enriched categories (rather than
functors), so this redefinition will do no harm.

\begin{eg}{egs:specific-enr-cats}

\item	\label{eg:spec-bim-ab}
Write $\Bim{\Ab}$ for the weak double category of rings and bimodules,
described in Example~\ref{egs:degens}\bref{eg:weak-double}. Any category $C$
enriched in the monoidal category \triple{\Ab}{\otimes}{\integers} gives rise
to a category $C'$ enriched in $\Bim{\Ab}$. This really just amounts to the
observations that if $a$ is an object of $C$ then the abelian group
\ehom{C}{a}{a} is in fact a ring, and that if $a$ and $b$ are objects of
$C$ then the abelian group \ehom{C}{a}{b} is both a left
\ehom{C}{b}{b}-module and a right \ehom{C}{a}{a}-module. Thus we take
$C'_0$ to be the set of objects of $C$, define \eend{C'}{a} to be
\ehom{C}{a}{a} with ring structure given by composition, and define
\ehom{C'}{a}{b} to be \ehom{C}{a}{b} with bimodule structure given by
composition. 

\item	\label{eg:spec-bim-set}
The additive structure played no essential role in the previous example. So,
write $\Bim{\Set}$ for the weak double category whose 0-cells are monoids,
whose horizontal 1-cells are sets with commuting left and right actions by
the monoids on either side, whose vertical 1-cells are monoid homomorphisms,
and the rest of whose structure is also defined analogously to that of
$\Bim{\Ab}$. (A rigorous version of this will be given in
section~\ref{sec:bim}.) Then \emph{any} category naturally gives rise to a
category enriched in $\Bim{\Set}$, since any category is naturally enriched
in the monoidal category \triple{\Set}{\times}{1}. Explicitly, if $C$ is a
(small) category then let $C'_0$ be the set of objects of $C$, let \eend{C'}{a}
be \homset{C}{a}{a} with monoid structure given by composition, and let
\ehom{C'}{a}{b} be \homset{C}{a}{b} with the obvious actions by \eend{C'}{a}
and \eend{C'}{b}. 

\item
As a variation on the previous example, let $C$ be any category, and again
define a category $C'$ enriched in $\Bim{\Set}$. This time, let \eend{C'}{a}
be the automorphism group $\fcat{Aut}(a)$ of $a$ in $C$, rather than the set
of all endomorphisms of $a$. Everything else is left the same.

\item	\label{eg:span}
Let \Span\ be the bicategory of spans (in \Set;
see~\ref{sec:two-dims}\bref{sec:spans}). A category enriched in \Span\
consists of
\begin{itemize}
\item a set $C_0$
\item for each $a\in C_0$, a set \eend{C}{a}
\item for each $a,b\in C_0$, a span
\[
\spn{\ehom{C}{a}{b}}{\eend{C}{a}}{\eend{C}{b}}
\]
\item for each $a,b,c\in C_0$, a function
\[
\comp: \ehom{C}{a}{b} \times_{\eend{C}{b}} \ehom{C}{b}{c} \go
\ehom{C}{a}{c}
\]
\item for each $a\in C_0$, a function
\[
\ids: \eend{C}{a} \go \ehom{C}{a}{a}
\]
\end{itemize}
such that these functions satisfy domain, codomain, associativity and
identity axioms. In other words, it is a category $D$, a set $I$, and a
function $D_0 \go I$, where $D_0$ is the set of objects of $D$. Here we have
taken $I=C_0$, $D_0 = \coprod_{a\in C_0} \eend{C}{a}$, $D_1 = \coprod_{a,b\in
C_0} \ehom{C}{a}{b}$, projection as the function $D_0 \go I$, and the obvious
category structure on $D$.

\item	\label{eg:span-actual}
For some actual categories enriched in \Span, we could take $D$ to be the
category of real manifolds and continuous functions, $I=\nat$, and $D_0 \go
I$ to be the function assigning dimension. Similarly, we could take $D$ to be
the category of finite-dimensional vector spaces over a fixed field, and
again $I=\nat$ and the dimension function. Or we could take a skeleton of
this: let $D$ be the category whose objects are the natural numbers, and with
\[
\homset{D}{m}{n} = M_{n,m} = \{n\times m\ \mr{matrices} \}
\]
and $D_0 \go I$ the identity function on \nat. More abstractly, we could take
$D$ to be any small category, $I$ to be either the set of objects of $D$ or
the set of isomorphism classes of objects of $D$, and the natural function
$D_0 \go I$ in each case. 

\item	\label{eg:spec-bim-span}
Let \Prof\ be the weak double category mentioned at the end of
Example~\ref{egs:degens}\bref{eg:weak-double}, whose underlying bicategory is
the `bicategory of profunctors'. A category enriched in \Prof\ consists of
	\begin{itemize}
	\item a set $C_0$
	\item for each $a\in C_0$, a category \eend{C}{a}
	\item for each $a,b\in C_0$, a profunctor $\ehom{C}{a}{b}:
	\eend{C}{a} \rMod \eend{C}{b}$
	\item for each $a,b,c\in C_0$, a morphism
	\[
	\ehom{C}{b}{c}\otimes\ehom{C}{a}{b} \go \ehom{C}{a}{c}
	\]
	of profunctors
	\item for each $a\in C_0$, a morphism
	\[	
	\eend{C}{a} \go \ehom{C}{a}{a}
	\]
	of profunctors, where on the left-hand side, \eend{C}{a} denotes the
	identity profunctor on \eend{C}{a};
	\end{itemize}
these morphisms of profunctors are to satisfy associativity and identity
laws. For instance:
	\begin{itemize}
	\item $C_0 = \nat$, \eend{C}{n} is the category of $n$-dimensional
	vector spaces over a fixed field, and the functor 
	\[
	\ehom{C}{m}{n}: \eend{C}{m}^{\op} \times \eend{C}{n} \go \Set
	\]
	sends a pair \pr{U}{W} to the set of linear maps $U\go W$.
	\item $C_0 = \nat$ and \eend{C}{n} is the multiplicative monoid
	$M_{n,n}$ of $n\times n$ matrices
	with coefficients in a fixed field, viewed as a category. A functor
	\[
	\eend{C}{m}^{\op} \times \eend{C}{n} \go \Set
	\]
	is a set with compatible left action by $M_{n,n}$ and right action by
	$M_{m,m}$; we thus take the profunctor \ehom{C}{m}{n} to be the set
	$M_{n,m}$ with actions by matrix multiplication.
	\item In the previous example, we could change \eend{C}{n} from
	$M_{n,n}$ to the general linear group $GL_{n}$ (or any other
	submonoid of $M_{n,n}$), and leave the rest of the definition the
	same. 
	\end{itemize}
The first two of these examples arise in a mechanical way from the examples
of categories enriched in \Span\ (\bref{eg:span-actual} above). We will see
later (\ref{sec:bim}) that $\Prof = \Bim{\Span}$ in an appropriate
sense, and that a category enriched in $V$ gives rise to a category enriched
in \Bim{V}. So far, $V$ has just been a bicategory, but the \fcat{Bim}
construction works for any \fc-multicategory $V$: indeed, this appears to be
its natural setting.

\item
Let $V$ be the bicategory $\Span\Cartpr$, where $T$ is a cartesian monad on a
cartesian category \Eee. A category enriched in $\Span\Cartpr$ consists of
	\begin{itemize}
	\item a set $C_0$
	\item for each $a\in C_0$, an object \eend{C}{a} of \Eee
	\item for each $a,b\in C_0$, a span
	\[
	\spn{\ehom{C}{a}{b}}{T\eend{C}{a}}{\eend{C}{b}}
	\]
	\item for each $a,b,c\in C_0$, a morphism
	\[
	\comp: \ehom{C}{b}{c}\of\ehom{C}{a}{b} \go \ehom{C}{a}{c}
	\]
	in \Eee
	\item for each $a\in C_0$, a morphism
	\[
	\ids: \eend{C}{a} \go \ehom{C}{a}{a}
	\]
	in \Eee,
	\end{itemize}
such that \comp\ and \ids\ are graph maps and satisfy associativity and
identity axioms.

\item
Let $V$ be the bicategory $\Span\pr{\Set}{\mr{free\ monoid}}$. We describe a
$V$-enriched category associated to any monoid.

Fix a monoid $M$. Denote by $\fcat{Sub}(M)$ the set of all submonoids of
$M$. Let $D$ be the category whose objects are pairs \pr{N}{X} with $N\in
\fcat{Sub}(M)$ and $X$ a left $N$-set, and in which an arrow $\pr{N}{X} \go
\pr{N'}{X'}$ consists of the inclusion $N\sub N'$ and a morphism $X\go X'$ of
$N$-sets. Note that if \range{\pr{N}{X_1}}{\pr{N}{X_k}} are all objects of
$D$ then we obtain an object \pr{N}{X_1 \times \cdots \times X_k} of $D$:
for $X_1 \times \cdots \times X_k$ is naturally a left $(N\times \cdots
\times N)$-set, and therefore becomes an $N$-set via the diagonal map $N\go
N\times \cdots \times N$. In other words, we put
\[
n.\bftuple{x_1}{x_k} = \bftuple{n.x_1}{n.x_k}.
\]

We have now described enough structure to make a $V$-enriched category. For
let $C_0 = \fcat{Sub}(M)$ and let 
%
%
\eend{C}{N} be the collection of left $N$-sets (or a small subcollection, if
one prefers); to define the span
\[
\spn{\ehom{C}{N}{N'}}{\eend{C}{N}^*}{\eend{C}{N'}},
\]
put
\[
\ehom{C}{N}{N'}\pr{\bftuple{X_1}{X_k}}{X'} = 
\homset{D}{\pr{N}{X_1 \times \cdots \times X_k}}{\pr{N'}{X'}}.
\]
Here \ust\ indicates free monoid, and $X_1 \times \cdots \times X_k$ is an
$N$-set as defined above. To describe \ehom{C}{N}{N'} in another way: if $N
\not\subseteq N'$ then $\ehom{C}{N}{N'} = \emptyset$, and if $N\sub N'$ then
an element of \ehom{C}{N}{N'} consists of a sequence \range{X_1}{X_k} of
$N$-sets ($k\geq 0$), an $N'$-set $X'$, and a function $X_1 \times \cdots
\times X_k \goby{f} X'$ such that
\[
f\bftuple{n.x_1}{n.x_k} = n.f\bftuple{x_1}{x_k}.
\]
Composition and identities work in an evident way.

\item	\label{eg:rel}

A category enriched in \fcat{Rel} (see~\ref{sec:two-dims}\bref{sec:rel})
consists of a set $I$, a preordered set $D$, and a map $\pi: D\go I$ of sets,
by just the same reasoning as we used for enrichment in \Span\
(part~\bref{eg:span}).

\item	\label{eg:glue}
A category enriched in \fcat{Glue} (see~\ref{sec:two-dims}\bref{sec:glue})
consists of $I$, $D$ and $\pi$ as in~\bref{eg:rel}, with the property that
for any $a\in D$ and $j\in I$, the sets
\begin{center}
$\{ b\in D \such b\leq a \mr{\ and\ } \pi(b)=j \}$,	\\
$\{ b\in D \such b\geq a \mr{\ and\ } \pi(b)=j \}$
\end{center}
each have at most one element.

\item	\label{eg:indexed-fam}
A set with an indexed family of subsets gives a category enriched in
\fcat{Glue}. For let $(C_i)_{i\in I}$ be a family of subsets of a set
$S$. Put $C_0 = I$, $\eend{C}{i} = C_i$, $\ehom{C}{i}{j} = C_i \cap C_j$, 
and take 
\begin{diagram}
\eend{C}{i}	&\lMonic	&\ehom{C}{i}{j}	
	&\rMonic	&\eend{C}{j}	\\
\end{diagram}
to be the obvious inclusions. Then 
\[
\ehom{C}{j}{k} \of \ehom{C}{i}{j} = C_i \cap C_j \cap C_k
\sub C_i \cap C_k = \ehom{C}{i}{k}
\]
and
\[
\eend{C}{i}=C_i=C_i \cap C_i=\ehom{C}{i}{i}
\]
so composition and identities can be defined (in a unique way), and
the axioms are satisfied. Or to put it in the way worked out
in~\bref{eg:glue}, the preordered set $D$ is the set $\coprod_{i\in I} C_i$
whose preorder is the equivalence relation induced by the function
$\coprod_{i\in I} C_i \go S$, and $\coprod_{i\in I} C_i \go I$ is
projection. 

\item
A more complicated version of the previous example: take a family
$(C_i)_{i\in I}$ of subsets of a set $S$ as before, but this time suppose $I$
has a preorder on it. Put $\eend{C}{i} = C_i$, and 
\[
\ehom{C}{i}{j} = 
\left\{
\begin{array}{ll}
C_i \cap C_j	&\mr{if\ }i\leq j	\\
\emptyset	&\mr{otherwise.}
\end{array}
\right.
\]
The axioms are still satisfied.  For the alternative formulation, the
preordered set $D$ is $\coprod_{i\in I} C_i$, with $\pr{i}{a} \leq
\pr{i'}{a'}$ if and only if $i\leq i'$ in $I$ and $a=a'$ as elements of
$S$. Thus the previous example corresponded to $I$ having the indiscrete
preorder.

\item
This example is in a different vein. Let $A$ be a nonempty
path-connected space, so that we can choose a basepoint $a_0$ and a path
$\gamma_a$ from $a_0$ to each point $a$, and let $V$ be the homotopy
bicategory of $A$ (see~\ref{sec:two-dims}\bref{sec:hty-bi}). Then we
obtain a category $C$ enriched in $V$:
\begin{itemize}
\item $C_0 = A$
\item $\eend{C}{a} = a$ for $a\in A$
\item $\ehom{C}{a}{b} = \gamma_b \of \gamma_a^*$, a path from $a$ to $b$,
where $\gamma_a^*$ is $\gamma_a$ run backwards
\item the composition 2-cell for \triple{a}{b}{c} is a homotopy from
\begin{figure}
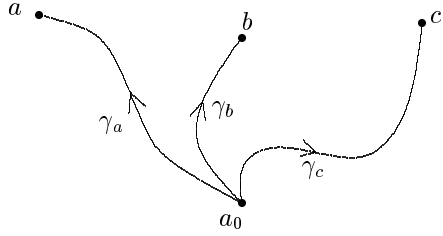

\centerline{\piccy{paths}}
\caption{Paths in $A$}
\label{fig:hty-pic}
\end{figure}
$\gamma_c \of \gamma_b^* \of \gamma_b \of \gamma_a^*$ to $\gamma_c \of
\gamma_a^*$, taken up to homotopy, and this comes from the obvious homotopy
from $\gamma_b^* \of \gamma_b$ to the constant path at $a_0$
(Figure~\ref{fig:hty-pic}).
\item the identity 2-cell $1_{\eend{C}{a}} \go \ehom{C}{a}{a}$ is the obvious
homotopy from the constant path at $a$ to $\gamma_a \of \gamma_a^*$ (taken up
to homotopy).
\end{itemize}
It is straightforward to check the associativity and identity axioms.

\item 
Finally, we give three related examples of a category enriched in a
commutative monoid (see~\ref{egs:enr-deg}\bref{eg:comm-mon}). For simplicity,
we stick to the case $r=0$. 

For the first example, let \fcat{Line} be the category of 1-dimensional real
vector spaces (or a small full subcategory, if one prefers). For every pair
\pr{l}{l'} of objects of \fcat{Line}, choose an isomorphism $\alpha_{ll'}: l
\go l'$. Let $A$ be the collection of objects of \fcat{Line}. Then for each
$l,l',l''\in A$, there is a unique real number \etrip{l}{l'}{l''} such that
$\alpha_{l'l''}\of\alpha_{ll'} = \etrip{l}{l'}{l''}.\alpha_{ll''}$, and for
each $l\in A$ there is a unique real number \esing{l} such that
$\esing{l}.\alpha_{ll}=I$. If we take $V$ to be the commutative monoid of
real numbers under multiplication then the axioms are satisfied.

For the second example, let $A$ be any set, let $V$ be any commutative
monoid, and let $\gamma: A\times A \go V$ be any function such that
$\gamma\pr{a}{b}$ is invertible for all $a,b\in A$. If we put
\begin{center}
$\etrip{a}{b}{c} = \gamma\pr{a}{b} + \gamma\pr{b}{c} - \gamma\pr{a}{c}$ \\
$\esing{a} = -\gamma\pr{a}{a}$ \\
\end{center}
then the axioms are satisfied and we obtain a category enriched in $V$.

For the third example, let $A$ be any subset of the plane. For each $a,b\in
A$ choose a smooth path $\gamma_{ab}$ from $a$ to $b$. Let \etrip{a}{b}{c} be
the area bounded by the loop formed by $\gamma_{ab}$, $\gamma_{bc}$, and
$\gamma_{ac}^*$ (i.e.\ $\gamma_{ac}$ run backwards). (Here, `area' might be
positive or negative: we choose an orientation on the plane so that, say, the
area bounded by an anticlockwise loop is positive.) Let $\esing{a} =
-(\mr{area\ bounded\ by\ }\gamma_{aa})$. Then we obtain a category enriched
in the abelian group of real numbers under addition. The axioms assert
obvious facts about area.

\end{eg}

We have now seen what enriched categories look like in elementary terms, and
looked at some examples. Next, we turn to the morphisms of enriched
categories: enriched functors.

\section{Enriched Functors: Elementary Description}

Let $V$ be an \fc-multicategory, and let $C$ and $D$ be categories enriched
in $V$. Recall from Definition~\ref{defn:enr-ftr} that a $V$-enriched functor
from $C$ to $D$ consists of a function $F_0: C_0 \go D_0$ together with a
transformation 
\begin{diagram}[width=2em,height=1em]
MI(C_0)	&	&\rTo^{MI(F_0)}	&	&MI(D_0)	\\
	&\rdTo(2,5)&		&\ \ldTo(2,5)&		\\
	&	&\ruNT>{F}	&	&		\\
	&\	&		&	&		\\
	&	&		&	&		\\
	&	&V		&	&		\\
\end{diagram}
of \fc-multicategories. Explicitly, this means that an enriched functor $C\go
D$ consists of:
\begin{itemize}
\item a function $F_0: C_0 \go D_0$
\item for each $a\in C_0$, a vertical 1-cell
\vslob{\eend{C}{a}}{F_a}{\eend{D}{F_0 a}} in $V$
\item for each $a,b\in C_0$, a 2-cell
\begin{diagram}
\eend{C}{a}	&\rTo^{\ehom{C}{a}{b}}		&\eend{C}{b}	\\
\dTo<{F_a}	&\Downarrow\,F_{ab}		&\dTo>{F_b}	\\
\eend{D}{F_0 a}	&\rTo_{\ehom{D}{F_0 a}{F_0 b}}	&\eend{D}{F_0 b}\\
\end{diagram}
in $V$, 
\end{itemize}
such that
\begin{eqnarray*}
\ &
\begin{diagram}[size=2em]
\eend{C}{a}	&\rTo^{\ehom{C}{a}{b}}		&\eend{C}{b}
&\rTo^{\ehom{C}{b}{c}}		&\eend{C}{c}	\\
\dTo<{F_a}	&\Downarrow F_{ab}		&\dTo<{F_b}	
&\Downarrow F_{bc}		&\dTo>{F_c}	\\
\eend{D}{F_0 a}	&\rTo_{\ehom{D}{F_0 a}{F_0 b}}	&\eend{D}{F_0 b}	
&\rTo_{\ehom{D}{F_0 b}{F_0 c}}	&\eend{D}{F_0 c}\\
\dTo<{1}	&				&\Downarrow \kappa^D_{F_0 a,
F_0 b, F_0 c}
&				&\dTo>{1}	\\
\eend{D}{F_0 a}	&				&\rTo_{\ehom{D}{F_0 a}{F_0 c}}
&				&\eend{D}{F_0 c}\\
\end{diagram}
\\
=&
\begin{diagram}[size=2em]
\eend{C}{a}	&\rTo^{\ehom{C}{a}{b}}		&\eend{C}{b}
&\rTo^{\ehom{C}{b}{c}}		&\eend{C}{c}	\\
\dTo<{1}	&				&\Downarrow \kappa^C_{abc}
&				&\dTo>{1}	\\
\eend{C}{a}	&				&\rTo_{\ehom{C}{a}{c}}	
&				&\eend{C}{c}	\\
\dTo<{F_a}	&			&\Downarrow F_{ac}	
&				&\dTo>{F_c}	\\
\eend{D}{F_0 a}	&				&\rTo_{\ehom{D}{F_0 a}{F_0 c}}
&				&\eend{D}{F_0 c}\\
\end{diagram}
\end{eqnarray*}
and
\[
\begin{diagram}[size=2em]
\eend{C}{a}	&\rEquals			&\eend{C}{a}	\\
\dTo<{F_a}	&=				&\dTo>{F_a}	\\
\eend{D}{F_0 a}	&\rEquals			&\eend{D}{F_0 a}\\
\dTo<{1}	&\Downarrow\kappa^D_{F_0 a}	&\dTo>{1}	\\
\eend{D}{F_0 a}	&\rTo_{\ehom{D}{F_0 a}{F_0 a}}	&\eend{D}{F_0 a}\\
\end{diagram}
\diagspace = \diagspace
\begin{diagram}[size=2em]
\eend{C}{a}	&\rEquals			&\eend{C}{a}	\\
\dTo<{1}	&\Downarrow\kappa^C_a		&\dTo>{1}	\\
\eend{C}{a}	&\rTo_{\ehom{C}{a}{a}}		&\eend{C}{a}	\\
\dTo<{F_a}	&\Downarrow F_{aa}		&\dTo>{F_a}	\\
\eend{D}{F_0 a}	&\rTo_{\ehom{D}{F_0 a}{F_0 a}}	&\eend{D}{F_0 a},\\
\end{diagram}
\]
where $\kappa^C$ indicates composition or identities in $C$, and
$\kappa^D$ in $D$.

Let us now look at what this means in some degenerate cases.

\begin{eg}{egs:deg-ftr}

\item \label{eg:bicat-enr-ftrs} 
Let $V$ be a bicategory and let $C$ and $D$ be $V$-enriched categories. A
$V$-enriched functor $F: C\go D$ consists of a function $F_0: C_0 \go D_0$,
satisfying $\eend{C}{a} = \eend{D}{F_0 a}$ for all $a\in C_0$ (since $V$ is
vertically discrete as an \fc-multicategory), with a 2-cell
\[
\eend{C}{a} = \eend{D}{F_0 a} 
\ctwo{\ehom{C}{a}{b}}{\ehom{D}{F_0 a}{F_0 b}}{F_{ab}}
\eend{C}{b} = \eend{D}{F_0 b}
\]
for each $a,b\in C_0$. 
These 2-cells are to be compatible with composition and identities in $C$ and
$D$. 

The requirement that $\eend{C}{a} = \eend{D}{F_0 a}$ seems rather unnatural,
and is too stringent in some practical cases.  This is one of the reasons why
we add in some vertical 1-cells to some of the bicategories we have been
considering: see~\ref{sec:specific-ftrs} below.

\item
Let $V$ be a monoidal category. Using~\bref{eg:bicat-enr-ftrs}, we see that a
functor $C\go D$ of $V$-enriched categories consists of a function $F_0: C_0
\go D_0$ and for each $a,b\in C_0$ a morphism $F_{ab}: \ehom{C}{a}{b} \go
\ehom{D}{F_0 a}{F_0 b}$ in $V$, such that
\begin{diagram}
\ehom{C}{b}{c} \otimes \ehom{C}{a}{b}	&\rTo^{F_{bc}\otimes F_{ab}}	
	&\ehom{D}{F_0 b}{F_0 c} \otimes \ehom{D}{F_0 a}{F_0 b}	\\
\dTo<{\kappa^C}		&	&\dTo>{\kappa^D}	\\
\ehom{C}{a}{c}	&\rTo_{F_{ac}}	&\ehom{D}{F_0 a}{F_0 c}	\\
\end{diagram}
and
\begin{diagram}
I		&\rTo^{1}	&I	\\
\dTo<{\kappa^C}	&		&\dTo>{\kappa^D}	\\
\ehom{C}{a}{a}	&\rTo_{F_{aa}}	&\ehom{D}{F_0 a}{F_0 a}	\\
\end{diagram}
commute. This is the traditional definition of a $V$-enriched functor.

\item
When $V$ is a plain multicategory, the description of a $V$-enriched functor
is just as when $V$ is a monoidal category, making the obvious
translation. 

\item
Let $V$ be an operad, and let $A$ and $A'$ be categories enriched in
$V$. Then an enriched functor $A\go A'$ consists of a function $f: A\go A'$
and for each $a,b\in A$ an element $f_{ab}$ of $V(1)$, satisfying the
equations
\begin{center}
$\etrip{fa}{fb}{fc}' \of (f_{ab}, f_{bc}) = f_{ac} \of \etrip{a}{b}{c}$	\\
$\esing{fa}' 	=	f_{aa} \of \esing{a}$.	
\end{center}

\item
Let $V$ be a commutative monoid with a specified invertible element $r$. Then
a map $A\go A'$ of categories enriched in \pr{V}{r} consists of functions $f:
A\go A'$ and $f_{\blob\blob}: A\times A \go V$, satisfying the equations
\begin{center}
$\etrip{fa}{fb}{fc}' + f_{ab} + f_{bc} = \etrip{a}{b}{c} + f_{ac}$	\\
$\esing{fa}' 	=	\esing{a} + f_{aa}$.	
\end{center}

\end{eg}

\section{Enriched Functors: Examples}	\label{sec:specific-ftrs}

We now give some more specific examples of enriched functors.

\begin{eg}{egs:specific-ftrs}

\item
Let $C$ and $C'$ be categories enriched in a bicategory $V$ with $C_0=C'_0=1$,
so that $C$ and $C'$ are just monads in $V$. Write them as
\tuplebts{X,t,\eta,\mu} and \tuplebts{X',t',\eta',\mu'} respectively. Then a
$V$-enriched functor $C\go C'$ is a `strict map of monads': that is, it's the
requirement that $X=X'$, together with a 2-cell
$
X \ctwo{t}{t'}{\phi} X
$
such that $\mu'\of (\phi *\phi) = \phi \of \mu$ and $\eta' = \phi \of
\eta$. 

\item
Given a 2-category $W$ we get a vertically discrete \fc-multicategory, just
as we have done for bicategories throughout. But we also get two more
\fc-multicategories, which are not in general vertically discrete, from the
two double categories $V$ and $V'$ of~\ref{sec:two-dims}\bref{dbl-cats}.

Note that all three of the \fc-multicategories we have constructed from
the 2-category $W$ are the same when we ignore the nonidentity vertical
1-cells, so a category enriched in $V$ or in $V'$ is just the same as in 
Example~\ref{egs:enr-deg}\bref{eg:enr-bi}. Let $C$ and $D$ be two categories
enriched in $W$. A $V$-enriched functor $C\go D$ is a function $F_0: C_0 \go
D_0$, plus for each $a\in C_0$ a 1-cell $F_a: \eend{C}{a} \go \eend{D}{F_0
a}$ in $W$, plus for each $a,b\in C_0$ a 2-cell
\begin{ntdiag}
\eend{C}{a}	&	&\rTo^{\ehom{C}{a}{b}}	&	&\eend{C}{b}	\\
		&	&			&\ 	&		\\
\dTo<{F_a}	&	&\ruNT>{\phi_{ab}}	&	&\dTo>{F_b}	\\
		&\ 	&			&	&		\\
\eend{D}{F_0 a}	&	&\rTo_{\ehom{D}{F_0 a}{F_0 b}}&	&\eend{D}{F_0 b}\\
\end{ntdiag}
in $W$, satisfying axioms stating compatibility with composition and
identities. In particular, suppose $C_0=D_0=1$, so that $C$ and $D$ are
monads in $W$. Then a $V$-enriched functor $C\go D$ is a monad functor (see
\cite{StFTM}), and a $V'$-enriched functor $C\go D$ is a monad opfunctor.

\item	\label{eg:new-span}
Previously, \Span\ had denoted a bicategory; we now redefine \Span\ to be a
certain weak double category. Since the underlying bicategory of this new
\Span\ is the old \Span, the previous explanation of categories enriched in
\Span\ remains valid.

The weak double category \Span\ is defined as follows. A 0-cell is a set, a
horizontal 1-cell is a span \spn{X}{A}{B}, a vertical 1-cell is a function
$A\go A'$ between sets, and a 2-cell is a function $h$ making the diagram
\begin{slopeydiag}
	&	&X	&	&	\\
	&\ldTo	&	&\rdTo	&	\\
A	&	&\dTo>{h}&	&B	\\
\dTo	&	&X'	&	&\dTo	\\
	&\ldTo	&	&\rdTo	&	\\
A'	&	&	&	&B'	\\
\end{slopeydiag}
commute. Horizontal composition of 1-cells is by pullback (as in the old
\Span), vertical composition of 1-cells is ordinary composition of functions,
and 2-cell composition works in the natural way. From now on, \Span\ will
denote this weak double category or the underlying \fc-multicategory (which
is what we really care about).

We saw in~\ref{egs:specific-enr-cats}\bref{eg:span} that a category
enriched in \Span\ consisted of a set $I$, for each $i\in I$ a set
\eend{D}{i}, for each $i,j\in I$ a span
\[
\spn{\ehom{D}{i}{j}}{\eend{D}{i}}{\eend{D}{j}},
\]
and functions for composition and identities. Put another way, it consists of
a category $D$, a set $I$, and a function $D_0 \go I$. If
\triple{D'}{I'}{D'_0\go I'} is another category enriched in \Span, then a
\Span-enriched functor between them consists of: 
\begin{itemize}
\item a function $f: I\go I'$
\item for each $i\in I$, a function \vslob{\eend{D}{i}}{F_i}{\eend{D'}{f(i)}}
\item for each $i,j\in I$ a function $F_{ij}$ making
\begin{slopeydiag}
	&	&\ehom{D}{i}{j}	&	&		\\
	&\ldTo	&		&\rdTo	&		\\
\eend{D}{i}&	&\dTo>{F_{ij}}	&	&\eend{D}{j}	\\
\dTo<{F_i}&	&\ehom{D'}{f(i)}{f(j)}&	&\dTo>{F_j}	\\
	&\ldTo	&		&\rdTo	&		\\
\eend{D'}{f(i)}&&		&	&\eend{D'}{f(j)}\\
\end{slopeydiag}
commute,
\end{itemize}
all compatible with the composition and identity functions. In other words, a
\Span-enriched functor from \triple{D}{I}{D_0\go I} to
\triple{D'}{I'}{D'_0\go I'} consists of a function $I\goby{f}I'$ and a
functor $D\goby{F}D'$ such that 
\begin{diagram}
D_0	&\rTo^{F_0}	&D'_0	\\
\dTo	&		&\dTo	\\
I	&\rTo_{f}	&I'	\\
\end{diagram}
commutes. Thus the category of \Span-enriched categories is the comma
category \commacat{\mr{ob}}{\Set}, where $\mr{ob}:\Cat\go\Set$ is the objects
functor.

\item
Similarly, \fcat{Rel} becomes a (strict) double category. There is 
at most one 2-cell with any given boundary; a 2-cell looks like
\begin{diagram}
X	&\rMod^{R}	&X'	\\
\dTo<{f}&		&\dTo>{f'}\\
Y	&\rMod_{S}	&Y',	\\
\end{diagram}
where $X, X', Y, Y'$ are sets, $R\sub X\times X'$ and $S\sub Y\times Y'$, and
$f$ and $f'$ are functions, all satisfying the condition that if
$\pr{x}{x'}\in R$ then $\pr{f(x)}{f'(x')}\in S$. Arguing as for \Span, the
category of \fcat{Rel}-enriched categories is the comma category
\commacat{U}{\Set}, where $U:\fcat{Preorders}\go\Set$ is the forgetful
functor.

\item
Since \fcat{Glue} is a sub-2-category of \fcat{Rel}, it also becomes a double
category. Suppose $(C_i)_{i\in I}$ is an indexed family of subsets of a set
$S$, and similarly $(C'_j)_{j\in I'}$ in $S'$. Suppose also that $f:I\go I'$
and $\phi: S\go S'$ are functions such that $\phi C_i \sub C'_{f(i)}$. Then
we get an enriched functor between the \fcat{Glue}-enriched categories
corresponding to these families
(see~\ref{egs:specific-enr-cats}\bref{eg:indexed-fam}), given by the
commuting square
\begin{diagram}
\coprod_{i}C_i	&\rTo		&\coprod_{j}C'_j	\\
\dTo		&		&\dTo		\\
I		&\rTo_{f}	&I'		\\
\end{diagram}
in which the function along the top row sends $a\in C_i$ to $\phi(a)\in
C'_{f(i)}$. 

\end{eg}

\section{Bimodules}	\label{sec:bim}

The bimodules construction has traditionally taken place in the context of
bicategories: given a bicategory \Bee, there is another bicategory \Bim{\Bee}
whose 0-cells are monads in \Bee\ and whose 1-cells are bimodules in \Bee\
(\cite{CKW}, \cite{Kos}). However, in order to do this one needs to assume
some special properties of \Bee, e.g.\ that \Bee\ locally has reflexive
coequalizers and that these are preserved by composition with any 1-cell. In
this section we describe a bimodules construction for \fc-multicategories
which extends the construction for bicategories and has the advantage that it
has no such technical restrictions on it. Having done this, we look at the
way in which a category enriched in an \fc-multicategory $V$ gives rise to a
category enriched in \Bim{V}.

Let $V$ be an \fc-multicategory. The \fc-multicategory \Bim{V} is defined as
follows:
\begin{description}
\item[0-cells]
A 0-cell of \Bim{V} is a \emph{monad} in $V$. That is, it is a 0-cell $x$ of
$V$ together with a horizontal 1-cell $x\goby{t}x$ and 2-cells
\[
\begin{diagram}
x	&\rTo^{t}	&x			&\rTo^{t}	&x	\\
\dTo<{1}&		&\Downarrow\,\mu	&		&\dTo>{1}\\
x	&		&\rTo_{t}		&		&x	\\
\end{diagram}
\diagspace
\begin{diagram}
x	&\rEquals		&x	\\
\dTo<{1}&\Downarrow\,\eta	&\dTo>{1}\\
x	&\rTo_{t}		&x	\\
\end{diagram}
\]
satisfying the usual monad axioms, $\mu\of\pr{\mu}{1_t} =
\mu\of\pr{1_t}{\mu}$ and $\mu\of\pr{\eta}{1_t} = 1 = \mu\of\pr{1_t}{\eta}$. 

\item[Horizontal 1-cells]
A horizontal 1-cell $(x,t,\eta,\mu) \rMod (x',t',\eta',\mu')$ consists of a
horizontal 1-cell $x\goby{f}x'$ in $V$ together with 2-cells
\[
\begin{diagram}
x	&\rTo^{t}	&x			&\rTo^{f}	&x'	\\
\dTo<{1}&		&\Downarrow\,\theta	&		&\dTo>{1}\\
x	&		&\rTo_{f}		&		&x'	\\
\end{diagram}
\diagspace
\begin{diagram}
x	&\rTo^{f}	&x'			&\rTo^{t'}	&x'	\\
\dTo<{1}&		&\Downarrow\,\theta'	&		&\dTo>{1}\\
x	&		&\rTo_{f}		&		&x'	\\
\end{diagram}
\]
satisfying the usual algebra axioms $\theta\of\pr{\eta}{1_f}=1$,
$\theta\of\pr{\mu}{1_f} = \theta\of\pr{1_t}{\theta}$, and dually for
$\theta'$, and the `commuting actions' axiom $\theta'\of\pr{\theta}{1_{t'}} =
\theta\of\pr{1_t}{\theta'}$.

\item[Vertical 1-cells]
A vertical 1-cell
\vslob{(x,t,\eta,\mu)}{}{(\hat{x},\hat{t},\hat{\eta},\hat{\mu})} 
in \Bim{V} is a vertical 1-cell \vslob{x}{p}{\hat{x}} in $V$ together with a
2-cell
\begin{diagram}
x		&\rTo^{t}		&x		\\
\dTo<{p}	&\Downarrow\,\omega	&\dTo>{p}	\\
\hat{x}		&\rTo_{\hat{t}}		&\hat{x}	\\
\end{diagram}
such that $\omega\of\mu = \hat{\mu}\of\pr{\omega}{\omega}$ and $\omega\of\eta
= \hat{\eta}\of p$. (The notation on the right-hand side of the second
equation is explained on page~\pageref{p:null-notation}.)

\item[2-cells]
A 2-cell 
%
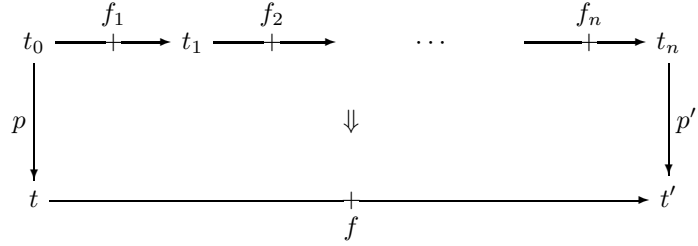
\begin{diagram}
t_0	&\rMod^{f_1}	&t_1	&\rMod^{f_2}	&\ 	&\cdots	
&\ 	&\rMod^{f_{n}}	&t_n	\\
\dTo<{p}&		&	&		&\Downarrow	&	
&	&		&\dTo>{p'}\\
t	&		&	&		&\rMod_{f}	&	
&	&		&t'	\\
\end{diagram}
%
in \Bim{V}, where $t$ stands for $(x,t,\eta,\mu)$, $f$ for
\triple{f}{\theta}{\theta'}, $p$ for \pr{p}{\omega}, and so on, consists of a
2-cell 
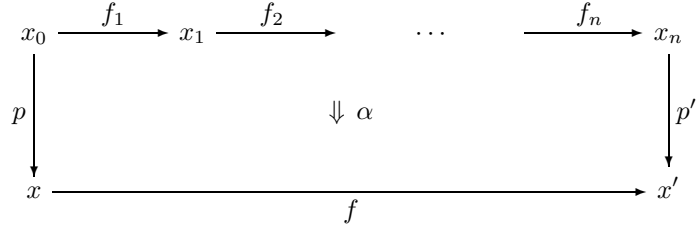
\begin{diagram}
x_0	&\rTo^{f_1}	&x_1	&\rTo^{f_2}	&\ 	&\cdots	
&\ 	&\rTo^{f_{n}}	&x_n	\\
\dTo<{p}&		&	&		&\Downarrow\,\alpha&	
&	&		&\dTo>{p'}\\
x	&		&	&		&\rTo_{f}	&	
&	&		&x'	\\
\end{diagram}
in $V$, satisfying the `external equivariance' axioms
%
%
\begin{eqnarray*}
\alpha\of(\theta_1,\range{1_{f_2}}{1_{f_n}}) 		&=&
\theta\of\pr{\omega}{\alpha}				\\
\alpha\of(\range{1_{f_1}}{1_{f_{n-1}}},\theta'_n)	&=&
\theta'\of\pr{\alpha}{\omega'}
\end{eqnarray*}
and the `internal equivariance' axioms
\[\hspace*{-12mm}
\alpha\of(\range{1_{f_1}}{1_{f_{i-2}}}, \theta'_{i-1}, 1_{f_{i}},
\range{1_{f_{i+1}}}{1_{f_n}})
=
\alpha\of(\range{1_{f_1}}{1_{f_{i-2}}}, 1_{f_{i-1}}, \theta_i,
\range{1_{f_{i+1}}}{1_{f_n}}) 
\]
for $2\leq i\leq n$.

\item[Composition and identities]
For both 2-cells and vertical 1-cells in \Bim{V}, composition is defined
directly from the composition in $V$, and identities similarly.

\end{description}

\begin{eg}{egs:bims}
\item
Let \Bee\ be a bicategory with the `special properties' mentioned in the
first paragraph of this section. Let $V$ be the \fc-multicategory coming from
\Bee. Then a 0-cell of \Bim{V} is a monad in \Bee, a horizontal 1-cell $t\go
t'$ is a \pr{t'}{t}-bimodule, and a 2-cell of the form
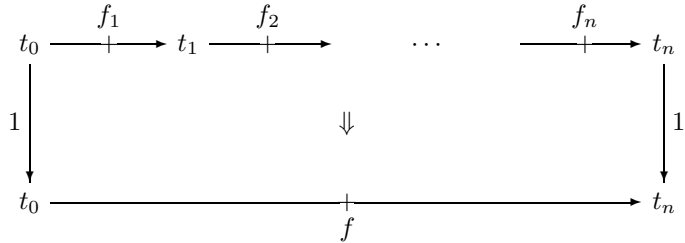
\begin{diagram}
t_0	&\rMod^{f_1}	&t_1	&\rMod^{f_2}	&\ 	&\cdots	
&\ 	&\rMod^{f_{n}}	&t_n	\\
\dTo<{1}&		&	&		&\Downarrow	&	
&	&		&\dTo>{1}\\
t_0	&		&	&		&\rMod_{f}	&	
&	&		&t_n	\\
\end{diagram}
is a map 
\[
f_n \otimes_{t_{n-1}}\cdots\otimes_{t_1} f_1 \go f
\] 
of \pr{t_n}{t_0}-bimodules, i.e.\ a 2-cell in the `traditional' bicategory
\Bim{\Bee}. Thus the \fc-multicategory coming from \Bim{\Bee} is the
vertically discrete part of \Bim{V}.

\item
Suppose $V$ comes from the monoidal category \pr{\Set}{\times}. Then a 0-cell
of \Bim{V} is a monoid, a horizontal 1-cell $M\rMod^{X} M'$ is a set $X$ with
commuting left action by $M'$ and right action by $M$, a vertical 1-cell
\vslob{M}{}{N} is a monoid homomorphism, and a 2-cell
\begin{diagram}
M_0	&\rMod^{X_1}	&M_1	&\rMod^{X_2}	&\ 	&\cdots	
&\ 	&\rMod^{X_{n}}	&M_n	\\
\dTo<{p}&		&	&		&\Downarrow	&	
&	&		&\dTo>{p'}\\
M	&		&	&		&\rMod_{X}	&	
&	&		&M'	\\
\end{diagram}
is a function $\phi: X_n \times\cdots\times X_1 \go X$ such that
\begin{eqnarray*}
\phi(\range{x_n}{x_2},x_1 m_0) 			&=& 
\phi\bftuple{x_n}{x_1}p(m_0)					\\
\phi(\ldots,x_i, m_{i-1} x_{i-1}, \ldots)	&=&
\phi(\ldots,x_i m_{i-1}, x_{i-1}, \ldots)\ \ \ (2\leq i\leq n)	\\
\phi(m_n x_n,\range{x_{n-1}}{x_1})		&=&
p'(m_n)\phi\bftuple{x_n}{x_1}.
\end{eqnarray*}
So \Bim{V} is the weak double category \Bim{\Set} of
\ref{egs:specific-enr-cats}\bref{eg:spec-bim-set}. 

\item
Similarly, by taking $V$ to be the monoidal category \pr{\Ab}{\otimes}, we
get the weak double category \Bim{\Ab} of
\ref{egs:degens}\bref{eg:weak-double}, made up of rings, bimodules, ring
homomorphisms, and bimodule maps under changes of base.

\item
Moving on from monoidal categories, let \Span\ be the weak double category of
\ref{egs:specific-ftrs}\bref{eg:new-span}, made up of sets, spans, functions,
and maps of spans. Then \Bim{\Span} is an \fc-multicategory with:
\begin{description}
\item[0-cells] Categories
\item[Horizontal 1-cells] Profunctors
\item[Vertical 1-cells] Functors
\item[2-cells] A 2-cell 
\begin{equation}	\label{eq:big-prof}
\begin{diagram}
C_0	&\rMod^{X_1}	&C_1	&\rMod^{X_2}	&\ 	&\cdots	
&\ 	&\rMod^{X_{n}}	&C_n	\\
\dTo<{F}&		&	&		&\Downarrow	&	
&	&		&\dTo>{F'}\\
C	&		&	&		&\rMod_{X}	&	
&	&		&C'	\\
\end{diagram}
\ \ \ 
\end{equation}
(where the $C$'s are categories, the $X$'s profunctors, and the $F$'s
functors) is a natural transformation
\begin{diagram}[width=2em,height=1em]
C_0^{\op}\times C_n&	&\rTo^{F^{\op}\times F'}&&C^{\op}\times C'\\
		&\rdTo(2,5)<{X_n \otimes\cdots\otimes X_1}
	&	&\ \ldTo(2,5)>{X}&	\\
		&	&\ruNT		&	&		\\
		&\	&		&	&		\\
		&	&		&	&		\\
		&	&\Set.		&	&		\\
\end{diagram}
In other words, it is a family of functions
\[
X_1\pr{c_0}{c_1}\times X_2\pr{c_1}{c_2} \times\cdots\times
X_n\pr{c_{n-1}}{c_n}
\go
X\pr{Fc_0}{F'c_n} 
\]
natural in $c_i\in C_i$.
\end{description}

\item	\label{eg:bim-cart-span}
In Example~\ref{egs:specific-ftrs}\bref{eg:new-span} we changed the
term \Span\ from meaning a certain bicategory to meaning a certain
\fc-multicategory. This process goes through without change for
$\Span\Cartpr$, so $\Span\Cartpr$ will now denote an \fc-multicategory, for
any cartesian \Cartpr. (It won't usually be a weak double category, though,
unless \Eee\ and $T$ have convenient properties.) Then \Bim{\Span\Cartpr}
has: 
\begin{description}
\item[0-cells] $T$-multicategories
\item[Horizontal 1-cell $C\go C'$:] span (\spn{X}{TC_0}{C'_0}), together with
maps (`actions') $X\of C_1 \go X$ and $C'_1\of X \go X$ satisfying some
axioms. We might call such an $X$ a \emph{profunctor} $C\rMod C'$.
\item[Vertical 1-cell $C\go D$:] functor $C\go D$ of $T$-multicategories
\item[2-cells] A 2-cell as at~\bref{eq:big-prof} is an arrow $X_n\of\cdots\of
X_1 \go X$ in \Eee\ (where the domain of this arrow is the 1-cell composite
in $\Span\Cartpr$), satisfying compatibility axioms for the actions by $C_i$,
$C$, and $C'$.
\end{description}

In the case of \Span, i.e.\ $\Cartpr=\Zeropr$, a 2-cell
\begin{diagram}
C	&\rMod^{1}	&C		\\
\dTo<{F}&\Downarrow	&\dTo>{G}	\\
D	&\rMod^{1}	&D		\\
\end{diagram}
is just a natural transformation $F\go G$ (by a Yoneda argument). The same is
in fact true for \Cartpr-multicategories: in \ref{defn:transf} we
defined a transformation $F\go G$ as an arrow $C_0 \goby{\alpha}D_1$ with
certain properties, but now we have an alternative description of it as an
arrow $C_1 \goby{\tilde{\alpha}} D_1$ with certain properties. In the case
$\Cartpr=\Zeropr$, $\tilde{\alpha}$ sends an arrow $f$ in $C$ to the diagonal
of the naturality square for $\alpha$ at $f$.

\end{eg}

To finish this section, we briefly discuss `change of base' and the fact that
from a category enriched in $V$ there arises a category enriched in
\Bim{V}. If $W\goby{G} W'$ is a functor between \Cartprp-multicategories then
any $W$-enriched multicategory becomes a $W'$-enriched multicategory just by
composition with $G$: thus there's a functor
\[
G_*: \Cartpr_{W}\hyph\Multicat \go \Cartpr_{W'}\hyph\Multicat.
\]

In particular, consider for any \fc-multicategory $V$ the forgetful functor
$U: \Bim{V} \go V$. This gives the functor
\[
U_*: \Cat_{\Bim{V}} \go \Cat_V,
\]
where $\Cat_W$ is the category of $W$-enriched categories. Then $U_*$ has a
right adjoint $Z$; this is not of the form $J_*$ for any $J$, but is the
functor $\Cat_V \go \Cat_{\Bim{V}}$ of which we saw examples
in~\ref{egs:specific-enr-cats}\bref{eg:spec-bim-ab}, \bref{eg:spec-bim-set},
\bref{eg:spec-bim-span}. On this occasion we omit the general definition of
$Z$, since it is indicated adequately by these examples and there is only one
sensible way of constructing a functor $\Cat_V \go \Cat_{\Bim{V}}$ for
general $V$.

\chapter{Other Enriched Multicategories}	\label{sec:other}

We have defined \Cartpr-multicategories enriched in an
\Cartprp-multicategory, but so far only looked at the case
$\Cartpr=\Zeropr$. In this chapter we look at some other cases. It should not
come as a surprise after Chapter~\ref{sec:cats} that \Cartprp-multicategories
are rather complicated structures for some of the usual examples of \Eee\ and
$T$; in the bulk of this chapter we therefore take $T$ to be the next most
simple case, the free monoid monad on $\Eee=\Set$. This does, however, give a
good idea of what $T_n$-multicategories enriched in a $T_{n+1}$-multicategory
look like for all $n$, and we discuss this (and the `full hierarchy') in the
final section.

A plain multicategory is a \pr{\Set}{\mr{free\ monoid}}-multicategory
(see~\ref{sec:pre-mti}). We may therefore speak of a plain multicategory
enriched in a \pr{\ust\hyph\Gph}{\fm}-multicategory, where \ust\ is the free
monoid monad on \Set, \ust-\Gph\ is the category of \ust-graphs
(\ref{sec:pre-mti}), and \fm\ is the free multicategory monad on
\ust-\Gph. In order to understand enriched plain multicategories, our first
task is therefore to see what an \fm-multicategory is.

\section{\fm-Multicategories}	\label{sec:fm-mti}

A \ust-graph $X$ consists of a set $X_0$ of objects $x$, together with a set
$X_1$ of arrows $\xi$ with specified domain and codomain, pictured as
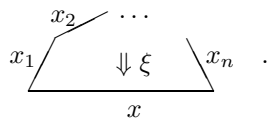
\begin{opetope}
	&	&	&\cnr 	&\ldots	&	&	&	\\
	&\cnr	&\ruLine(2,1)^{x_2}&&	&	&\cnr 	&	\\
\ruLine(1,2)<{x_1}&&	&	&\Downarrow \xi&&	&\rdLine(1,2)>{x_n\ \ \ .}\\
\cnr	&	&	&	&\rLine_{x}&	&	&\cnr	\\
\end{opetope}
In the graph $\fm(X)$, the set of objects is just $X_0$ again, but an arrow
is a formal gluing of arrows in $X$, such as
\begin{equation}	\label{eq:glued-twos}
\piccy{labpd}.
\end{equation}
We can define $\fm(X) = (\spaan{X'_1}{X_0^*}{X_0}{\dom}{\cod})$
inductively by:
\begin{itemize}
\item if $x\in X_0$ then $1_x\in X'_1$; $\dom(1_x) = (x)$ and $\cod(1_x) = x$
\item if $\xi\in X_1$ and $\theta_1, \ldots, \theta_n \in X'_1$ with
$\dom(\xi) = \bftuple{\cod(\theta_1)}{\cod(\theta_n)}$, then
$\xi\abftuple{\theta_1}{\theta_n}\in X'_1$; domain and codomain are given by
\begin{eqnarray*}
\dom(\xi\abftuple{\theta_1}{\theta_n})	&=	
&\mu_{X_0}\bftuple{\dom(\theta_1)}{\dom(\theta_n)},	\\
\cod(\xi\abftuple{\theta_1}{\theta_n})	&=	&\cod(\xi).
\end{eqnarray*}
\end{itemize}
Here $1_x$ and $\xi\abftuple{\theta_1}{\theta_n}$ are formal
expressions; the definition is just a special case of the free multicategory
construction described in the Appendix. Informally, it is clear how the
functor \fm\ acts on morphisms, and what its monad structure is, so we omit
the details here.

An \fm-graph looks like
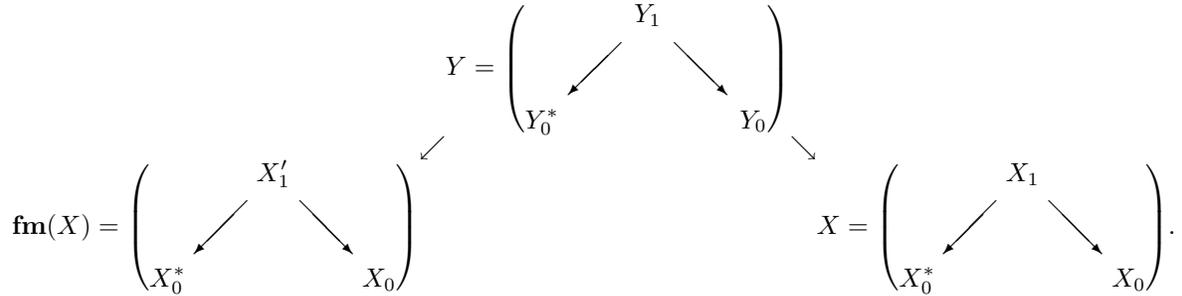
\begin{diagram}
	&	&Y=\pointyspn{Y_1}{Y_0^*}{Y_0}	&	&	\\
	&\swarrow	&		&\searrow	&	\\
\fm(X)=\pointyspn{X'_1}{X_0^*}{X_0}&&&	&X=\pointyspn{X_1}{X_0^*}{X_0}.	\\
\end{diagram}
Interpret the whole graph as a 3-dimensional structure this time: $x\in X_0$
is a 1-cell, $\xi\in X_1$ is a horizontal 2-cell, $y\in Y_0$ is a vertical
2-cell, and $\eta\in Y_1$ is a 3-cell, as in the picture
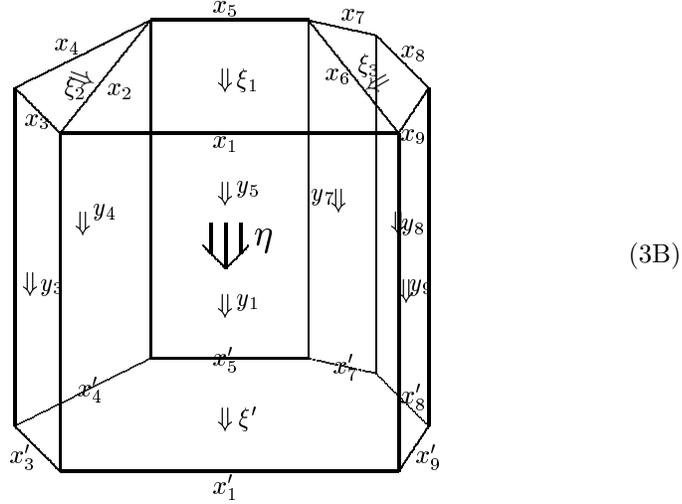
\begin{equation}	\label{eq:three-cell}
\piccy{3cell}
\end{equation}
or, drawn another way,
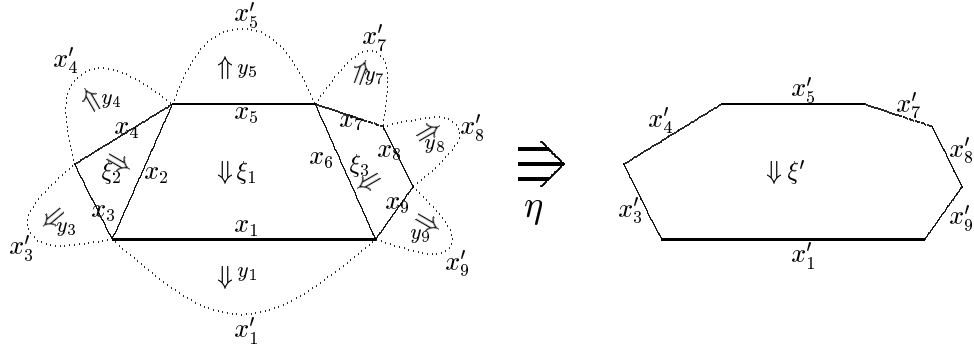
\[
\piccy{3cellnet}.
\]
An \fm-multicategory is an \fm-graph together with identities and
composition. This firstly means that \spn{Y_0}{X_0}{X_0} has the structure of
a category, i.e.\ that the 1-cells and vertical 2-cells are respectively the
objects and arrows of a category. More significantly, it means that 3-cells
can be composed vertically. The identities can be portrayed as
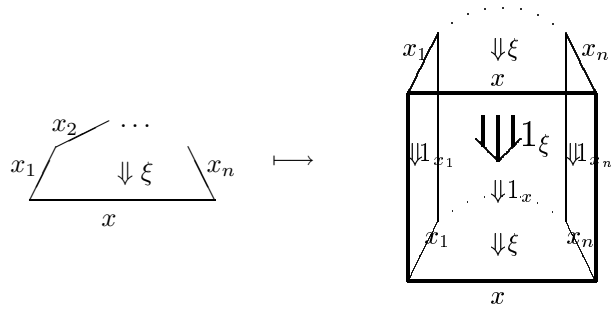
\[
\begin{opetope}
		&	&	&\cnr 	&\ldots	&	&	&	\\
		&\cnr	&\ruLine(2,1)^{x_2}&&	&	&\cnr 	&	\\
\ruLine(1,2)<{x_1}&	&	&	&\Downarrow \xi&&	
&\rdLine(1,2)>{x_n}\\
\cnr		&	&	&\rLine_{x}&	&	&	&\cnr\\
\end{opetope}
\diagspace\goesto\diagspace
\begin{array}{c}\piccy{fmmtiid}\end{array}
\]
and composition looks like
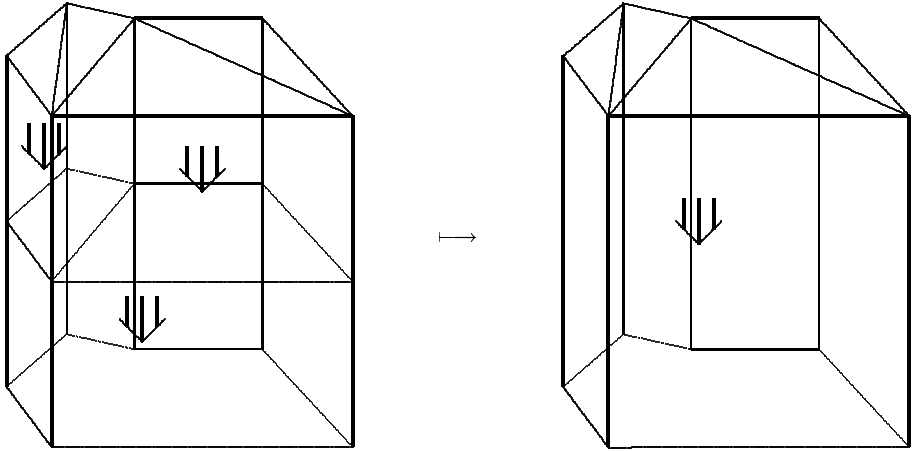
\[
\piccy{fmmticomp}.
\]

We describe a few degenerate cases, as we did for \fc-multicategories.

\begin{eg}{egs:degen-fm}

\item	\label{eg:bim-fm}
There are \fm-multicategories in which diagrams of horizontal 2-cells are
`representable', in a sense which we do not make precise but which is
analogous to the sense in which a string
$x_0\goby{\xi_1}\cdots\goby{\xi_n} x_n$ of horizontal 1-cells in a weak
double category is `represented' by $\xi_n\of\cdots\of\xi_1$. For
instance, there is an \fm-multicategory where a 1-cell is a ring, a
horizontal 2-cell
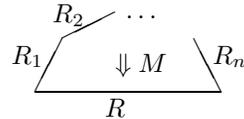
\begin{opetope}
		&	&	&\cnr 	&\ldots	&	&	&	\\
		&\cnr	&\ruLine(2,1)^{R_2}&&	&	&\cnr 	&	\\
\ruLine(1,2)<{R_1}&	&	&	&\Downarrow M&&	
&\rdLine(1,2)>{R_n}\\
\cnr		&	&	&\rLine_{R}&	&	&	&\cnr\\
\end{opetope}
is an abelian group $M$ with commuting left action by $R$ and
right actions by $R_1, \ldots, R_n$ (which we shall call an
\emph{$(R;R_1,\ldots,R_n)$-module}), a vertical 2-cell from $R$ to $R'$ is a
ring homomorphism, and a 3-cell inside
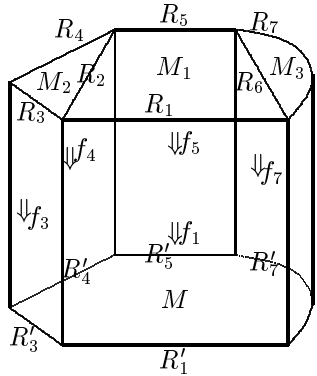
\[
\piccy{bim3cell}
\]
is a multilinear map $\phi: M_1, M_2, M_3 \go M'$ of abelian groups,
satisfying `internal compatibility' axioms
\begin{eqnarray*}
\phi\triple{m_1 r_2}{m_2}{m_3}	&=	&\phi\triple{m_1}{r_2 m_2}{m_3}	\\
\phi\triple{m_1 r_6}{m_2}{m_3}	&=	&\phi\triple{m_1}{m_2}{r_6 m_3}
\end{eqnarray*}
and `external compatibility' axioms
\begin{eqnarray*}
\phi\triple{r_1 m_1}{m_2}{m_3}	&=&f_1(r_1).\phi\triple{m_1}{m_2}{m_3}	\\
\phi\triple{m_1}{m_2 r_3}{m_3}	&=&\phi\triple{m_1}{m_2}{m_3}.f_3(r_3)	\\
\phi\triple{m_1}{m_2 r_4}{m_3}	&=&\phi\triple{m_1}{m_2}{m_3}.f_4(r_4)	\\
\phi\triple{m_1 r_5}{m_2}{m_3}	&=&\phi\triple{m_1}{m_2}{m_3}.f_5(r_5)	\\
\phi\triple{m_1}{m_2}{m_3 r_7}	&=&\phi\triple{m_1}{m_2}{m_3}.f_7(r_7).	\\
\end{eqnarray*}
There is an $(R_1;R_3,R_4,R_5,R_7)$-module $M$ with the property that 3-cells
$\phi$ as illustrated correspond naturally to maps $M\go M'$, where `maps'
means module maps with respect to the change of base
$(f_1;f_3,f_4,f_5,f_7)$. This $M$ is simply the quotient of the abelian group
tensor $M_1\otimes
M_2\otimes M_3$ by the relations
\begin{eqnarray*}
m_1 r_2 \otimes m_2 \otimes m_3	&\sim	&m_1 \otimes r_2 m_2 \otimes m_3\\
m_1 r_6 \otimes m_2 \otimes m_3	&\sim	&m_1 \otimes m_2 \otimes r_6 m_3
\end{eqnarray*}
(corresponding to the internal compatibility axioms above). The existence of
such an $M$ is what was meant by `representability'.

\item	\label{eg:fm-vd}
Suppose $(\spn{Y_0}{X_0}{X_0}) = (\spaan{X_0}{X_0}{X_0}{1}{1})$, so that all
vertical 2-cells are identities. Then the 3-cells are shaped like 3-opetopes
(for which terminology, see~\ref{sec:opetopes}): by collapsing the
vertical 2-cells (say of~\bref{eq:three-cell}), one obtains a 3-dimensional
figure with one flat face on the bottom and several (necessarily curved)
faces on the top. The graph of the \fm-multicategory can now be drawn as
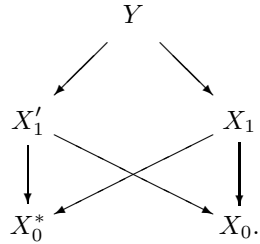
\begin{slopeydiag} 
	&	&Y	&	&	\\
	&\ldTo	&	&\rdTo	&	\\
X'_1	&	&	&	&X_1	\\
\dTo	&\rdTo(4,2)&	&\ldTo(4,2)&\dTo\\
X^*_0	&	&	&	&X_0.	\\
\end{slopeydiag} 
We call such an \fm-multicategory \emph{vertically discrete}. An example of
such a structure is a tricategory with just one 0-cell, also known as a
monoidal bicategory.

Whether or not an \fm-multicategory is vertically discrete, we will use the
convention that identity vertical 2-cells will not be drawn at all, as in
diagram~\bref{eq:no-walls}. 

\item	\label{eg:V-T2}
Suppose that the underlying vertical category is the terminal category, i.e.\
that $X_0=Y_0=1$. Then the \fm-multicategory consists of some objects $\xi$
(previously called horizontal 2-cells) of shape
\begin{equation}	\label{eq:two-ope}
\begin{opetope}
		&	&	&\cnr 	&\ldots	&	&	&	\\
		&\cnr	&\ruLine(2,1)&	&	&	&\cnr 	&	\\
\ruLine(1,2)	&	&	&	&	&	&	&\rdLine(1,2)\\
\cnr		&	&	&\rLine	&	&	&	&\cnr\\
\end{opetope}
\diagspace
\mbox{($n$ upper edges)}
\end{equation}
for each $n$, plus some arrows (previously called 3-cells) $\eta$, where the
domain of an arrow is a diagram of objects pasted together, and the codomain
is a single object with the same shape as the boundary of the domain, as in
\begin{equation}	\label{eq:arrow-in-two-mti}
\piccy{t23cell}.
\end{equation}
These arrows can be composed, and the composition obeys associative and
identity laws. In other words, it is a $T_2$-multicategory (see \cite[p.\
66]{SHDCT}).

\item	\label{eg:smc-fm}
Any symmetric monoidal category \pr{W}{\otimes} gives rise to a
$T_2$-multicategory $V$ in which the pastings of objects are
representable. Explicitly, let $V$ be the $T_2$-multicategory in which an
object of shape~\bref{eq:two-ope} is an object of $W$ (for any $n$), and in
which an arrow $\eta$ as in~\bref{eq:arrow-in-two-mti} is a morphism
\[
\xi_1 \otimes \xi_2 \otimes \xi_3 \otimes \xi_4 \otimes \xi_5
\go \xi
\]
in $W$. Composition and identities are defined in the natural way; since the
ordering of the $\xi_i$'s in the domain of~\bref{eq:arrow-in-two-mti} has no
special properties, we will need to use the symmetry isomorphisms in $W$ when
we are defining composites. Choosing different orderings only changes $V$ up
to isomorphism.

\end{eg}

\clearpage
\section{Enriched Plain Multicategories: Elementary Description}

We have now seen what an \fm-multicategory is. The
next question: if $A$ is a set, what is $MI(A)$? Firstly, $I(A)$ has graph
\spn{A^* \times A}{A^*}{A}. Then, $\fm(I(A)) = (\spn{B}{A^*}{A})$, where an
element of $B$ is a diagram like
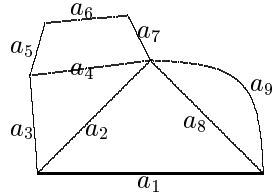
\[
\piccy{semilabpd1}
\]
($x_i\in A$; cf.\ diagram~\bref{eq:glued-twos}). There's a natural map $B\go
A^* \times A$ specifying domain in the first coordinate and codomain in the
second, and $M(I(A))$ has graph
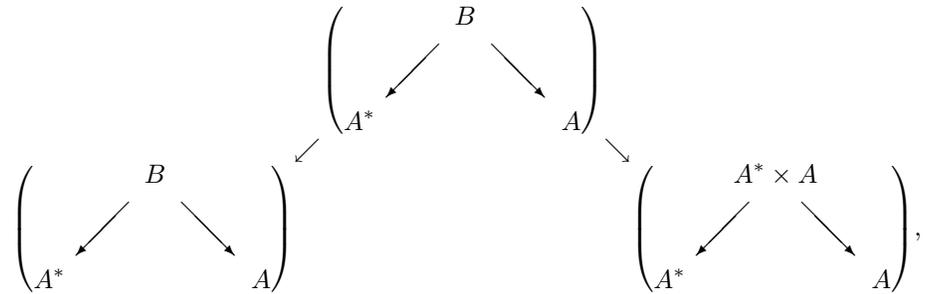
\begin{diagram}
	&	&\pointyspn{B}{A^*}{A}	&	&	\\
	&\swarrow	&		&\searrow	&	\\
\pointyspn{B}{A^*}{A}&&&	&\pointyspn{A^* \times A}{A^*}{A},	\\
\end{diagram}
which can also be written
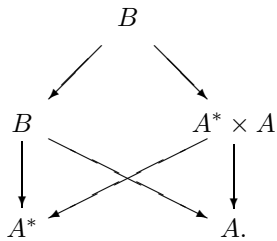
\begin{slopeydiag} 
	&	&B	&	&	\\
	&\ldTo	&	&\rdTo	&	\\
B	&	&	&	&A^*\times A	\\
\dTo	&\rdTo(4,2)&	&\ldTo(4,2)&\dTo\\
A^*	&	&	&	&A.	\\
\end{slopeydiag} 
Thus a horizontal 1-cell is an element of $A$; there's precisely one horizontal
2-cell
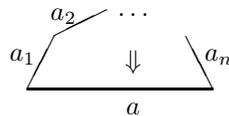
\[
\begin{opetope}
		&	&	&\cnr 	&\ldots	&	&	&	\\
		&\cnr	&\ruLine(2,1)^{a_2}&&	&	&\cnr 	&	\\
\ruLine(1,2)<{a_1}&	&	&	&\Downarrow&	&	
&\rdLine(1,2)>{a_n}\\
\cnr		&	&	&	&\rLine_{a}&	&	&\cnr	\\
\end{opetope}
\]
for each $a_1, \ldots, a_n, a \in A$; the only vertical 2-cells are
identities (i.e.\ $MI(A)$ is vertically discrete); and there's precisely one
3-cell
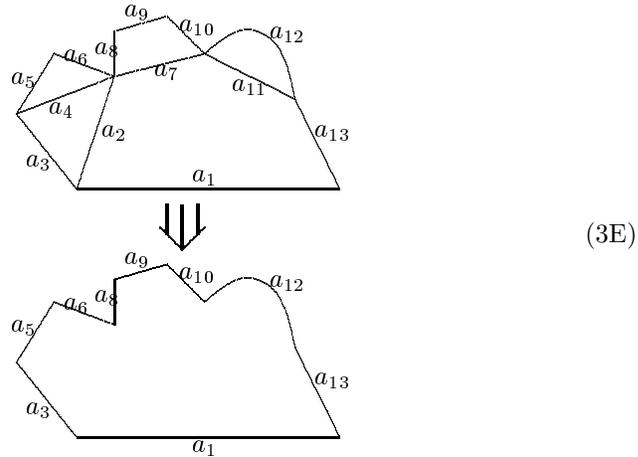
\begin{equation}	\label{eq:no-walls}
\piccy{mia3cell}
\end{equation}
(and similarly for any other pasting diagram of objects in the domain), using
the diagrammatic convention of \ref{egs:degen-fm}\bref{eg:fm-vd}. There is
only one possible way that composition and identities can be defined.

Finally, we can say what a plain multicategory enriched in a
\pr{\ust\hyph\Gph}{\fm}-multicategory $V$ is: namely,
\begin{itemize}
\item a set $C_0$
\item for each $a\in C_0$, a 1-cell \eend{C}{a} of $V$
\item for each $\range{a_1}{a_n}, a \in C_0$, a horizontal 2-cell
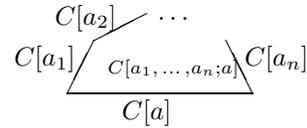
\begin{opetope}
		&	&	&\cnr 	&\ldots	&	&	&	\\
		&\cnr	&\ruLine(2,1)^{\eend{C}{a_2}}&&	&	&\cnr 	
&	\\
\ruLine(1,2)<{\eend{C}{a_1}}&	&	&	
&{\scriptstyle C[\range{a_1}{a_n};a]}&	&	&\rdLine(1,2)>{\eend{C}{a_n}}\\
\cnr		&	&	&\rLine_{\eend{C}{a}}&	&	&	&\cnr\\
\end{opetope}
of $V$
\item for each diagram like
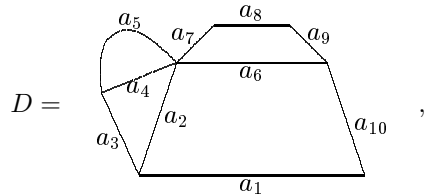
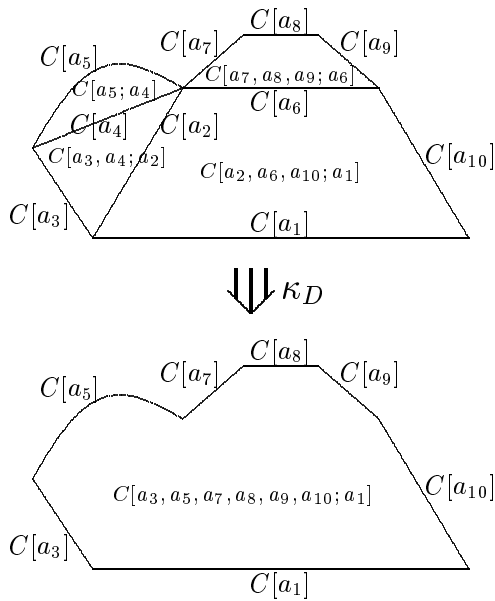
\[
D = \begin{array}{c}\piccy{semilabpd2}\end{array},
\]
a 3-cell like
\[
\piccy{comp3cell}
\]
\end{itemize}
such that the $\kappa_D$'s are closed under composition and identities, in a
way similar to that for enriched categories
(page~\pageref{p:kappas-closed}). 

We have become vague when discussing diagrams `like $D$', not wanting to get
too involved in the details of the explicit definition of \fm. By induction
we could prove that the definition of $V$-enriched plain multicategory is
equivalent to the following `biased' version (cf.\
page~\pageref{p:finitary}): a plain multicategory enriched in $V$ consists of
\begin{itemize}
\item a set $C_0$
\item for each $a\in C_0$, a 1-cell \eend{C}{a} of $V$
\item for each $\range{a_1}{a_n}, a \in C_0$, a horizontal 2-cell
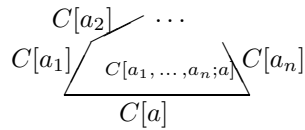
\begin{opetope}
		&	&	&\cnr 	&\ldots	&	&	&	\\
		&\cnr	&\ruLine(2,1)^{\eend{C}{a_2}}&&	&	&\cnr 	
&	\\
\ruLine(1,2)<{\eend{C}{a_1}}&	&	&	
&{\scriptstyle C[\range{a_1}{a_n};a]}&	&	&\rdLine(1,2)>{\eend{C}{a_n}}\\
\cnr		&	&	&\rLine_{\eend{C}{a}}&	&	&	&\cnr\\
\end{opetope}
of $V$
\item for each $a, \range{a_1}{a_n},
\range{\range{a_1^1}{a_1^{k_1}}}{\range{a_n^1}{a_n^{k_n}}} \in C_0$, a 3-cell
\begin{equation}	\label{eq:comp-three}
\piccy{biasedcomp3}
\end{equation}
in $V$
\item for each $a\in C_0$, a 3-cell
\[
\piccy{biasedid3}
\]
in $V$
\end{itemize}
satisfying associativity and identity axioms.

Before moving on to specific examples, we take a brief look at how the
definition of enriched plain multicategory reads when $V$ is in some way
degenerate. As for enriched categories, we might as well assume that $V$ is
vertically discrete.

\begin{enumerate}	\label{egs:enr-deg-fm}
\item
Suppose that $V$ is vertically trivial, i.e.\ a $T_2$-multicategory
(\ref{egs:degen-fm}\bref{eg:V-T2}). Then a plain multicategory $C$ enriched
in $V$ consists of
\begin{itemize}
\item a set $C_0$
\item for each $\range{a_1}{a_n},a\in C_0$, an object 
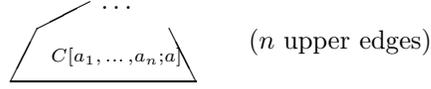
\[
\begin{opetope}
		&	&	&\cnr 	&\ldots	&	&	&	\\
		&\cnr	&\ruLine(2,1)&	&	&	&\cnr 	&	\\
\ruLine(1,2)	&	&	&	
&{\scriptstyle C[\range{a_1}{a_n};a]}	&	&	&\rdLine(1,2)\\
\cnr		&	&	&\rLine	&	&	&	&\cnr\\
\end{opetope}
\diagspace
\mbox{($n$ upper edges)}
\]
of $V$
\item for each $a, \range{a_1}{a_n},
\range{\range{a_1^1}{a_1^{k_1}}}{\range{a_n^1}{a_n^{k_n}}} \in C_0$, an arrow
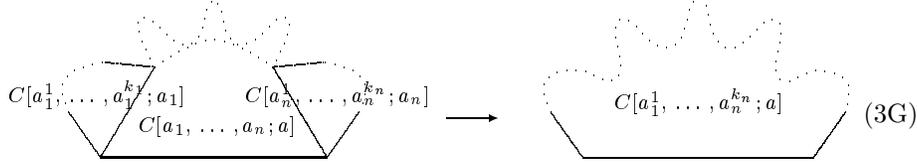
\begin{equation}	
\mbox{\hspace*{-25mm}}\piccy{t2comp3}
\end{equation}
in $V$
\item for each $a\in C_0$, an arrow
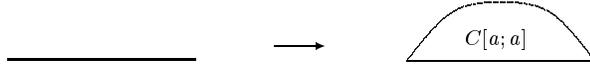
\[
\piccy{t2id3}
\]
in $V$,
\end{itemize}
satisfying the inevitable associativity and identity axioms.

\item	\label{eg:V-from-smc}
Suppose $V$ comes from a symmetric monoidal category $W$
(\ref{egs:degen-fm}\bref{eg:smc-fm}). Then a plain multicategory $C$
enriched in $V$ consists of
\begin{itemize}
\item a set $C_0$
\item for each $\range{a_1}{a_n},a\in C_0$, an object $C[\range{a_1}{a_n};a]$
of $W$
\item for each $a, a_i, a_i^j$, a `composition' arrow 
\begin{eqnarray*}
C[\range{a_1}{a_n};a] \otimes 
C[\range{a_1^1}{a_1^{k_1}};a_1] \otimes \cdots \otimes
C[\range{a_n^1}{a_n^{k_n}};a_n] \\
\go C[\range{a_1^1}{a_n^{k_n}};a]
\end{eqnarray*}
in $W$
\item for each $a$, an `identity' arrow $I\go C[a;a]$
\end{itemize}
satisfying associativity and identity.

A plain multicategory enriched in the symmetric monoidal category
\pr{\Set}{\times} is therefore just the same thing as a plain multicategory. 

\item	\label{eg:enr-smc}
Topologists are used to considering plain operads enriched in a symmetric
monoidal category $W$ (e.g.\ \pr{\fcat{Spaces}}{\times}): that is, the case of
\bref{eg:V-from-smc} where $C_0 = 1$. (They call them just
`(non-symmetric) operads in $W$'.)  Such a structure consists of a sequence
$(C(n))_{n\in\nat}$ of objects of $W$, together with morphisms
\begin{eqnarray*}
C(n)\otimes C(k_1)\otimes\cdots\otimes C(k_n) &\go &C(k_1 + \cdots + k_n) \\
I &\go &C(1)
\end{eqnarray*}
satisfying the usual axioms.

\end{enumerate}

\section{Enriched Plain Multicategories: Examples}	\label{sec:epme}

We finish the material on enriched plain multicategories with three
examples. In Chapter~\ref{sec:relaxed} we will meet another important
example: relaxed multicategories are plain multicategories enriched in
the $T_2$-multicategory \ftrcat{\fcat{TR}^{\op}}{\Set}.

\begin{eg}{egs:epms}
\item	\label{eg:str-mon-two}
We will be particularly interested in plain operads enriched in the symmetric
monoidal category \pr{\Cat}{\times}. These are exactly the same things as
\pr{\Cat}{\ust}-operads, where \ust\ is the free strict monoidal category
monad given by $A^* = \coprod_{n\in\nat}A^n$. They are also the same as
$T_2$-structured categories. Both of these equivalences are established in
\cite[IV.2]{SHDCT}. 

\item 
Modules over a fixed commutative ring $R$ form a plain multicategory $C$
enriched in the symmetric monoidal category \Ab\ of abelian groups: the
objects are $R$-modules, and $C[\range{X_1}{X_n};X]$ is the abelian group of
$R$-module homomorphisms $X_1\otimes\cdots\otimes X_n \go X$.

\item
Let $V$ be the \fm-multicategory of \ref{egs:degen-fm}\bref{eg:bim-fm}, and
take any plain multicategory $C$ enriched in the symmetric monoidal category
\Ab. Then there arises a plain multicategory $C'$ enriched in $V$. (This is
analogous to the situation one level down, as described in
\ref{egs:specific-enr-cats}\bref{eg:spec-bim-ab}, but we do not attempt a
general result.) So, $C'_0=C_0$; \eend{C'}{a} is the 1-cell, i.e.\ ring,
\ehom{C}{a}{a}; and \ehom{C'}{\range{a_1}{a_n}}{a} is the abelian group
\ehom{C}{\range{a_1}{a_n}}{a} acted on by \ehom{C}{a}{a} on the left and
\ehom{C}{a_i}{a_i} on the right, both actions being by composition in $C$.

\end{eg}

\section{The Hierarchies}	\label{sec:hier}

The process $\Cartpr \goesto \Cartprp$ can be applied indefinitely. In
particular, if we start with \Zeropr\ then we obtain a hierarchy \Hpn{n} of
monads on categories:
\begin{itemize}
\item $\Hpn{0}=\Zeropr$
\item $\Hpn{n+1} = \pr{P_n \hyph\Gph}{\mbox{free $P_n$-multicategory}}.$
\end{itemize}
(The persistence property of suitability, stated in
Theorem~\ref{thm:free-main}, makes this possible.) We can therefore discuss
$P_n$-multicategories enriched in $P_{n+1}$-multicategories, for all $n$.

\begin{figure}
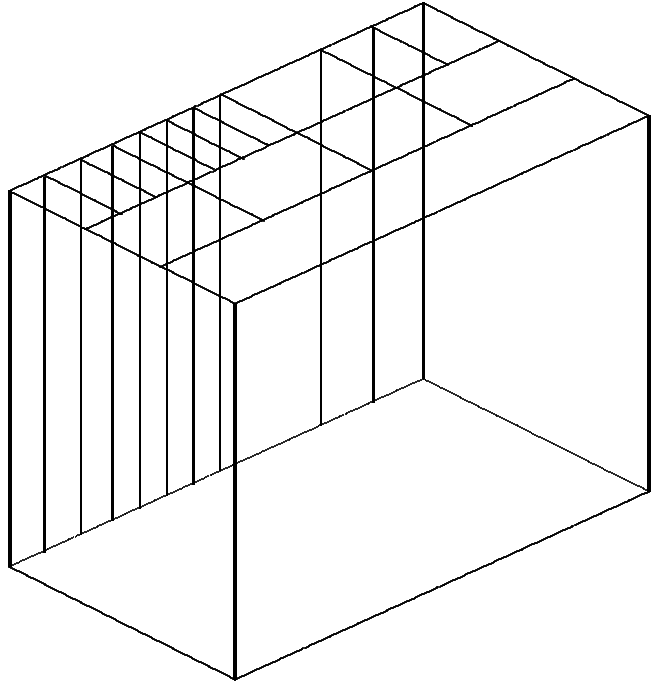

\centerline{\piccy{p23cell}}
\caption{The shape of a 3-cell in a $P_2$-graph}
\label{fig:p-two}
\end{figure}

As might be expected from the complexity of the structures in
Chapter~\ref{sec:cats}, which was just the case $n=0$, this rapidly becomes
difficult in at least visual terms. For instance, take $n=1$. Then  $P_1 =
\fc$, so we are interested in \fc-multicategories enriched in a
$P_2$-multicategory $V$. The graph of a $P_2$-multicategory looks like
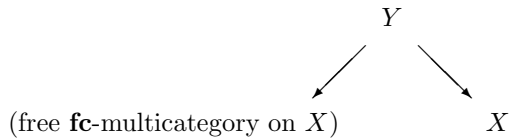
\begin{slopeydiag}
	&	&Y	&	&	\\
	&\ldTo	&	&\rdTo	&	\\
\mbox{(free \fc-multicategory on $X$)}	
&	&	&	&X	\\
\end{slopeydiag}
where $X$ and $Y$ are \fc-graphs (i.e.\ objects of $\Eee_2$): thus a
$P_2$-graph is defined by $2^3=8$ sets and various functions between them. We
picture a $P_2$-graph as some kind of 3-dimensional cubical structure; a
3-cell (without labels) is drawn in Figure~\ref{fig:p-two}.

More tractable than the full hierarchy \Hpn{n} is the restricted hierarchy
\Cpn{n}, which was constructed so that
\begin{itemize}	\label{p:Cpn-property}
\item $\Cpn{0}=\Zeropr$
\item $\Cpn{n+1} = \pr{\mbox{$T_n$-graphs on 1}}{\mbox{free $T_n$-operad}}$
\end{itemize}
(see~\ref{sec:opetopes}). Since ($T_n$-graphs on 1) is a subcategory of $T_n
\hyph\Gph$ and (free $T_n$-operad) is the restriction of (free
$T_n$-multicategory) to this subcategory, we can talk about
$T_n$-multicategories enriched in $T_{n+1}$-multicategories, for any $n$. Put
another way, a $T_{n+1}$-multicategory is a special kind of
$T'$-multicategory, where $\Cartpr=\Cpn{n}$; this is the special case we
considered in~\ref{egs:enr-deg}\bref{eg:deg-plain} (categories enriched in
plain multicategories) and~\ref{egs:degen-fm}\bref{eg:V-T2} (plain
multicategories enriched in $T_2$-multicategories).

We finish this chapter with some remarks on symmetric monoidal categories. In
Chapter~\ref{sec:prelims} we observed that any monoidal category yields a
plain multicategory (canonically, up to isomorphism), and
in~\ref{egs:degen-fm}\bref{eg:smc-fm} we also saw that any symmetric monoidal
category yields a $T_2$-multicategory. In fact, the method of
\ref{egs:degen-fm}\bref{eg:smc-fm} suggests that a symmetric monoidal
category should yield a $T_n$-multicategory for all $n$. This is indeed the
case:
\begin{propn}	\label{propn:smc-tn}
Let $n\in\nat$ and let $W$ be a symmetric monoidal category. Then there is an
associated $T_n$-multicategory $V$, defined canonically up to isomorphism.
\end{propn}
\textbf{Sketch proof}\ 
By induction on $n$, using the explicit free multicategory construction in
the Appendix, we show that for each $n$ there is a (canonical) natural
transformation
\begin{ntdiag}
\Set/S_n	&	&\rTo^{T_n}	&	&\Set/S_n	\\
		&	&		&\ 	&		\\
\dTo<{U}	&	&\ldNT>{\theta}	&	&\dTo>{U}	\\
		&\ 	&		&	&		\\
\Set		&	&\rTo_{M}	&	&\Set,		\\
\end{ntdiag}
where $U$ is the forgetful functor and $M$ is the free commutative monoid
functor. This pair \pr{U}{\theta} is, in fact, a monad opfunctor $\Cpn{n} \go
\pr{\Set}{M}$. Intuitively, what $\theta$ does is to send a diagram like
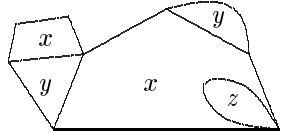
\[
\piccy{regionlabpd}
\]
($n=2$) to the set-with-multiplicities $\{x,x,y,y,z\}$.

We now follow the strategy of \ref{egs:degen-fm}\bref{eg:smc-fm}. The
objects-object of $V$ is to be \bktdvslob{W_0 \times S_n}{}{S_n}: so for any
$s\in S_n$, an object of $V$ over $s$ is just an object of $W$. By choosing
a function $M(W_0) \go W_0$ which takes a set-with-multiplicities of objects
and tensors them in some order, we obtain via $\theta$ a `tensor' map
\[
T_n \bktdvslob{W_0 \times S_n}{}{S_n} 
\go \bktdvslob{W_0 \times S_n}{}{S_n} .
\]
Performing a similar operation for arrows, we get a $T_n$-multicategory $V$,
which can be thought of as some kind of `weak $T_n$-structured
category'. Different choices of the tensor function $M(W_0) \go W_0$ only
affect $V$ up to isomorphism.\begin{flushright}$\Box$\end{flushright}

An explanation of this result can be given in terms of the `periodic table'
(\cite{BDC}, \cite{Sim}), which displays for each $0\leq k\leq n$ what kind
of a structure an $n$-category is when it only has one 0-cell, one 1-cell,
\ldots, one $k$-cell. The table suggests that for any $n\geq 1$, a symmetric
monoidal category gives rise to an $n$-category with only one $d$-cell for
$d\leq n-2$. (When $n=1$, of course, we don't even need `monoidal', and when
$n=2$ we don't need `symmetric'.) But according to the opetopic definition of
$n$-category given in \cite{BD} or \cite{HMP}, an $n$-category with such a
degeneracy is a $T_{n-1}$-multicategory with certain universality
properties. Putting these together, a symmetric monoidal category gives rise
to a $T_{n-1}$-multicategory for any $n\geq 1$, that is, to a
$T_n$-multicategory for any $n\in\nat$.

\chapter{Relaxed Multicategories}	\label{sec:relaxed}

\begin{figure}[t]
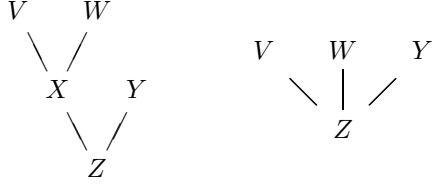

\[
\begin{diagram}[height=1.5em,width=1.5em]
V	&	&W	&	\\
	&\rt{1}\lt{1}&	&	\\
	&X	&	&Y	\\
	&	&\rt{1}\lt{1}&	\\
	&	&Z	&	\\
\end{diagram}
\diagspace\diagspace
\begin{diagram}[height=1.5em,width=1.5em]
V	&	&W	&	&Y	\\
	&\rt{2}	&\dn	&\lt{2}	&	\\
	&	&Z	&	&	\\
\end{diagram}
\]
\caption{Composing singular functions}
\label{fig:sing-trees}
\end{figure}

Relaxed multicategories can be thought of, according to~\cite{Bor}, as
multicategories in which the morphisms might have some sort of
singularity. In a genuine multicategory, two morphisms $V,W\goby{f} X$ and
$X, Y \goby{g} Z$ have a composite $V,W,Y \go Z$. But imagine now that the
objects of the multicategory are spaces of some kind and that an arrow
$\range{X_1}{X_n} \go X$ is a function $X_1 \times\cdots\times X_n \go X$
which might have
singularities of an `allowable' kind. Then there's a composite $g\of (f\times
1): V\times W\times Y \go Z$ of $f$ and $g$, but it might have a singularity
of a more severe kind---compare the fact that the composite of two
meromorphic functions need not be meromorphic. So $g\of (f\times 1)$ lies in
the set of functions $V\times W\times Y \go Z$ `with singularity of type $
\begin{tree}
\node	&	&\node	&	\\
	&\rt{1}\lt{1}&	&	\\
	&\node	&	&\node	\\
	&	&\rt{1}\lt{1}&	\\
	&	&\node	&	\\
\end{tree}
$
'; this is a larger set than that of functions $V\times W\times Y \go Z$ with
singularities of the allowable kind, which we call `singularities of type 
$
\begin{tree}
\node	&	&\node	&	&\node	\\
	&\rt{2}	&\dn	&\lt{2}	&\ \ \ \mbox{'}	\\
	&	&\node	&	&	\\
\end{tree}
$
(see Figure~\ref{fig:sing-trees}).
In general, for each $n$-leafed tree $\tau$ and sequence of objects
$\range{X_1}{X_n},X$ there's a set \relhom{\Hom}{\tau}{X_1}{X_n}{X} of
functions with singularity type $\tau$, and if $\tau' \go \tau$ is a map of
trees then there's a map
\[
\relhom{\Hom}{\tau}{X_1}{X_n}{X} \go \relhom{\Hom}{\tau'}{X_1}{X_n}{X}
\]
---in the present case, an inclusion, corresponding to the tree map
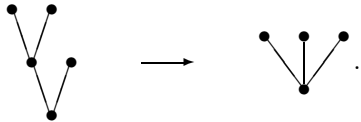
\[
\begin{tree}
\node	&	&\node	&	\\
	&\rt{1}\lt{1}&	&	\\
	&\node	&	&\node	\\
	&	&\rt{1}\lt{1}&	\\
	&	&\node	&	\\
\end{tree}
\diagspace\go\diagspace
\begin{tree}
\node	&	&\node	&	&\node	\\
	&\rt{2}	&\dn	&\lt{2}	&\ \ \ \mbox{.}	\\
	&	&\node	&	&	\\
\end{tree}
\]
The whole structure is called a relaxed multicategory
(Definition~\ref{defn:rel-mti}).

The first occurrence of this definition, or something close, was apparently
in Beilinson and Drinfeld's \cite{BeDr}. This is cited in Soibelman's
\cite{Soi}, where the notion seems to have been adapted somewhat. Both papers
use the term `pseudo-monoidal category'; Borcherds, working independently,
called them relaxed multicategories. The definition of the category of trees
is different in each of \cite{Bor}, \cite{Soi} (where the definition traces
back to \cite{KMGWC} and \cite{KMQCP}) and the present work. My understanding
is that the definition in \cite{Bor} is \emph{ad hoc} enough that it could be
replaced by the one used here, but I do not know whether this is also the
case for \cite{Soi}.

In~\ref{sec:basics} we define relaxed multicategories and give some
examples. Section~\ref{sec:rel-via-enr} shows that a relaxed multicategory is
just a plain multicategory enriched in the $T_2$-multicategory \tropset, by
which means we define relaxed $T_n$-multicategories for all $n$ (which are
relaxed multicategories when $n=1$). We pay particular attention to `relaxed
categories', the case $n=0$
(section~\ref{sec:rel-cats}). Sections~\ref{sec:rel-mon-cats}
and~\ref{sec:fun} present dual ways in which relaxed multicategories can
arise, which we will call `modifying the domain' and `modifying the codomain'
respectively.

\section{Basics}	\label{sec:basics}

\begin{defn}	\label{defn:rel-mti}
A \emph{relaxed multicategory} $C$ consists of:
\begin{itemize}
\item a set $C_0$ of objects
\item for each $\tau\in\TR{n}$ and $\range{a_1}{a_n},a\in C_0$, a
set \multihom{C_{\tau}}{\range{a_1}{a_n}}{a}
\item for each map $\tau' \goby{\phi} \tau$ in \TR{n}, a function
\begin{eqnarray*}
\multihom{C_{\tau}}{\range{a_1}{a_n}}{a}
&\go &\multihom{C_{\tau'}}{\range{a_1}{a_n}}{a}	\\
f	&\goesto	&f.\phi
\end{eqnarray*}
\item for each $\tau\in\TR{n}$, $\tau_1\in\TR{k_1}$, \ldots,
$\tau_n\in\TR{k_n}$, a `composition' function
\begin{eqnarray*}
\multihom{C_{\tau_1}}{\range{a_1^1}{a_1^{k_1}}}{a_1} \times \cdots
\times 
\multihom{C_{\tau_n}}{\range{a_n^1}{a_n^{k_n}}}{a_n} \times
\multihom{C_{\tau}}{\range{a_1}{a_n}}{a} \\
\go
\multihom{C_{\tau\of\bftuple{\tau_1}{\tau_n}}}{\range{a_1^1}{a_n^{k_n}}}{a},
\end{eqnarray*}
$\tuplebts{\range{f_1}{f_n},f} \goesto f\of\bftuple{f_1}{f_n}$
\item for each $a\in C_0$, an `identity' element $1_a$ of
\multihom{C_{\utree}}{a}{a},
\end{itemize}
subject to the following axioms:
\begin{itemize}
\item $f.(\phi\of\phi')=(f.\phi).\phi'$ and $f.1=f$, for any $f$, $\phi'$ and
$\phi$ for which these make sense
\item composition is associative, and the `identity' elements are two-sided
identities for composition (i.e.\ $f\of\bftuple{1}{1} = f = 1\of f$)
\item compatibility:
\[
(f\of\bftuple{f_1}{f_n}).(\phi\of\bftuple{\phi_1}{\phi_n}) = 
(f.\phi)\of\bftuple{f_1.\phi_1}{f_n.\phi_n}
\]
for any $f\in \multihom{C_{\tau}}{\range{a_1}{a_n}}{a}$, $f_i \in
\multihom{C_{\tau_i}}{\range{a_i^1}{a_i^{k_i}}}{a_i}$,
$\tau' \goby{\phi} \tau$ in \TR{n}, and $\tau'_i \goby{\phi_i} \tau_i$ in
\TR{k_i} ($1\leq i\leq n$). 
\end{itemize}
\end{defn}

\begin{eg}{egs:rel-mtis}
\item	\label{eg:triv-rel-mti}
If $E$ is any monoidal category, then putting $C_0=E_0$ and
\[
\relhom{C}{\tau}{a_1}{a_n}{a} = \homset{E}{a_1\otimes\cdots\otimes a_n}{a}
\]
for all $\tau\in\TR{n}$ makes $C$ into a relaxed multicategory.

\item	\label{eg:monoid}
Let $\tau\in\TR{n}$, and number the $n$ leaves of $\tau$ as \range{1}{n},
working from left to right. Let $h_i(\tau)$ be the height of the $i$th leaf:
that is, the length of the path from the leaf down to the root. Let $M$ be a
monoid. Then there is a relaxed multicategory whose objects are sets, and
with 
\[
\relhom{\Hom}{\tau}{X_1}{X_n}{X} =
\prod_{1\leq i\leq n} \homset{\Set}{M^{h_i(\tau)} \times X_i}{X}.
\]
To see what the maps between homsets induced by maps of trees are, let
$\tau'\go\tau$ be a map in \TR{n} and let $1\leq i \leq n$. Then there is an
induced map $M^{h_i(\tau')} \go M^{h_i(\tau)}$, as can be seen informally by
thinking of a tree map as contracting internal edges and expanding nodes (see
\cite[IV.3]{SHDCT}), and will be shown in
\ref{egs:rel-mtis-from-mons}\bref{eg:new-monoid}. This in turn induces the
map of homsets, and composition and identities are defined in a natural
way. (Examples \bref{eg:monoid}--\bref{eg:R-algebra} are explained further in
\ref{egs:rel-mtis-from-mons}\bref{eg:new-monoid}--\bref{eg:new-R-algebra}.)

\item 	\label{eg:exceptions} 
Let $v(\tau)$ denote the number of \emph{internal vertices} of a tree $\tau$:
that is, those vertices which are not leaves. Any map $\tau'\go\tau$ induces
a map from the set of internal vertices of $\tau'$ to the set of internal
vertices of $\tau$, as may again be seen informally, and again will be gone
into in more detail later
(\ref{egs:rel-mtis-from-mons}\bref{eg:together-exceptions}). We also have the
identity
\[
v(\tau\of\bftuple{\tau_1}{\tau_n}) =
v(\tau) + v(\tau_1) + \cdots + v(\tau_n).
\]
These facts together mean that there is a relaxed multicategory whose objects
are sets, and with
\[
\relhom{\Hom}{\tau}{X_1}{X_n}{X} = 
X^{v(\tau)} \times \prod_{1\leq i \leq n} \homset{\Set}{X_i}{X}.
\]

\item	\label{eg:R-algebra}
If $R$ is a commutative ring and $A$ an $R$-algebra, then there is a
relaxed multicategory whose objects are $R$-modules, and with
\[
\relhom{\Hom}{\tau}{X_1}{X_n}{X} =
\homset{\Hom_R}{X_1 \otimes\cdots\otimes X_n \otimes A^{\otimes v(\tau)}}{X}.
\]

\item	\label{eg:rep}
The motivating example of a relaxed multicategory in \cite{Bor} is as
follows. A \emph{vertex group} over a commutative ring $R$ consists, roughly,
of a Hopf algebra $H$ over $R$ and a two-sided $H$-module $K$ which also has
the structure of an algebra over $H^* = \homset{\Hom_R}{H}{R}$. Borcherds
defines a relaxed multicategory $\mb{Rep}(G)$ of representations of any
vertex group $G=\pr{H}{K}$, and then defines a \emph{$G$-vertex algebra} to
be a commutative monoid in $\mb{Rep}(G)$ (in a suitable sense). For a
certain choice of $G$, the standard kind of vertex algebra is the
same as a $G$-vertex algebra (see \cite{Sny}). 

\end{eg}

In \ref{sec:fun} we will look more closely at how the relaxed multicategory
of~\bref{eg:rep} arises, and see that the general method is `dual' to the
general method in examples~\bref{eg:monoid}--\bref{eg:R-algebra}.

\section{Relaxation via Enrichment}	\label{sec:rel-via-enr}

In this section we show that relaxed multicategories are a completely natural
idea (`no artificial ingredients') from the point of view of general
multicategory theory. Specifically, we first show that relaxed
multicategories are just plain multicategories enriched in the
$T_2$-multicategory \tropset, and then we observe that \tropset\ is
`completely natural'. This enables us to give a definition of relaxed
$T_n$-multicategory for any $n\in\nat$, the familiar case being $n=1$.

Our first task is to define \tropset. We have already seen
(\ref{egs:degen-fm}\bref{eg:smc-fm}) that any symmetric monoidal category,
and in particular \pr{\Set}{\times}, is naturally a $T_2$-multicategory. We
have also met the $T_2$-structured category \fcat{TR}, and by swapping \dom\
and \cod\ we obtain an opposite $T_2$-structured category,
$\fcat{TR}^{\op}$. (Thus the category $\fcat{TR}^{\op}(n)$ is the opposite of
the category \TR{n}.) Then \tropset\ is the exponential in
$T_2\hyph\Multicat$, the existence of which is guaranteed by:
\begin{thm}	\label{thm:exp}
Let $S$ be a set and $T$ a cartesian monad on $\Set/S$. Then any
$T$-structured category (or rather, its underlying $T$-multicategory) is
exponentiable in $T\hyph\Multicat$. 
\end{thm}
\pf\ Largely omitted. Any \pr{\Set/S}{T}-multicategory $B$ has an underlying
category object $|B|$ in $\Set/S$, so for each $s\in S$ there is a category
$|B|(s)$. If $A$ is a $T$-structured category, the objects over $s\in S$ of
the exponential \ftrcat{A}{B} are the functors $|A|(s) \go |B|(s)$. We can
then construct the arrows and the composition; for the latter, we need $A$ to
be a structured category rather than just a multicategory.
Further hints are given at the beginning of~\ref{sec:rel-cats}.\done

\emph{Remark:}\/ We will only use the exponential \ftrcat{A}{B} when $A$ is
small, even though $B$ might be large. This adds a measure of safety to our
otherwise cavalier attitude to size. 

Specifically, we need to know what the \pr{\Set/\nat}{T_2}-multicategory
\tropset\ looks like. Following the description of $T_2$-multicategories in
\cite[p.\ 66]{SHDCT}:
\begin{itemize}
\item an object over $n$ is a functor $\TR{n}^{\op} \go \Set$
\item an arrow 
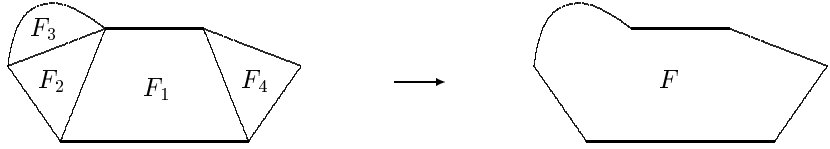
\[
\piccy{tropsetarrow}
\]
(for instance) is a family of functions
\[
F_1\tau_1 \times
F_2\tau_2 \times
F_3\tau_3 \times
F_4\tau_4
\go
F(\tau_1\of(\tau_2\of(1,\tau_3),1,\tau_4)),
\]
one for each $\tau_1\in\TR{3}, \tau_2\in\TR{2}, \tau_3\in\TR{1},
\tau_4\in\TR{2}$, which is natural in the $\tau_i$'s
\item composition and identities are as of functions, i.e.\ come from the
\pr{\Set/\nat}{T_2}-multicategory \Set.
\end{itemize}

We can now demonstrate that a plain multicategory enriched in \tropset\ is
the same thing as a relaxed multicategory. For the former consists of
\begin{itemize}
\item a set $C_0$
\item for each $\range{a_1}{a_n},a\in C_0$, a functor 
\[\begin{array}{rrcl}
C_{\_}[\range{a_1}{a_n};a]: 	&\TR{n}^{\op} 	&\go 	&\Set,	\\
	&\tau	&\goesto	&C_{\tau}[\range{a_1}{a_n};a]
\end{array}\]
\item for each $a, a_1, \ldots, a_n, a_1^1, \ldots, a_n^{k_n} \in C_0$, a
family of functions 
\begin{eqnarray*}
C_{\tau_1}[\range{a_1^1}{a_1^{k_1}};a_1] \times\cdots\times
C_{\tau_n}[\range{a_n^1}{a_n^{k_n}};a_n] \times
C_{\tau}[\range{a_1}{a_n};a] \\
 \go 
C_{\tau\of\bftuple{\tau_1}{\tau_n}}[\range{a_1^1}{a_n^{k_n}};a],
\end{eqnarray*}
one for each $\tau\in\TR{n}, \tau_1\in\TR{k_1}, \ldots, \tau_n\in\TR{k_n}$,
which is natural in the $\tau_i$'s and $\tau$
\item for each $a\in C_0$, a function $1\go C_{\utree}[a;a]$, i.e.\ an
element of $C_{\utree}[a;a]$,
\end{itemize}
such that associativity and identity axioms hold. This is exactly what a
relaxed multicategory is.

Three generalizations of the notion of $T_1$-multicategories enriched in
\tropset\ now present themselves.

Firstly, \Set\ could be changed to any other symmetric monoidal category. By
changing \pr{\Set}{\times} to \pr{\Ab}{\otimes} we obtain what Borcherds
calls a relaxed multi\emph{linear} category: thus the homsets
$C_{\tau}[\range{a_1}{a_n};a]$ are not just sets but abelian groups.

Secondly, \fcat{TR} could be changed to any other $T_2$-structured category
\cat{T}. Then we obtain what Soibelman calls a `\cat{T}-pseudo monoidal
category' in \cite{Soi}.

Thirdly, we have been discussing enriched $T_1$-multicategories, but there is
nothing special here about the number 1: it could be replaced by any $n$, as
follows.

Recall that \fcat{TR} is not just any old $T_2$-structured category, but in
fact the free such on the terminal $T_2$-multicategory: in the terminology
of~\ref{sec:opetopes}, $\fcat{TR}=\PD{2}$. So a relaxed multicategory is a
$T_1$-multicategory enriched in the $T_2$-multicategory
\ftrcat{\PD{2}^{\op}}{\Set}. Recall too that a symmetric monoidal category
gives not just a $T_2$-multicategory, but a $T_n$-multicategory for any $n$
(Proposition~\ref{propn:smc-tn}). We may therefore generalize:
\begin{defn}	\label{defn:rel-tn-mti}
Let $n\in\nat$. A \emph{relaxed $T_n$-multicategory} is a
$T_n$-multicategory enriched in the $T_{n+1}$-multicategory
\ftrcat{\PD{n+1}^{\op}}{\Set}.
\end{defn} 
Thus a relaxed $T_1$-multicategory is a relaxed multicategory. By thinking in
terms of pasting diagrams rather than trees, the basic idea for $n\geq 2$
becomes apparent. The case $n=0$ is also illuminating, and this is the
subject of the next section.

\section{Relaxed Categories}	\label{sec:rel-cats} 

Since a $T_0$-multicategory is a category, relaxed $T_0$-multicategories will
be called, perhaps rather disturbingly, \emph{relaxed categories}. These are
categories enriched in the plain multicategory \dopset, so we first of all
need to know what \dopset\ is; but in order to do \emph{this} we make a short
digression on Theorem~\ref{thm:exp}.

Let $A$ be any strict monoidal category and $B$ any plain
multicategory. Theorem~\ref{thm:exp} guarantees that the exponential
\ftrcat{A}{B} of multicategories exists. One might think that an object of
\ftrcat{A}{B} would be a multicategory map $A\go B$, but this is not the
case. For multicategory maps $A\go B$ correspond to multicategory maps
$1\go\ftrcat{A}{B}$, and a map from 1 (the terminal multicategory) to a
multicategory $C$ is a `monoid' in $C$, in other words, an object $c$ of $C$
together with arrows $c,c\go c$ and $\cdot \go c$ satisfying the usual
axioms. (When $C$ is a monoidal category, this is just B\'{e}nabou's
observation that a lax functor $1\go C$ is a monoid in $C$:
\cite[5.4.2]{Ben}.) So \ftrcat{A}{B} has the property that a \emph{monoid} in
it is a multicategory map $A\go B$.

However, if we define $I$ to be the multicategory with one object and just
one arrow (the identity), then a map $I\go C$ of multicategories is the same
as an object of $C$. Thus the objects of \ftrcat{A}{B} are the multicategory
maps $I\times A \go B$: equivalently, they are the functors $|A|\go|B|$,
where $|\cdot|$ indicates the underlying category of a plain multicategory
obtained by ignoring all non-unary arrows.

An object of the multicategory \dopset\ is, therefore, a functor
$\Delta^{\op} \go \Set$. We can also compute that an arrow $\range{F_1}{F_n}
\goby{\phi} F$ is a family of functions
\[
F_1 k_1 \times\cdots\times F_n k_n
\goby{\phi_{\range{k_1}{k_n}}}
F(k_1 + \cdots + k_n)
\] 
($\range{k_1}{k_n}\in\nat$) which are natural in the $k_i$'s: that is,
if $k'_i \goby{p_i} k_i$ is an arrow in $\Delta$ for each $i$, then 
\begin{diagram}
F_1 k_1 \times\cdots\times F_n k_n	&\rTo^{\phi_{\range{k_1}{k_n}}}
&F(k_1 + \cdots + k_n)		\\
\dTo<{F_1 p_1 \times\cdots\times F_n p_n}&
&\dTo>{F(p_1 + \cdots + p_n)}	\\
F_1 k'_1 \times\cdots\times F_n k'_n	&\rTo_{\phi_{\range{k'_1}{k'_n}}}
&F(k'_1 + \cdots + k'_n)	\\
\end{diagram}
commutes. Composition and identities are got by composing the functions
$\phi_{\range{k_1}{k_n}}$. 

(In fact, \dopset\ is not just a multicategory but
a monoidal category, although we don't need to know this. The monoidal
structure is the Day tensor product. Note in particular that it is not the
cartesian product, so relaxed categories are not the same as the simplicially
enriched categories sometimes considered in homotopy theory.)


A relaxed category $C$ consists, therefore, of
\begin{enumerate}
\item a set $C_0$
\item for each $a,b\in C_0$ and $n\in\nat$, a set \ehom{C_n}{a}{b}
\item 	\label{part:restr}
for each map $n'\go n$ in $\Delta$, a function $\ehom{C_n}{a}{b} \go
\ehom{C_{n'}}{a}{b}$ 
\item for each $a,b,c,m,n$, a `composition' function
\[
\ehom{C_m}{a}{b} \times \ehom{C_n}{b}{c}
\go
\ehom{C_{m+n}}{a}{c}
\]
\item for each $a$, an `identity' element of \ehom{C_0}{a}{a}
\end{enumerate}
such that the assignment in~\bref{part:restr} is functorial, the composition
and identities obey associativity and identity laws, and the following
compatibility rule holds: if $m'\go m$ and $n'\go n$ in $\Delta$ then
\begin{diagram}
\ehom{C_m}{a}{b} \times \ehom{C_n}{b}{c}	&\rTo	
&\ehom{C_{m+n}}{a}{c}	\\
\dTo	&	&\dTo	\\
\ehom{C_{m'}}{a}{b} \times \ehom{C_{n'}}{b}{c}	&\rTo	
&\ehom{C_{m'+n'}}{a}{c}	\\
\end{diagram}
commutes.

A relaxed category can be thought of, roughly, as an ordinary category in
which each map has a degree (a natural number) and composing maps adds the
degrees. 

\begin{eg}{egs:rel-cats}
\item	\label{eg:cat-to-rel}
Any category $D$ gives two different relaxed categories, $C$ and $C'$. Both
have the same object-set as $D$, and then
\begin{eqnarray*}
\ehom{C_n}{a}{b}	&=	&\homset{D}{a}{b} \mbox{\ for all $n$,}\\
\ehom{C'_n}{a}{b}	&=	&
\left\{
\begin{array}{ll}
\homset{D}{a}{b}&\mbox{if $n=0$}	\\
\emptyset	&\mr{otherwise.}
\end{array}
\right.
\end{eqnarray*}

\item 	\label{eg:craig}
(This example is due to Craig Snydal.) Let $\cat{C}_0$ be a/the collection of
abelian groups, and define a relaxed category \cat{C} by 
\[
\ehom{\cat{C}_n}{A}{B} = \homset{\Ab}{A}{B[\range{x_1}{x_n}]},
\]
where $B[\range{x_1}{x_n}]$ is the abelian group of polynomials in $n$
commuting variables with coefficients in $B$. A map $n'\go n$ in $\Delta$
induces a map $B[\range{x_1}{x_n}] \go B[\range{y_1}{y_{n'}}]$ of abelian
groups, and so induces the required map of the enriched homsets. (The exact
formula for this can be recovered from~\bref{eg:comonad-rel} below.) For
composition, suppose $f\in\homset{\Ab}{A}{B[\range{x_1}{x_n}]}$ and
$g\in\homset{\Ab}{B}{C[\range{y_1}{y_m}]}$: then the composite $g\of f$ in
\cat{C} is the composite
\[
A \goby{f} B[\range{x_1}{x_n}] \goby{g[\range{x_1}{x_n}]}
C[\range{y_1}{y_m}][\range{x_1}{x_n}] \iso C[\range{z_1}{z_{m+n}}]
\]
in \Ab.

\item	\label{eg:craig-dual}
Dually, there is a relaxed category \cat{C'} whose objects are again abelian
groups, and with
\[
\ehom{\cat{C'}_n}{A}{B} = \homset{\Ab}{A\odblbkt\range{x_1}{x_n}\cdblbkt}{B}
\]
where $A\odblbkt\range{x_1}{x_n}\cdblbkt$ is the abelian group of formal
power series.

\item 	\label{eg:comonad-rel}
Let $D$ be a category with a comonad $G$ on it: then there's an associated
relaxed category $C$ with $C_0=D_0$,
\[
\ehom{C_n}{a}{b} = \homset{D}{a}{G^{n}b},
\]
composition defined by $G$ being a functor, and the maps $\ehom{C_n}{a}{b}
\go \ehom{C_{n'}}{a}{b}$ defined by $G$ being a
comonad. Example~\bref{eg:craig} is exactly this, with $G$ the comonad
$B\goesto B[x]$ on \Ab. This \emph{is} a comonad: for \nat, like any other
set, is a comonoid in \pr{\Set}{\times}, so the copower functor
$\nat\times\dashbk : \Ab\go\Ab$ has the structure of a comonad, and
$\nat\times B = B[x]$.

\item	\label{eg:monad-rel-cat}
Dually, if $D$ is a category with a monad $T$ on it, then there arises a
relaxed category $C'$ with $C'_0=D_0$ and
\[
\ehom{C'_n}{a}{b} = \homset{D}{T^{n}a}{b}.
\]
Example~\bref{eg:craig-dual} exhibits this when $T$ is the monad $A\goesto
A\odblbkt x\cdblbkt = A^{\nat}$ on \Ab.

\end{eg}

\section{Relaxed Monoidal Categories}	\label{sec:rel-mon-cats}

A relaxed multicategory is like a multicategory, except that composition
behaves less uniformly than usual (in the sense of the introduction to this
chapter). Similarly, a relaxed monoidal category will be like a monoidal
category except that the tensor (written \atsr) behaves less uniformly than
usual. As the name suggests, any relaxed monoidal category has an underlying
relaxed multicategory: just as we got a multicategory $C$ from a monoidal
category $D$ by putting
\[
\multihom{C}{\range{a_1}{a_n}}{a} = \homset{D}{a_1 \otimes\cdots\otimes
a_n}{a}, 
\]
we will similarly get a relaxed multicategory $C$ from a relaxed monoidal
category $D$ by putting
\[
\relhom{C}{\tau}{a_1}{a_n}{a} = 
\homset{D}{\atsr_{\tau}\bftuple{a_1}{a_n}}{a}.
\]

A relaxed monoidal category is a category $D$ together with $n$-fold
tensor functors $\atsr_n: D^n \go D$ for each $n$, but these are \emph{not}
required to patch together isomorphically: instead, we just have maps like
\begin{eqnarray}	\label{eq:rel-mon-cat-maps}
\obt \obt a_1 \atsr a_2\cbt \atsr a_3\cbt &\go &\obt a_1 \atsr a_2 \atsr
a_3\cbt ,\nonumber\\ 
\obt a_1 \atsr \obt a_2 \atsr a_3 \atsr \obt \, \cbt \cbt
\cbt &\go&\obt \obt a_1\cbt \atsr \obt a_2 \atsr a_3\cbt \cbt ,
\end{eqnarray}
where
\begin{eqnarray*}	
\atsr_n\bftuple{a_1}{a_n}	&\mbox{ is written as }	
&[a_1 \atsr \cdots \atsr {a_n}] \mbox{\ ($n\geq 2$)},	\\
\atsr_1(a) 			&\mbox{ is written as }	&[a],	\\
\atsr_0				&\mbox{ is written as }	&[\,].
\end{eqnarray*}
The maps in~\bref{eq:rel-mon-cat-maps} come from maps of trees whose leaves
are notionally labelled by the $a_i$'s: e.g.\ in the first case,
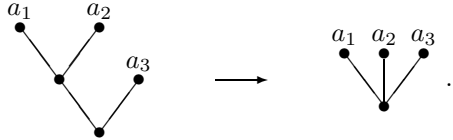
\[	
\begin{tree}
\nl{a_1}&	&	&	&\nl{a_2}&	&	\\
	&\rt{2}	&	&\lt{2}	&	&	&	\\
	&	&\node	&	&	&	&\nl{a_3}\\
	&	&	&\rt{2}	&	&\lt{2}	&	\\
	&	&	&	&\node	&	&	\\
\end{tree}
\diagspace \go \diagspace
\begin{tree}
\nl{a_1}&	&\nl{a_2}&	&\nl{a_3}\\
	&\rt{2}	&\dn	&\lt{2}	&	\\
	&	&\node	&	&	\\
\end{tree}
\ .
\]
However, relaxed monoidal categories can be defined without mention of trees:
\begin{defn} 
Let \triple{\stbk}{\eta}{\mu} be the free strict monoidal category monad on
\Cat. A \emph{relaxed monoidal category} consists of a category $D$ together
with a functor $\atsr: D^* \go D$ and natural transformations
\[
\begin{ntdiag}
D^{**}		&	&\rTo^{\mu_D}	&	&D^*		\\
		&	&		&\ 	&		\\
\dTo<{\atsr^*}	&	&\ruNT>{\gamma}	&	&\dTo>{\atsr}	\\
		&\ 	&		&	&		\\
D^*		&	&\rTo_{\atsr}	&	&D		\\
\end{ntdiag}
\diagspace
\begin{ntdiag}
D	&		&\rTo^{\eta_D}	&		&D^*		\\
	&\rdTo(4,4)<{1}	&		&\ruNT>{\iota}	&		\\
	&		&		&		&\dTo>{\atsr}	\\
	& 		&		&		&		\\
	&		&		&		&D		\\
\end{ntdiag}
\]
satisfying axioms looking like associativity and identity: that is,
\[
\begin{diagram}
\atsr\of\atsr^{*}\of\atsr^{**}	&\rTo^{\gamma*1}&\atsr\of\mu_{D}\of\atsr^{**}
&\rEquals	&\atsr\of\atsr^{*}\of\mu_{D^*}	\\
\dTo<{1*\gamma^*}	&	&	&	&\dTo>{\gamma*1}\\
\atsr\of\atsr^{*}\of\mu^*_{D}	&\rTo_{\gamma *1}
&\atsr\of\mu_{D}\of\mu^*_{D}	&\rEquals
&\atsr\of\mu_{D}\of\mu_{D^*}\\
\end{diagram}
\]
\[
\begin{diagram}
1_{D}\of\atsr	&\rTo^{\iota*1}	&\atsr\of\eta_{D}\of\atsr
&\rEquals	&\atsr\of\atsr^{*}\of\eta_{D^*}	\\
		&\rdTo(4,2)<{1}	&
&		&\dTo>{\gamma*1}		\\
		&		&
&		&\atsr\of\mu_{D}\of\eta_{D^*}	\\
\end{diagram}
\diagspace
\begin{diagram}
\atsr\of\atsr^{*}\of\eta_{D}^*	&\lTo^{1*\iota^*}	&\atsr\of 1_{D}^*\\
\dTo<{\gamma*1}			&\ldTo>{1}		&		\\
\atsr\of\mu_{D}\of\eta^*_{D}	&			&		\\
\end{diagram}
\]
all commute.
\end{defn}

So $\gamma$ assigns maps
\[
[[a_1^1 \atsr\cdots\atsr a_1^{k_1}] \atsr\cdots\atsr 
[a_n^1 \atsr\cdots\atsr a_n^{k_n}]]
\go
[a_1^1 \atsr\cdots\atsr a_n^{k_n}]
\]
($n,k_i \in\nat$), and $\iota$ assigns maps $a\go [a]$.

\begin{eg}{egs:rel-mon-cats}
\item	\label{eg:triv-rel-mon}
A relaxed monoidal category in which $\gamma$ and $\iota$ are isomorphisms is
just a monoidal category (of the unbiased kind---see \ref{sec:struc}).

\item	\label{eg:recipe}
Here is a general recipe for constructing relaxed monoidal categories. Let
\pr{D}{\otimes} be a genuine monoidal category, and \Mnd\ a monad on $D$
in which $T$ has the structure of a lax monoidal functor. Write
$\otimes\bftuple{a_1}{a_n}$ as $\langle a_1 \otimes\cdots\otimes 
a_n\rangle$. Now define
\[
[a_1 \atsr\cdots\atsr a_n] = T\langle a_1 \otimes\cdots\otimes a_n \rangle,
\]
$\iota_a$ as
\[
a \goby{\eta_a} Ta \goiso T\langle a\rangle = [a],
\]
and $\gamma$ by
\begin{eqnarray*}
\lefteqn{[[a_1^1 \atsr\cdots\atsr a_1^{k_1}] \atsr\cdots\atsr 
[a_n^1 \atsr\cdots\atsr a_n^{k_n}]]}\\
&=& T\langle T\langle a_1^1 \otimes\cdots\otimes a_1^{k_1} \rangle
\otimes\cdots\otimes T\langle a_n^1 \otimes\cdots\otimes a_n^{k_n} \rangle\rangle \\
&\go& T^2\langle \langle a_1^1 \otimes\cdots\otimes a_1^{k_1} \rangle
\otimes\cdots\otimes \langle a_n^1 \otimes\cdots\otimes a_n^{k_n} \rangle\rangle \\
&\goiso& T^2\langle a_1^1 \otimes\cdots\otimes a_n^{k_n} \rangle	\\
&\go& T\langle a_1^1 \otimes\cdots\otimes a_n^{k_n} \rangle	\\
&=& [a_1^1 \atsr\cdots\atsr a_n^{k_n}],
\end{eqnarray*}
where the first map comes from $T$ being lax monoidal, the second comes
from a coherence isomorphism in $D$, and the third is a component of
$\mu$.

\item
Let $D$ be a category with finite coproducts: then \pr{D}{+} is a monoidal
category, and any endofunctor on $D$ naturally has the structure of a lax
monoidal functor. So by~\bref{eg:recipe}, any monad \Mnd\ on $D$ gives a
relaxed monoidal category with $[a_1 \atsr\cdots\atsr a_n] = T(a_1 +\cdots +
a_n)$.

\item	\label{eg:rel-mon-monoid}
In~\bref{eg:recipe}, take the monoidal category \pr{\Set}{+}, a monoid $M$,
and the monad $T=M\times\dashbk$. Then $T$ actually preserves $+$ up to
coherent isomorphism (i.e.\ is a strong monoidal functor), so
\[
[X_1 \atsr\cdots\atsr X_n] = M\times (X_1 +\cdots + X_n)
\iso (M\times X_1) +\cdots + (M\times X_n).
\]

\item 	\label{eg:gen-mon}
Suppose \pr{D}{\otimes} is a \emph{symmetric} monoidal category, $M$ is a
monoid in \pr{D}{\otimes}, and $T=M\otimes\dashbk$. The monoid $M$
corresponds to a lax monoidal functor $1\goby{\bar{M}}D$, and $T$ is then the
composite
\[
D \goiso 1\times D \goby{\bar{M}\times 1} D\times D \goby{\otimes} D.
\]
But \pr{D}{\otimes} being symmetric means that $\otimes$ is a monoidal
functor, so $T$ has the structure of a lax monoidal functor. So we get a
relaxed monoidal category \pr{D}{\atsr} with
\[
[a_1 \atsr\cdots\atsr a_n] = M\otimes a_1 \otimes\cdots\otimes a_n.
\]

\item	\label{eg:rel-mon-exceptions}
As a concrete instance of the last example, take $\pr{D}{\otimes} =
\pr{\Set}{+}$ and $M=1$ to get a relaxed multicategory whose objects are sets
and with
\[
[X_1 \atsr\cdots\atsr X_n] = X_1 +\cdots + X_n + 1.
\]

\end{eg}

We now give an alternative definition of relaxed monoidal category, more
suitable for generalizing to other dimensions. If \pr{D}{\atsr} is a relaxed
monoidal category then for each $\tau\in\TR{n}$ there is a functor
$\atsr_{\tau}: D^n\go D$, defined inductively by
\begin{itemize}
\item $\atsr_{\utree} = 1_D$
\item $\atsr_{\abftuple{\tau_1}{\tau_n}} = \atsr_n \of (\atsr_{\tau_1}
\times\cdots\times \atsr_{\tau_n})$.
\end{itemize}
(We are using the inductive definition of trees given in~\cite[p.\ 13,
70]{SHDCT}.) Moreover, if $\tau'\go\tau$ is a map of trees then we get a
natural transformation $\atsr_{\tau'} \go \atsr_{\tau}$, by similar inductive
means. So we have found a functor
\[
\atsr_{\_}: \TR{n} \go \ftrcat{D^n}{D}
\]
for each $n$. These fit together nicely: if $\END(D)$ is the obvious
\pr{\Cat}{\ust}-operad with 
\[
(\END(D))(n) = \ftrcat{D^n}{D},
\]
then we have a map 
\[
\atsr_{\_}: \fcat{TR} \go \END(D)
\]
of \pr{\Cat}{\ust}-operads ($=$ \Cat-enriched operads:
see~\ref{egs:epms}\bref{eg:str-mon-two}). Conversely, the functors $\atsr_n$
can be recovered from $\atsr_{\_}$ by taking $\tau$ to be the $n$-leafed tree 
$
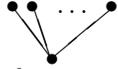
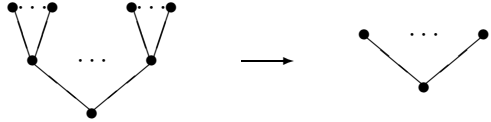
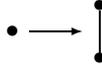
\begin{tree}
\node	&\node	&	&\ldots	&	&\node	\\
	&\rt{2}	&\rt{1}	&	&\lt{3}	&	\\
	&	&\node	&	&	&	\\
\end{tree},
$
and the transformations $\gamma$ and $\iota$ can be recovered by considering
the tree maps
\[
\begin{tree}
\node &\cdots &\node & & & &\node &\cdots &\node \\
 &\rt{1}\lt{1}& & & & & &\rt{1}\lt{1} & \\
 &\node & & &\cdots & & &\node & \\
 & &\rt{3} & & & &\lt{3} & & \\
 & & & &\node & & & & \\
\end{tree}
\diagspace\go\diagspace
\begin{tree}
\node & & &\cdots & & &\node \\
 &\rt{3} & & & &\lt{3} & \\
 & & &\node & & & \\
\end{tree}
\]
and
\[
\node \go
\begin{tree}
\node \\ \dn \\ \node \\
\end{tree}
\]
respectively. Hence:
\begin{propn}
A relaxed monoidal category is a category $D$ together with a map $\fcat{TR}
\go \END(D)$ of \pr{\Cat}{\ust}-operads. \done
\end{propn}

The data we have just assembled is enough to show that every relaxed monoidal
category yields a relaxed multicategory. For if $\atsr_{\_}: \fcat{TR} \go
\END(D)$ is a relaxed monoidal category, then there is a relaxed
multicategory $C$ defined by
\begin{itemize}
\item $C_0 = D_0$
\item $\relhom{C}{\tau}{a_1}{a_n}{a} =
\homset{D}{\atsr_{\tau}\bftuple{a_1}{a_n}}{a}$ 
\item if $\tau'\go\tau$ is a map of trees then 
\[
\relhom{C}{\tau}{a_1}{a_n}{a} \go \relhom{C}{\tau'}{a_1}{a_n}{a}
\]
is induced by $\atsr_{\tau'}\go\atsr_{\tau}$
\item if $f: \atsr_{\tau}\bftuple{a_1}{a_n} \go a$ and $f_i:
\atsr_{\tau_i}\bftuple{a_i^1}{a_i^{k_i}} \go a_i$ are maps in $D$ (with
$\tau\in\TR{n}, \tau_i\in\TR{k_i}, 1\leq i\leq n$), then their composite in
$C$ is the composite
\begin{eqnarray*}
\lefteqn{\atsr_{\tau\of\bftuple{\tau_1}{\tau_n}}\bftuple{a_1^1}{a_n^{k_n}}}\\
&\goby{1}&\atsr_{\tau} \bftuple{\atsr_{\tau_1}
\bftuple{a_1^1}{a_1^{k_1}}}{\atsr_{\tau_n} \bftuple{a_n^1}{a_n^{k_n}}} \\
&\goby{\atsr_{\tau}\bftuple{f_1}{f_n}}&
\atsr_{\tau}\bftuple{a_1}{a_n}	\\
&\goby{f}&
a	
\end{eqnarray*}
in $D$.
\item if $a\in C_0$ then the identity on $a$ is $\atsr_{\utree}(a) \goby{1}
a$. 
\end{itemize}
We call $C$ the \emph{underlying} relaxed multicategory of $D$.

In fact, all but one of our examples~\ref{egs:rel-mtis} of relaxed
multicategories come from relaxed monoidal categories (and the exception
is handled in section~\ref{sec:fun}):
\begin{eg}{egs:rel-mtis-from-mons}
\item 
If we take the relaxed multicategory $C$ arising trivially from a monoidal
category $E$ (see~\ref{egs:rel-mtis}\bref{eg:triv-rel-mti}), and the relaxed
monoidal category $D$ arising trivially from $E$
(see~\ref{egs:rel-mon-cats}\bref{eg:triv-rel-mon}), then $C$ is the relaxed
multicategory underlying $D$.

\item	\label{eg:new-monoid}
The relaxed multicategory of~\ref{egs:rel-mtis}\bref{eg:monoid} underlies
the relaxed monoidal category $D$
of~\ref{egs:rel-mon-cats}\bref{eg:rel-mon-monoid}. For by a simple induction
on $\tau$,
\[
\atsr_{\tau}\bftuple{X_1}{X_n} =
T^{h_1(\tau)}X_1 + \cdots + T^{h_n(\tau)}X_n
\]
(just on the grounds that $T=(M\times\dashbk)$ preserves the monoidal structure
$+$), giving the formula of~\ref{egs:rel-mtis}\bref{eg:monoid}.

\item	\label{eg:together-exceptions}
Similarly, \ref{egs:rel-mtis}\bref{eg:exceptions} comes
from~\ref{egs:rel-mon-cats}\bref{eg:rel-mon-exceptions}. In the
relaxed monoidal category we have 
\[
\atsr_{\tau}\bftuple{X_1}{X_n} =
X_1 +\cdots + X_n + v(\tau),
\]
so in the underlying relaxed multicategory,
\begin{eqnarray*}
\relhom{\Hom}{\tau}{X_1}{X_n}{X} 
&= &\homset{\Set}{X_1 +\cdots + X_n +v(\tau)}{X} \\
&= &X^{v(\tau)} \times \prod_{1\leq i \leq n} \homset{\Set}{X_i}{X}.
\end{eqnarray*}

\item	\label{eg:new-R-algebra}
With just the same reasoning as in~\bref{eg:together-exceptions},
Example~\ref{egs:rel-mtis}\bref{eg:R-algebra} comes
from \ref{egs:rel-mon-cats}\bref{eg:gen-mon}.

\end{eg}

We have seen how a relaxed multicategory arises from a category $C$ together
with a map $\fcat{TR}\go\END(C)$ of \pr{\Cat}{\ust}-operads. We now give an
analogue one level down, for relaxed categories; this seems as if it should
generalize quite easily to relaxed $T_n$-multicategories for all $n$, but we
do not attempt this here. 

So: \Cpn{1} is \pr{\Set}{\mr{free\ monoid}}, and by taking internal
categories we obtain the category \Cat\ and the free strict monoidal category
monad \ust\ on it. One level down, \Cpn{0} is \Zeropr, and taking internal
categories gives the category \Cat\ with the identity monad on it. The
analogue of $\fcat{TR}=\PD{2}$ is $\Delta=\PD{1}$, and so the analogue of a
relaxed monoidal category is a category $D$ together with a map
$\Delta\go\END(D)$ of \pr{\Cat}{\id}-operads. Now \pr{\Cat}{\id}-operads are
just strict monoidal categories, and $\END(D)$ is being used to denote the
familiar monoidal category \ftrcat{D}{D}, so what we have is just a category
$D$ with a monad on it. Our argument by analogy therefore suggests: given a
category $D$ with a monad $T$ on it, we should obtain a relaxed category $C$
by putting $C_0=D_0$ and
\[
\ehom{C_n}{a}{b} = \homset{D}{T^n a}{b}.
\]
And indeed, this is precisely
Example~\ref{egs:rel-cats}\bref{eg:monad-rel-cat}.

\section{Modifying the Codomain}	\label{sec:fun}

In the previous section we defined relaxed multicategories and relaxed
categories $C$ by, respectively, 
\begin{eqnarray*}
\relhom{C}{\tau}{a_1}{a_n}{a}
&=&
\homset{D}{\atsr_{\tau}\bftuple{a_1}{a_n}}{a},	\\
\ehom{C_n}{a'}{a}
&=&
\homset{D}{T^n a'}{a},
\end{eqnarray*}
where $D$ is a category and $\atsr_{\tau}$ and $T^n$ are functions with
suitable properties. Thus we obtain relaxed structures by modifying the
domain. This section is the dual: relaxed structures will be obtained by
modifying the codomain. More specifically, given a plain multicategory $D$ we
define a relaxed multicategory $C$ by $C_0=D_0$ and
\begin{equation}	\label{eq:rel-mti}
\relhom{C}{\tau}{a_1}{a_n}{a} =
\multihom{D}{\range{a_1}{a_n}}{Q_{\tau}a},
\end{equation}
where $Q_{\_}$ is some suitable family of functions. Similarly, given a
category $D$ we define a relaxed category $C$ by $C_0 = D_0$ and
\begin{equation}	\label{eq:rel-cat}
\ehom{C_n}{a'}{a} = \homset{D}{a'}{Q_n a}.
\end{equation}
The question is, then: what must $Q_{\_}$ be?

Let us take relaxed categories first. In order for \ehom{C_{\_}}{a'}{a} to be
a simplicial set for fixed $a'$ and $a$, $Q_{\_}a$ must be a functor
$\Delta^{\op} \go D$. The only sensible way of defining composites will be
if we have a map $Q_{n'}Q_{n}a \go Q_{n'+n}a$ for each $n'$, $n$ and $a$: for
then, if $f\in\homset{D}{a'}{Q_n a}$ and $f'\in\homset{D}{a''}{Q_{n'}a'}$,
the composite $f'\of f$ in $C$ can be defined as the composite
\[
a'' \goby{f'} Q_{n'}a' \goby{Q_{n'}f} Q_{n'}Q_{n}a \go Q_{n'+n}a
\]
in $D$. Here we have used the expression $Q_{n'}f$, so $Q_{n'}$ should
be a functor $D\go D$; thus $Q_{\_}$ is a functor
$\Delta^{\op}\go\ftrcat{D}{D}$. To get identities, we will also need a map
$a\go Q_0 a$ for each $a$. 

We have now argued that a functor $Q: \Delta^{\op} \go \ftrcat{D}{D}$,
together with maps $Q_{n'}Q_{n}a \go Q_{n'+n}a$ and $a\go Q_0 a$ satisfying
some axioms, will give us a relaxed category via~\bref{eq:rel-cat}. These
axioms assert precisely that $Q$ is a lax monoidal functor, which is the same
as a map of plain multicategories (see~\ref{sec:struc}). Hence:
\begin{thm}	\label{thm:fun-zero}
Let $D$ be a category and $Q_{\_}: \Delta^{\op} \go \ftrcat{D}{D}$ a map of
plain multicategories. Then there is a relaxed category $C$ with $C_0 = D_0$
and $\ehom{C_n}{a'}{a} = \homset{D}{a'}{Q_n a}$.\done
\end{thm}

\begin{eg}{egs:fun-zero}
\item
If $Q_{\_}$ is constant with value $1_D$ then the resulting relaxed category
is the $C$ of \ref{egs:rel-cats}\bref{eg:cat-to-rel}.

\item 
Comonads on a category $D$ correspond to \emph{strict} monoidal functors
$\Delta^{\op}\linebreak \go \ftrcat{D}{D}$, so for any comonad $G$ on $D$ there
is a relaxed category $C$ with $\ehom{C_n}{a'}{a} = \homset{D}{a'}{G^n
a}$. This is Example~\ref{egs:rel-cats}\bref{eg:comonad-rel}.

\item
Let $R$ be a commutative ring and $L$ a coalgebra over $R$, i.e.\ a comonoid
in the monoidal category \pr{R\hyph\fcat{Mod}}{\otimes_R} of left
$R$-modules. Then $R$-modules can be made to form a relaxed category as
follows. 

Firstly, there is a lax monoidal functor
\[
\overline{\otimes_{\integers}}: 
R\hyph\fcat{Mod} \go \ftrcat{R\hyph\fcat{Mod}}{R\hyph\fcat{Mod}}
\]
which is the transpose of the functor
\[
\otimes_{\integers}:
R\hyph\fcat{Mod} \times R\hyph\fcat{Mod} \go R\hyph\fcat{Mod}.
\]
(Note that we are using $\otimes_{\integers}$ rather than $\otimes_R$, so
$\overline{\otimes_{\integers}}$ is genuinely only a \emph{lax} monoidal
functor.)  Secondly, the coalgebra $L$ `is' the monoidal functor
$L^{\otimes_{R}(-)}: \Delta^{\op} \go R\hyph\fcat{Mod}$, sending $n$
to $L^{\otimes_{R} n}$. Hence
there is a lax monoidal functor
\[
\overline{\otimes_{\integers}} \of L^{\otimes_{R}(-)}: 
\Delta^{\op} \go \ftrcat{R\hyph\fcat{Mod}}{R\hyph\fcat{Mod}},
\]
giving a relaxed category $C$ whose objects are the $R$-modules. Explicitly,
a map $X'\go X$ in $C$ of degree $n$ is a homomorphism $X' \go
X\otimes_{\integers}L^{\otimes_R n}$ of left $R$-modules. 

\end{eg}

Moving up a level, we will demonstrate:
\begin{thm}	\label{thm:fun-one}
Let $D$ be a plain multicategory and $Q_{\_}: \fcat{TR}^{\op} \go
\ftrcat{D}{D}$ a map of $T_2$-multicategories. Then there is a relaxed
multicategory $C$ with $C_0 = D_0$ and 
\[
\relhom{C}{\tau}{a_1}{a_n}{a} = \multihom{C}{\range{a_1}{a_n}}{Q_{\tau}a}.
\]
\end{thm}
For this to make sense we have to define a $T_2$-multicategory \ftrcat{D}{D}
for any plain multicategory $D$ (although we have not done quite enough
formal work to do this with absolute rigour).

Firstly, recall the \fm-multicategory
of~\ref{egs:degen-fm}\bref{eg:bim-fm}, made up of rings,
$(R;\range{R_1}{R_n})$-modules, and maps between them. There is similarly an
\fm-multicategory $V$ in which a 1-cell is a category (in place
of a ring) and a horizontal 2-cell
\begin{opetope}
		&	&	&\cnr 	&\ldots	&	&	&	\\
		&\cnr	&\ruLine(2,1)^{A_2}&&	&	&\cnr 	&	\\
\ruLine(1,2)<{A_1}&	&	&	&\Downarrow &&	
&\rdLine(1,2)>{A_n}\\
\cnr		&	&	&	&\rLine_{A}&	&	&\cnr\\
\end{opetope}
is a profunctor $A_1 \times\cdots\times A_n \rMod A$ (in place of a
module). Suppose we fix a category $A$ and consider the substructure of
$V$ in which the only 1-cell allowed is $A$ and the only vertical
2-cell allowed is the identity on $A$. This is a $T_2$-multicategory
(see~\ref{egs:degen-fm}\bref{eg:V-T2}) which we call $\Prof(A)$. An object of
$\Prof(A)$ over $n$ is thus a profunctor $A^n \rMod A$, and arrows are, in a
sense analogous to the modules example, morphisms of profunctors.

Now we can define the $T_2$-multicategory \ftrcat{D}{D}, for any plain
multicategory $D$. Write $|D|$ for the underlying category of $D$. An object
of \ftrcat{D}{D} over $n$ is just a functor $|D|\go|D|$, for any $n$. For
each $n$ and $F: |D| \go |D|$ there is a profunctor $\hat{F}: |D|^n \rMod
|D|$, given as the composite
\[
\begin{array}{rcccl}
(|D|^n)^{\op}\times |D| &\goby{1\times F}&(|D|^n)^{\op}\times |D|
&\goby{\Hom_D}	&\Set	\\
(\range{a_1}{a_n};a)	&\goesto	&(\range{a_1}{a_n};Fa)
&\goesto&\multihom{D}{\range{a_1}{a_n}}{Fa}.
\end{array}
\]
An arrow
\begin{equation}	\label{eq:arr-in-end}
\piccy{endarrow}
\end{equation}
in \ftrcat{D}{D} is then an arrow
\[
\piccy{endarrowhat}
\]
in $\Prof(|D|)$, and similarly composition and identities: in other words,
the $T_2$-multicategory structure of \ftrcat{D}{D} is transported from that
of $\Prof(|D|)$ via the functions
\begin{eqnarray*}
\{\mr{functors\ } |D| \go |D| \} &\go 
&\{\mr{objects\ of\ } \Prof(|D|) \mr{\ over\ } n \} \\
F&\goesto&\hat{F}.
\end{eqnarray*}

\begin{figure}
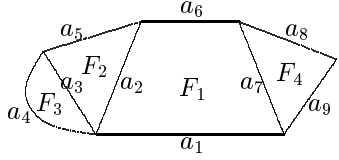

\centerline{\piccy{mental}}
\caption{Mental picture for placing of variables $a_i$}
\label{fig:var-placing}
\end{figure}

In elementary terms, a map $\phi$ as in~\bref{eq:arr-in-end} consists of a
rule
\[
\begin{array}{c}
a_2,a_6,a_7 \go F_1 a_1	
\ \ \ 
a_3,a_5 \go F_2 a_2
\ \ \ 
a_4 \go F_3 a_3
\ \ \ 
a_8,a_9 \go F_4 a_7	
\\
\hrulefill\\
a_4,a_5,a_6,a_8,a_9 \go Fa_1
\end{array}
\]
where \range{a_1}{a_9} are objects of $D$ (Figure~\ref{fig:var-placing}).
So given some arrows in $D$ as above the line, we get an arrow as below the
line; this assignment is to be compatible with composition with arrows in
$|D|$, in a sense suggested by the compatibility axioms
of~\ref{egs:degen-fm}\bref{eg:bim-fm}. 

We can now sketch a proof of Theorem~\ref{thm:fun-one}. Take a plain
multicategory $D$ and a map $Q_{\_}: \fcat{TR}^{\op} \go \ftrcat{D}{D}$. Put
$C_0 = D_0$. Fix $\range{a_1}{a_n},a \in C_0$. For any $\tau\in\TR{n}$ we
have a functor $Q_{\tau}: |D| \go |D|$, and we define
\[
\relhom{C}{\tau}{a_1}{a_n}{a} = \multihom{D}{\range{a_1}{a_n}}{Q_{\tau}a}.
\]
A map $\tau'\go\tau$ of trees induces a morphism $\widehat{Q_{\tau}} \go
\widehat{Q_{\tau'}}$ of profunctors, or equivalently a natural transformation
$Q_{\tau} \go Q_{\tau'}$. Thus we obtain a function
\[
\relhom{C}{\tau}{a_1}{a_n}{a} \go
\relhom{C}{\tau'}{a_1}{a_n}{a},
\]
as required.

This has covered all the structure of $C$ for fixed objects
$\range{a_1}{a_n},a$; next we turn to composition. Suppose we have trees
$\tau\in\TR{n}$, $\tau_1\in\TR{k_1}$, \ldots, $\tau_n\in\TR{k_n}$. The
canonical map
\[
\piccy{canonical}
\]
in the $T_2$-multicategory $\fcat{TR}^{\op}$ is sent by $Q_{\_}$ to a map
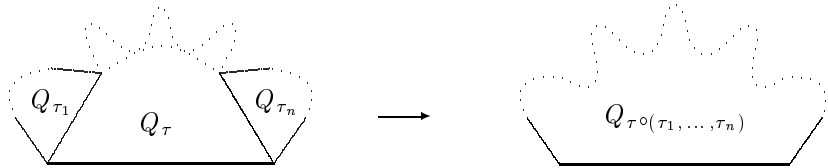
\[
\piccy{Qarrow}
\]
in \ftrcat{D}{D}. This map is a family of functions
\begin{eqnarray*}
\multihom{D}{\range{a_1}{a_n}}{Q_{\tau}a} \times
\multihom{D}{\range{a_1^1}{a_1^{k_1}}}{Q_{\tau_1}a_1} \times \cdots
\times 
\multihom{D}{\range{a_n^1}{a_n^{k_n}}}{Q_{\tau_n}a_n}\\
\go
\multihom{D}{\range{a_1^1}{a_n^{k_n}}}{Q_{\tau\of\bftuple{\tau_1}{\tau_n}}a},
\end{eqnarray*}
one for each $a,\range{a_1}{a_n},\range{a_1^1}{a_n^{k_n}}$, and so provides
composition in $C$. Identities go through similarly. We then just have to
check the axioms, and this is straightforward.

Theorem~\ref{thm:fun-one} was developed with representations of vertex groups
in mind (Example \ref{egs:rel-mtis}\bref{eg:rep}). Since the Hopf algebra $H$
is in particular a (not necessarily commutative) ring, there is a plain
multicategory $H\hyph\fcat{Mod}$ in which the objects are $H$-modules and
the arrows are multilinear maps. From $G=\pr{H}{K}$, a new $H$-module
$\mi{Fun}_{\tau}\pr{G^n}{A}$ is constructed for each $H$-module $A$ and tree
$\tau\in\TR{n}$. There is then a relaxed multicategory whose objects are
$H$-modules, and with
\[
\relhom{\Hom}{\tau}{A_1}{A_n}{A} = 
\multihom{H\hyph\fcat{Mod}}{\range{A_1}{A_n}}{\mi{Fun}_{\tau}\pr{G^n}{A}}.
\]
(This relaxed multicategory is not quite $\mb{Rep}(G)$, as there are are some
further subtleties involving $G$-invariance; it contains $\mb{Rep}(G)$ as
a substructure.) Whether or not $\mi{Fun}$ really does determine a
multicategory map 
\[
\fcat{TR}^{\op} \go 
\ftrcat{H\hyph\fcat{Mod}}{H\hyph\fcat{Mod}}
\] 
has probably yet to be verified.

\appendix

\chapter{Free Multicategories}

In this appendix we define `suitability' and sketch proofs of
Theorems~\ref{thm:free-main} and~\ref{thm:free-fixed}. First, we need some
terminology.

Let \Eee\ be a category with pullbacks, \scat{I} a small category, $D:
\scat{I}\go\Eee$ a functor for which a colimit exists, and $(D(I)\go
Z)_{I\in \scat{I}}$ a colimit cone. We say that the colimit is
\emph{stable under pullback} if for any map $Z'\go Z$ in \Eee, the cone
$(D'(I)\go Z')_{I\in \scat{I}}$ is a colimit cone; here $D'$ and the new cone
are obtained by pullback, so that
\begin{minidiagram}
D'	&\rTo	&D	\\
\dTo	&	&\dTo	\\
Z'	&\rTo	&Z	\\
\end{minidiagram}
is a pullback square in the functor category \ftrcat{\scat{I}}{\Eee}.

The morphisms $k_I$ in a colimit cone $(D(I)\goby{k_I} Z)_{I\in \scat{I}}$
will be called the \emph{coprojections} of the colimit, and in particular we
say that the colimit of $D$ `has monic coprojections' to mean that each $k_I$
is monic.

A category will be said to have \emph{disjoint finite coproducts} if it has
finite coproducts, these coproducts have monic coprojections, and for any
pair $A, B$ of objects, the square
\begin{minidiagram}
0	&\rTo	&B	\\
\dTo	&	&\dTo	\\
A	&\rTo	&A+B	\\
\end{minidiagram}
is a pullback.

Let $\omega$ be the natural numbers with their usual ordering. A \emph{nested
sequence} in a category \Eee\ is a functor $\omega\go\Eee$ in which the image
of every morphism of $\omega$ is monic. In other words, it is a diagram
\[
A_0 \monic A_1 \monic \cdots
\]
in \Eee, where as usual \monic\ indicates a monic. Note that a functor which
preserves pullbacks also preserves monics, so it makes sense for such a
functor to `preserve colimits of nested sequences'. Similarly, it makes sense
to say that colimits of nested sequences commute with pullbacks (by an easy
calculation), where `commute' is used in the same sense as when we say that
filtered colimits commute with finite limits in \Set.

Recall (from \cite{GOM} or \cite[Chapter I]{SHDCT}) that a category is called
\emph{cartesian} if it has all finite limits; a monad \Mnd\ is called
\emph{cartesian} if $T$ preserves pullbacks and the squares expressing the
naturality of $\eta$ and of $\mu$ are all pullbacks; the phrase `\Cartpr\ is
cartesian' is used to mean that $T$ is a cartesian monad on a cartesian
category \Eee. (Here, as elsewhere, we abuse language by calling a monad by
the name of its functor part.)

Finally, we can give the definition of suitability. A category \Eee\ is
\emph{suitable} if it satisfies
\begin{description}
\item[C1] \Eee\ is cartesian
\item[C2] \Eee\ has disjoint finite coproducts which are stable under
pullback
\item[C3] \Eee\ has colimits of nested sequences; these commute with
pullbacks and have monic coprojections.
\end{description}
A monad \Mnd\ is \emph{suitable} if it satisfies
\begin{description}
\item[M1] \Mnd\ is cartesian
\item[M2] $T$ preserves colimits of nested sequences.
\end{description}
We say that \Cartpr\ is \emph{suitable} when \Mnd\ is a suitable monad on a
suitable category \Eee.

We now sketch a proof of the main theorem,~\ref{thm:free-main}, on the
formation of free multicategories, which for convenience is re-stated here.
\begin{trivlist} \item 
\textbf{Theorem \ref{thm:free-main}}\ \itshape
Let \Cartpr\ be suitable. Then the forgetful functor
\[
\Cartpr\hyph\Multicat \goby{U} \Eeep = \Cartpr\hyph\Gph
\]
has a left adjoint, the adjunction is monadic, and if $T'$ is the resulting
monad on \Eeep\ then \Cartprp\ is suitable.
\end{trivlist}

\noindent\textbf{Sketch proof}
We proceed in four steps:
\begin{enumerate}
\item 	\label{ftr}
construct a functor $F: \Eeep \go \Cartpr\hyph\Multicat$
\item 	\label{adjn}
construct an adjunction between $F$ and $U$
\item 	\label{primesuit}
check that \Cartprp\ is suitable
\item 	\label{monadic}
check that the adjunction is monadic.
\end{enumerate}

Each step goes roughly as follows:
\begin{enumerate}
\item \emph{Construct a functor $F: \Eeep \go \Cartpr\hyph\Multicat$}\\ 
Let $X$ be a $T$-graph. Define for each $n$ a graph
\spaan{A_n}{TX_0}{X_0}{d_n}{c_n}, by
\begin{itemize}
\item $A_0=X_0$, $d_0=\eta_{X_0}$ and $c_0=1$
\item $A_{n+1} = X_0 + X_1\of A_n$, where $X_1\of A_n$ is the 1-cell
composite in $\Cartpr\hyph\Span$, with the obvious choices of $d_{n+1}$ and
$c_{n+1}$.
\end{itemize}
Define for each $n$ a map $A_n \goby{i_n} A_{n+1}$, by
\begin{itemize}
\item $i_0: X_0 \go X_0 + X_1 \of X_0$ is first coprojection
\item $i_{n+1} = 1_{X_0} + (1_{X_1} * i_n)$.
\end{itemize}
Then the $i_n$'s are monic, and by taking $A$ to be the colimit of 
\[
A_0 \rMonic^{i_0} A_1 \rMonic^{i_1} \cdots
\]
we obtain a graph \spn{A}{TX_0}{X_0}. This graph naturally has the structure
of a multicategory: the identities map $X_0 \go A$ is just the colimit
coprojection $A_0 \monic A$, and composition comes from maps $A_m \of A_n \\
\go A_{m+n}$ which piece together to give a map $A\of A \go A$. The latter
construction needs many of the axioms for suitability.

We have now described what effect $F$ is to have on objects, and extension to
morphisms is straightforward.

(Incidentally, the colimit of the nested sequence of $A_n$'s appears, in
light disguise, as the recursive description of the free plain multicategory
monad \fm\ in~\ref{sec:fm-mti}: $A_n$ is the set of formal expressions which
can be obtained from the first clause and up to $n$ applications of the
second clause.)

\item \emph{Construct an adjunction between $F$ and $U$}\\ 
We do this by constructing unit and counit transformations and verifying the
triangle identities. Both transformations are the identity on the object of
objects, so we only need to define them on the object of arrows. For the unit
$\eta'$, if $X\in\Eeep$ then $\eta'_X: X_1 \go A$ is the composite
\[
X_1 \goiso X_1 \of X_0 \rMonic^{\mr{copr}_2} X_0 + X_1 \of X_0 = A_1
\rMonic A.
\]
For the counit $\epsln'$, let $C\in\Cartpr\hyph\Multicat$. Write $A$ and
$A_n$ for the objects used in the construction of the free multicategory on
$U(C)$, as if $X=U(C)$ in part~\bref{ftr}. Define for each $n$ a map
$\epsln'_{C,n}: A_n \go C_1$ by
\begin{itemize}
\item $\epsln'_{C,0} = (A_0 \goby{=} C_0 \goby{\ids} C_1)$
\item $\epsln'_{C,n+1} = (C_0 + C_1 \of A_n \goby{1+1*\epsln'_{C,n}} 
C_0 + C_1 \of C_1 
\goby{\left(\begin{array}{c}\scriptstyle \ids\\ 
\scriptstyle \comp\end{array}\right)} 
C_1)$,
\end{itemize}
and there is a unique $\epsln'_C: A \go C_1$ such that $\epsln'_{C,n} =
(A_n\monic A \goby{\epsln'_C} C_1)$ for all $n$. 

\item \emph{Check that \Cartprp\ is suitable}\\
This is quite routine.

\item \emph{Check that the adjunction is monadic}\\
We apply the Monadicity Theorem by checking that $U$ creates coequalizers for
$U$-absolute coequalizer pairs. This can be done quite separately from the
rest of the proof, and again is quite routine.
\done 
\end{enumerate}

Theorem~\ref{thm:free-fixed}, the fixed-object version of the theorem, is
proved in much the same way. Indeed, if one has already
proved~\ref{thm:free-main} then much of~\ref{thm:free-fixed} is an easy
consequence. The most substantial difference between the two cases is that
the inclusion functor $\Eeep_S \go \Eeep$ does not preserve coproducts
(although it does preserve pullbacks and colimits of nested sequences).

Finally, we need to see that \Zeropr\ is suitable. The only part of this
which is not quite standard is \textbf{C3}, and this is an easy
calculation. (The fact that pullbacks commute with colimits of nested
sequences is also a special case of the fact that finite limits commute with
filtered colimits in \Set.) More generally, axioms \textbf{C} hold for any
presheaf category (since they hold for \Set), and axioms \textbf{M} for any
finitary cartesian monad.

\end{document}